\documentclass[11pt]{article}
\usepackage{amsmath}
\usepackage{amssymb}
\usepackage{amsfonts}
\usepackage{color}
\usepackage[colorlinks=true]{hyperref}
\usepackage{epsfig}
\usepackage{srcltx}
\usepackage{multirow}
\usepackage{graphicx}

\def\crd{\color{red}}
\def\pen{{\hbox{\scriptsize\rm pen}}}
\def\ds{{\hbox{\scriptsize\rm DS}}}
\def\cP{{\cal P}}

\def\cV{{\cal V}}

\def\reg{{\hbox{\scriptsize\rm reg}}}
\def\las{{\hbox{\scriptsize\rm las}}}
\def\qed{\hfill $\Box$}

\def\las{\hbox{\scriptsize\rm lasso}}
\def\Argmin{\mathop{\hbox{\rm Argmin}}}
\def\bR{{\mathbf{R}}}
\def\bO{{\mathbf{O}}}
\def\bH{{\mathbf{H}}}

\def\Risk{{\,\hbox{\rm Risk}}}
\def\sign{{\hbox{\rm sign}}}
\def\Prob{{\hbox{\rm Prob}}}

\def\RIP{{\hbox{\rm RIP}}}
\def\RE{{\hbox{\rm RE}}}

\def\Opt{{\mathop{\hbox{\rm Opt}}}}
\def\sign{{\mathop{\hbox{\rm sign}}}}

\def\cR{{\cal R}}

\newcommand{\bbr}{{\Bbb{R}}}

\def\ErfInv{{\,\hbox{\rm erfinv}}}
\newcommand{\cN}{{\cal N}}

\newcommand{\cU}{{\cal U}}

\newcommand{\half}{ \mbox{\small$\frac{1}{2}$}}
\newcommand{\four}{ \mbox{\small$\frac{1}{4}$}}

\def\eps{\varepsilon}

\def\e{\epsilon}

\def\wh#1{\widehat{#1}}

\def\Prob{{\hbox{\rm Prob}}}

\newcommand{\epr}{\hfill\hbox{\hskip 4pt
                \vrule width 5pt height 6pt depth 1.5pt}\vspace{0.5cm}\par}

\newcommand{\be}{\begin{eqnarray}}
\newcommand{\ee}[1]{\label{#1}\end{eqnarray}}
\newcommand{\nn}{\nonumber \\}
\newcommand{\ese}{\end{eqnarray*}}
\newcommand{\bse}{\begin{eqnarray*}}
\newcommand{\rf}[1]{~(\ref{#1})}
\def\Conv{{\hbox{\rm Conv}}}

\newtheorem{lemma}{Lemma}
\newtheorem{observation}{Observation}
\newtheorem{proposition}{Proposition}
\newtheorem{corollary}{Corollary}
\newtheorem{algorithm}{Algorithm}

\newcommand{\an}[2]{{\color{green} #1}{ \color{magenta} #2}}
\begin{document}\title{Accuracy guarantees for $\ell_1$-recovery
\thanks{Research of the second author was supported by the Office of Naval Research grant \# N000140811104 and the NSF grant DMS-0914785.}}

\author{Anatoli Juditsky \\
              LJK, Universit\'e J. Fourier, B.P. 53,
              38041 Grenoble Cedex 9, France \\
              {\tt Anatoli.Juditsky@imag.fr}
           \and
           Arkadi Nemirovski\\
              Georgia Institute of Technology,
              Atlanta, Georgia 30332, USA \\
              {\tt nemirovs@isye.gatech.edu}
}

\maketitle
\begin{abstract}
We discuss two new methods of recovery of sparse signals from noisy observation based on $\ell_1$-minimization. While they are closely related to the well-known techniques such as Lasso and Dantzig  Selector, these estimators come with {\em efficiently verifiable guaranties of performance}. By optimizing these bounds with respect to the method parameters we are able to construct the estimators which possess better statistical properties than the commonly used ones.

We link our performance estimations to the well known results of Compressive Sensing and justify our proposed approach with an oracle inequality which links the properties of the recovery algorithms and the best estimation performance when the signal support is known. We also show how the estimates can be computed using the Non-Euclidean Basis Pursuit algorithm.
\end{abstract}
\section{Introduction}\label{intro}\label{Intro}
Recently several methods of estimation and selection which refer to the $\ell_1$-minimization received much attention in the statistical literature. For instance, {\em Lasso estimator}, which is the $\ell_1$-penalized least-squares method is probably the most studied (a theoretical analysis of the Lasso estimator is provided in, e.g., \cite{tsyb1,sara2,tsyb2, sara1,chzhang,tzhang,loun,mein}, see also the references cited therein). Another, closely related to the Lasso, statistical estimator is the {\em Dantzig Selector} \cite{Candes07,tsyb1,kolch09,loun}.  To be more precise, let us consider the estimation problem as follows.
Assume that an observation
\begin{equation}\label{nun}
y=Ax+\sigma \xi\in\bR^m
\end{equation}
is available,
where $x\in \bR^n$ is an unknown signal and $A\in \bR^{m\times n}$ is a known {\sl sensing matrix}.
We suppose that $\sigma \xi$ is a Gaussian disturbance with $\xi\sim N(0,I_m)$ (i.e., $\xi=(\xi_1,....,\xi_n)^T$,
where  $\xi_i$ are independent normal r.v.'s with zero mean and unit variance), and $\sigma>0$ is a known
deterministic noise level. Our focus is on the recovery of unknown signal $x$.

The Dantzig Selector estimator $\wh{x}_\ds$ of the signal $x$ is defined as follows \cite{Candes07}:
\[
\wh{x}_\ds(y)\in \Argmin_{v\in \bR^n}\{\|v\|_1\,|\;\|A^T(Av-y)\|_\infty\le \rho\}
\]
where $\rho=O\left(\sigma\sqrt{\ln n}\right)$ is the algorithm's parameter. Since $\wh{x}_\ds$ is obtained as a solution of an linear  program, it is very attractive by its low computational cost. Accuracy bounds for this estimator are readily available.
For instance, a well known result about this estimator (cf. \cite[Theorem 1.1]{Candes07}) is that if $\rho=O\left(\sigma\sqrt{\ln(n/\e)}\right)$ then
\[
\|\wh{x}_\ds(y)-x\|_2\le K\sigma\sqrt{s\log(n\e^{-1})}
\]
with probability $1-\e$ if $a)$ the signal $x$ is $s$-sparse, i.e. has at most $s$ non-vanishing components, and $b)$ the sensing matrix $A$ with unit columns possesses the {\em Restricted Isometry Property} $\RIP(\delta,k)$ with parameters $0<\delta<{1\over 1+\sqrt{2}}$ and $k\ge 3s$.
\footnote{Recall that $\RIP(\delta,k)$, called also {\em uniform uncertainty
principle}, means that for any $v\in \bR^n$ with at most $k$ nonzero entries,
\[
(1 - \delta) \|v\|_2\le \|Av\|_2^2\le (1 + \delta) \|v\|_2
 \]
 This property essentially
requires that every set of columns of $A$ with cardinality less than $k$ approximately behaves like
an orthonormal system.}
Further, in this case one has $K={C(1-\delta)^{-1}}$, where $C$ is a moderate absolute constant. This result is quite impressive,
in part due to the  fact (see, e.g. \cite{CRT1,CRT2}) that there exist $m\times n$ random matrices, with $m<n$,
which possess the $\RIP$ with  probability close to $1$, $\delta$ close to zero and the value of $k$ as large as $O\left({m \ln^{-1}(n/m)}\right)$.
Similar performance guarantees are known for {\sl Lasso recovery}
$$
\wh{x}_{\las}(y)\in\Argmin_{v\in\bR^n}\left\{\|v\|_1+\varkappa\|Av-y\|_2^2 \right\},
$$
with properly chosen penalty parameter $\varkappa$.
{The available accuracy bounds for Lasso and Dantzig Selector rely upon the Restricted Isometry Property or less restrictive assumptions about the sensing matrix, such as Restricted Eigenvalue \cite{tsyb1} or Compatibility \cite{sara2} conditions (a complete overview of those and several other assumptions
with description of how they relate to each other is provided in \cite{sara1}). However, these assumptions cannot be verified efficiently.}
The latter implies that there is currently no way to provide any guaranties  (e.g., confidence sets) of the performance of the proposed procedures. A notable exception from this rule is the {\em Mutual Incoherence} assumption (see, e.g. \cite{DE,EB,GN} and \cite{tzhang} for the case of, respectively, deterministic and random observation noise) which can be used to compute the accuracy bounds for recovery algorithms: a matrix $A$ with columns of unit $\ell_2$-norm and {\em mutual incoherence} $\mu(A)$ possesses $\RIP(\delta,k)$ with $\delta=(m-1)\mu(A)$.\footnote{
The {\em mutual incoherence} $\mu(A)$ of a sensing matrix $A=[A_1,...,A_n]$ is computed according to
\[
\mu(A)=\max_{i\neq j}{|A_i^TA_j|\over A_i^TA_i}.
\]
Obviously, the mutual incoherence can be easily computed even for large matrices. }
Unfortunately, the latter relation implies that $\mu(A)$ should be very small to certify the possibility
of accurate $\ell_1$-recovery of non-trivial sparse signals,
so that  performance guarantees based on mutual incoherence are very conservative. This ``theoretical observation'' is supported by numerical experiments -- the practical guarantees which may be obtained using the mutual incoherence are generally quite poor even for the problems with nice theoretical properties (cf. \cite{JNCS,JKKN}).
\par
 Recently the authors have proposed a new approach for efficient computing of upper and lower bounds on the
 ``level of goodness'' of a sensing matrix $A$, i.e. the maximal $s$ such that the $\ell_1$-recovery of {\em all} signals with no more than $s$
 non-vanishing components is accurate in the case where the measurement noise vanishes (see \cite{JNCS}).
 In the present paper we aim to use the related {\em verifiable sufficient conditions} of ``goodness'' of a sensing
 matrix $A$ to provide efficiently computable bounds for the error of $\ell_1$ recovery procedures in the case when the observations are affected
 by random noise.
\par
The main body of the paper is organized as follows:
\begin{enumerate}
\item We start with Section \ref{conditionhgamma} where we formulate the sparse recovery problem and introduce our core assumption -- a verifiable condition $\bH_{s,\infty}(\kappa)$ linking matrix $A\in\bR^{m\times n}$ and a {\sl contrast matrix} $H\in\bR^{m\times n}$. In Sections \ref{sectregrec}, \ref{sectpenrec} we present two recovery routines with contrast matrices:
\begin{itemize}
\item {\em regular recovery}: $$ \wh{x}_{\reg}(y)\in \Argmin\limits_{v\in \bR^n} \{\|v\|_1: \|H^T(Av-y) \|_\infty\le \rho \},$$
\item {\em penalized recovery}: $$\wh{x}_{\pen}(y)\in \Argmin\limits_{v\in \bR^n} \{\|v\|_1+\theta s\|H^T(Av-y)\|_\infty \},$$ ($s$ is our guess for the number of nonzero entries in the true signal,  $\theta>0$ is the penalty parameter)
\end{itemize}
along with their performance guarantees under condition $\bH_{s,\infty}(\kappa)$ with $\kappa<1/2$, that is, explicit upper bounds on the confidence levels of the recovery errors $\|\widehat{x}-x\|_p$. The novelty here is that our bounds are of the form
\begin{equation}\label{conditionbH}
\Prob\left\{\|\widehat{x}-x\|_p\leq O\left(s^{1/p}\sigma\sqrt{\ln(n/\epsilon)}\right)\;\hbox{for every $s$-sparse signal $x$ and all $1\leq p\leq\infty$}\right\}\geq1-\epsilon
\end{equation} (with hidden factors in $O(\cdot)$ independent of $\epsilon,\sigma$), and {\em are valid in the entire range $1\leq p\leq \infty$ of values of $p$}. Note that similar error bounds for Dantzig Selector and Lasso are only known for $1\leq p\leq 2$, whatever be the assumptions on ``essentially nonsquare'' matrix $A$.

\item Our interest in condition $\bH_{s,\infty}(\kappa)$  stems from the fact that this condition, in contrast to the majority of the known sufficient conditions for the validity of $\ell_1$-based sparse recovery (e.g., Restricted Isometry/Eigenvalue/Compatibility), is efficiently verifiable. Moreover, it turns out that one can efficiently optimize the error bounds of the associated with this verifiable condition regular/penalized recovery routines over the contrast matrix $H$. The related issues are considered in Section \ref{matrixh}. In Section \ref{boundingomega}
we provide some additional justification of the condition $\bH$, in particular, by linking it with the Mutual Incoherence and Restricted Isometry properties. This, in particular, implies that the condition $\bH_{s,\infty}(\kappa)$ with, say, $\kappa={1\over 3}$ associated with  randomly selected $m\times n$ matrices $A$ is feasible, with probability approaching 1 as $m,n$ grow, for $s$ as large as $O(\sqrt{m/\ln(n)})$. We also establish {\sl limits of performance}
of the condition, specifically, show that unless $A$ is nearly square, $\bH_{s,\infty}(\kappa)$ with $\kappa<1/2$ can be feasible only when $s\leq O(1)\sqrt{m}$, meaning that the tractability of the condition has a heavy price: when designing and validating $\ell_1$ minimization based sparse recovery routines, this condition can be useful only in a severely restricted range of the sparsity parameter $s$.
\item In Section \ref{secextent} we show that the condition $\bH_{s,\infty}(\kappa)$ is the strongest (and seemingly the only verifiable one) in a natural family of conditions $\bH_{s,q}(\kappa)$ linking a sensing and a contrast matrix; here $s$ is the number of nonzeros in the sparse signal to be recovered $q\in[1,\infty]$. We demonstrate that when a contrast matrix $H$ satisfies $\bH_{s,q}(\kappa)$ with $\kappa<1/2$, the associated regular and penalized $\ell_1$ recoveries admit error bounds similar to (\ref{conditionbH}), but now in the restricted range $1\leq p\leq q$ of values of $p$. We demonstrate also that feasibility of $\bH_{s,q}(\kappa)$ with $\kappa<1/2$ implies instructive (although slightly worse than those in (\ref{conditionbH})) error bounds
    for the Dantzig Selector  and Lasso recovering routines.
\item In Section \ref{sectnumres}, we present numerical results on comparison of regular/penalized $\ell_1$ recovery with the Dantzig Selector and Lasso algorithms. The conclusion suggested by these preliminary numerical results is that {\sl when the former procedures are applicable} (i.e., when the techniques of Section \ref{matrixh} allow to build a ``not too large'' contrast matrix satisfying the condition $\bH_{s,\infty}(\kappa)$ with, say, $\kappa=1/3$), {\sl our procedures outperform significantly the Dantzig Selector and work exactly as well as the Lasso algorithm with ``ideal'' (unrealistic in actual applications) choice of the regularization parameter\footnote{With ``theoretically optimal,'' rather than ``ideal,'' choice of the regularization parameter in Lasso,  this algorithm is essentially worse than our algorithms utilizing the contrast matrix.}.}
\item In the concluding Section \ref{sectmatpur} we present a ``Non-Euclidean Matching Pursuit algorithm'' (similar to the one presented in \cite{JKKN}) with the same performance characteristics as those of regular/penalized $\ell_1$ recoveries; this algorithm, however, does not require optimization and can be considered as a computationally cheap alternative to $\ell_1$ recoveries, especially in the case when one needs to process a series of recovery problems with common sensing matrix.
\end{enumerate}
All proofs are placed in the Appendix.
\section{Accuracy bounds for $\ell_1$-Recovery Routines}\label{ell1rout}\label{Ell1Recovery}

\subsection{Problem statement}\label{conditionhgamma}
\paragraph{Notation.}
For a vector $x\in \bR^n$ and $1\le s\le n$ we denote $x^s$ the vector obtained from $x$ by setting to $0$ all but the $s$ largest in magnitude entries of $x$. Ties, if any, could be resolved arbitrarily; for the sake of definiteness assume that among entries of equal magnitudes, those with smaller indexes have priority (e.g., with $x=[1;2;2;3]$ one has $x^2=[0;2;0;3]$). $\|x\|_{s,p}$ stands for the usual $\ell_p$-norm of $x^s$ (so that $\|x\|_{s,\infty}=\|x\|_\infty$).
We say that a vector $z$ is {\sl $s$-sparse} if it has at most  $s$ nonzero entries. Finally, for a set $I\subset\{1,...,n\}$ we denote by $J$ its complement $\{1,...,n\}\backslash I$;  given $x\in\bR^n$, we denote by $x_I$ the vector obtained from $x$ by zeroing the entries with indices outside of $I$, so that $x=x_I+x_{J}$.\par
Given a norm $\nu(\cdot)$ on $\bR^m$ and a matrix $H=[h_1,...,h_N]\in\bR^{m\times N}$, we set $\nu(H)=\max\limits_{i\leq N}\nu(h_i)$.

\paragraph{The problem.}
We consider an observation $y\in \bR^m$
\be
y=Ax+u+\sigma \xi,
\ee{obs}
where $x\in \bR^n$ is an unknown signal and $A\in \bR^{m\times n}$ is {\em the sensing matrix}. We suppose that $\sigma \xi$ is a Gaussian disturbance, where $\xi\sim N(0,I_m)$ (i.e., $\xi=(\xi_1,....,\xi_n)^T$ with independent normal random variables $\xi_i$ with zero mean and unit variance),  $\sigma>0$ being known, and $u$ is a nuisance parameter known to belong to a given {\sl uncertainty set} $\cU\subset\bR^m$ which we will suppose to be convex, compact and symmetric w.r.t. the origin. Our goal is to recover $x$ from $y$, provided that $x$ is ``nearly $s$-sparse.'' Specifically, we consider the sets
 \bse
 X(s,\upsilon)=\{x\in\bR^n:\|x-x^s\|_1\leq\upsilon\}
  \ese
 of signals which admit $s$-sparse approximation of $\|\cdot\|_1$-accuracy $\upsilon$. Given $p$, $1\leq p\leq\infty$, and a confidence level $1-\epsilon$, $\epsilon\in(0,1)$, we quantify a recovery routine --- a Borel function $\bR^m\ni y\mapsto\widehat{x}(y)\in\bR^n$ --- by its
 worst-case, over $x\in X(s,\upsilon)$, confidence interval, taken w.r.t. $\|\cdot\|_p$-norm of the error. Specifically, we define the {\sl risks} of a recovery routine as
\[
\Risk_p(\widehat{x}(\cdot)|\epsilon,\sigma,s,\upsilon) =\inf\left\{\delta:\Prob\{\xi:\exists x\in X(s,\upsilon),u\in \cU:
\|\widehat{x}(Ax+\sigma\xi+u)-x\|_p>\delta\}\leq\epsilon\right\}.
\]
Equivalently: $\Risk_p(\widehat{x}(\cdot)|\epsilon,\sigma,s,\upsilon)\leq\delta$ if and only if there exists a set $\Xi$ of ``good'' realizations of $\xi$ with $\Prob\{\xi\in\Xi\}\geq1-\epsilon$ such that whenever $\xi\in\Xi$, one has $\|\widehat{x}(Ax+\sigma\xi+u)-x\|_p\leq\delta$ {\sl for all $x\in X(s,\upsilon)$ and all $u\in\cU$.}
\paragraph{Norm $\nu(\cdot)$.} Given $\e$ and $\sigma>0$ let us denote
\be
\nu(v)=\nu_{\e,\sigma,\cU}(v)=\sup_{u\in \cU}u^Tv+\sigma\sqrt{2\ln (n/\e)}\|v\|_2.
\ee{nunorm}
Since $\cU$ is convex, closed and symmetric with respect to the origin, $\nu(\cdot)$ is a norm.
Let  $\nu_*$ be the norm on $\bR^n$ conjugate to $\nu$:
\bse
\nu_*(u)=\max\limits_{v}\{v^Tu:\nu(v)\leq 1\}.
\ese
\paragraph{Conditions $\bH(\gamma)$ and $\bH_{s,\infty}(\kappa)$.}
Let $\gamma=(\gamma_1,...,\gamma_n)\in\bR^n_+$. Given $A\in\bR^{m\times n}$, consider the following condition on a matrix $H=[h_1,...,h_n]\in \bR^{m\times n}$:
\begin{quote} {\em $\bH(\gamma)$:
for all $x\in \bR^n$ and $1\le i\leq n$ one has}
\be
|x_i|\le |h^T_iAx|+\gamma_i\|x\|_1.
\ee{cond00}
\end{quote}
Now let $s$ be a positive integer and $\kappa>0$. Given $A\in \bR^{m\times n}$, we say that a matrix $H=[h_1,...,h_n]\in \bR^{m\times n}$ satisfies condition $\bH_{s.\infty}(\kappa)$ \footnote{The reason for this cumbersome, at the first glance, notation will become clear later, in Section \ref{secextent}.}, if
\begin{equation}\label{Hsinfty}
\forall x\in\bR^n: \|x\|_\infty\leq \|H^TAx\|_\infty +s^{-1}\kappa\|x\|_1.
\end{equation}
The conditions we have introduced are closely related to each other:
\begin{lemma}\label{newlemma} If $H$ satisfies $\bH(\gamma)$, then $H$ satisfies $H_{s,\infty}(s\|\gamma\|_\infty)$, and ``nearly vice versa:'' given $H\in\bR^{m\times n}$ satisfying $\bH_{s,\infty}(\kappa)$, one can build efficiently a matrix $H'\in\bR^{m\times n}$ satisfying $\bH(\gamma)$ with $\gamma={\kappa\over s}[1;...;1]$ (i.e., $\kappa=s\|\gamma\|_\infty$) and such that  the columns of $H'$ are convex combinations of the columns of $H$ and $-H$, so that $\nu(H')\leq\nu(H)$ for every norm $\nu(\cdot)$ on $\bR^m$.
\end{lemma}
\subsection{Regular $\ell_1$ Recovery}\label{sectregrec}
 In this section we discuss the properties of
 the {\em regular  $\ell_1$-recovery} $\wh{x}_{\reg}$ given by:
\be
\wh{x}_{\reg}=\wh{x}_{\reg}(y)\in \Argmin\limits_{v\in \bR^n} \{\|v\|_1: |h_i^T(Av-y)|\le \rho_i,\;\;i=1,...,n \},
\ee{regular}
 where $y $ is as in \rf{obs}, $h_i$, $i=1,...,n$, are some vectors in $\bR^m$ and $\rho_i>0$, $i=1,...,n$. We refer to the matrix $H=[h_1,...,h_n]$ as to the {\sl contrast matrix} underlying the recovering procedure.
\par
The starting point of our developments is the following
\begin{proposition}\label{prop10} Given an $m\times n$ sensing matrix $A$, noise intensity $\sigma$, uncertainty set $\cU$ and  a tolerance $\epsilon\in(0,1)$, let the matrix
$H=[h_1,...,h_n]$ from {\rm (\ref{regular})} satisfy the condition $\bH(\gamma)$ for some $\gamma\in\bR^n_+$, and let $\rho_i$ in {\rm (\ref{regular})}
satisfy the relation
\begin{equation}\label{rhoi}
\rho_i\geq\nu_i:=\nu(h_i),\,i=1,...,n
\end{equation}
where $\nu(\cdot)$ is given by {\rm (\ref{nunorm})}. Then there exists a set $\Xi\subset\bR^m$, $\Prob\{\xi\in\Xi\}\geq1-\epsilon$, of ''good'' realizations of $\xi$ such that
\par
{\rm (i)}
 Whenever $\xi\in\Xi$, for every $x\in\bR^n$, every $u\in\cU$ and every subset $I\subset\{1,...,n\}$ such that
\begin{equation}\label{suchthat}
\gamma_{I}:=\sum_{i\in I}\gamma_i<\half,
\end{equation}
the regular $\ell_1$-recovery $\wh{x}_\reg$ given by {\rm (\ref{regular})} satisfies:
\begin{eqnarray}
(a)&\|\wh{x}_{\reg}(Ax+\sigma\xi+u)-x\|_1&\leq {2\|x_{J}\|_1+2\rho_I+2\nu_I\over 1-2\gamma_{I}};\vspace{0.2cm}\nn
(b)&|\left[\wh{x}_{\reg}(Ax+\sigma\xi+u)-x\right]_i|&\le\rho_i+\nu_i+\gamma_i\|\wh{x}_{\reg}(y)-x\|_1\label{r100}\\
&&\leq \rho_i+\nu_i+\gamma_i{2\|x_{J}\|_1+2\rho_I+2\nu_I\over 1-2\gamma_{I}},\,i=1,...,n,\nonumber
\end{eqnarray}
where
$\rho_I=\sum_{i\in I}\rho_i$  and $\nu_{I}=\sum_{i\in I} \nu_i.$
\par
{\rm (ii)} In particular, when setting
\begin{equation}\label{whensetting}
\begin{array}{l}
\widehat{\rho}_s=\|[\rho_1;\,...;\,\rho_n]\|_{s,1},\;\;\widehat{\nu}_s=\|[\nu(h_1);...;\nu(h_n)]\|_{s,1},\;\;
\widehat{\gamma}_s=\|[\gamma_1;\,...;\,\gamma_n]\|_{s,1},\\
\widehat{\rho}=\widehat{\rho}_1=\max_i\rho_i,\;\;\nu(H)=\widehat{\nu}_1=\max_i\nu(h_i),\;\;\widehat{\gamma}=\widehat{\gamma}_1=\max_i\gamma_i,
\end{array}
\end{equation}
and assuming $\widehat{\gamma}_s<\half$, for every $x\in\bR^n$, $\xi\in\Xi$ and $u\in \cU$ it holds
\bse
\|\wh{x}_{\reg}(Ax+\sigma\xi+u)-x\|_1&\leq &2{\|x-x^s\|_1+\widehat{\rho}_s+\widehat{\nu}_s\over 1-2\widehat{\gamma}_s}\leq 2{\|x-x^s\|_1\over
1-2\widehat{\gamma}_s}+2s{\widehat{\rho}+\nu(H)\over
1-2\widehat{\gamma}_s};\\
\|\wh{x}_{\reg}(Ax+\sigma\xi+u)-x\|_\infty&\le&2\widehat{\gamma}{\|x-x^s\|_1\over 1-2\widehat{\gamma}_s}+
{[1+2s\widehat{\gamma}-2\widehat{\gamma}_s][\widehat{\rho}+\nu(H)]\over 1-2\widehat{\gamma}_s}
\ese
\par (iii) Finally, assuming $s\widehat{\gamma}<1/2$, for every $\xi\in {\crd\Xi}$, $x\in\bR^n$ and $u\in\cU$ one has
\be
\begin{array}{rcl}
\|\wh{x}_{\reg}(Ax+\sigma\xi+u)-x\|_1&\leq& 2{\|x-x^s\|_1\over
1-2s\widehat{\gamma}}+2s{\widehat{\rho}+\nu(H)\over
1-2s\widehat{\gamma}};\\
\|\wh{x}_{\reg}(Ax+\sigma\xi+u)-x\|_\infty&\leq& s^{-1}{\|x-x^s\|_1\over 1-2s\widehat{\gamma}}+
{\widehat{\rho}+\nu(H)\over 1-2s\widehat{\gamma}}.
\end{array}
\ee{therr}
\end{proposition}
The following result is an immediate corollary of Proposition \ref{prop10}:
\begin{lemma}\label{allp}
Under the premise of Proposition \ref{prop10}, assume that  $\widehat{\gamma}_s<\half$. Then for all $1\le p\le \infty$ and $\upsilon\geq0$:
\begin{equation}\label{rl1}
\Risk_p(\wh{x}_{\reg}(\cdot)|\e,\sigma,s,\upsilon)\le  {2\over1-2\widehat{\gamma}_s}\big[
\upsilon+\widehat{\rho}_s+\widehat{\nu}_s\big]^{{1\over p}}\big[\widehat{\gamma}\upsilon+[\half-\widehat{\gamma}_s][\widehat{\rho}+\nu(H)]+\widehat{\gamma}[\widehat{\nu}_s+\widehat{\rho}_s]\big]^{{p-1\over p}}
\end{equation}
(for notation, see (\ref{whensetting})).
Further, if $s\widehat{\gamma}<1/2$, we have also
\begin{equation}\label{also}
1\leq p\leq\infty\Rightarrow \Risk_p(\wh{x}_{\reg}(\cdot)|\e,\sigma,s,\upsilon)\le
{(2s)^{1\over p}\over 1-2s\widehat{\gamma}}(s^{-1}\upsilon+\widehat{\rho}+\nu(H)).
\end{equation}
\end{lemma}
The next statement is similar to the cases of $\kappa:=s\widehat{\gamma}<1/2$ in Proposition \ref{prop10} and Lemma \ref{allp}; the difference is that
now we assume that $H$ satisfies $\bH_{s,\infty}(\kappa)$, which, by Lemma \ref{newlemma}, is a weaker requirement on $H$ than to satisfy $\bH(\gamma)$ with $s\widehat{\gamma}\equiv s\|\gamma\|_\infty=\kappa$.
\begin{proposition}\label{prop101}
Given an $m\times n$ sensing matrix $A$, noise intensity $\sigma$, uncertainty set $\cU$ and  a tolerance $\epsilon\in(0,1)$, let the matrix
$H=[h_1,...,h_n]$ from {\rm (\ref{regular})} satisfy the condition $\bH_{s,\infty}(\kappa)$ for some $\kappa<1/2$, and let $\rho_i$ in (\ref{regular})
satisfy the relation (\ref{rhoi}).  Then there exists a set $\Xi\subset\bR^m$, $\Prob\{\xi\in\Xi\}\geq1-\epsilon$, of ''good'' realizations of $\xi$ such that
whenever $\xi\in\Xi$, for every $x\in\bR^n$ and every $u\in\cU$ one has
\be
\begin{array}{rcl}
\|\wh{x}_{\reg}(Ax+\sigma\xi+u)-x\|_1&\leq& 2{\|x-x^s\|_1\over
1-2\kappa}+2s{\widehat{\rho}+\nu(H)\over
1-2\kappa};\\
\|\wh{x}_{\reg}(Ax+\sigma\xi+u)-x\|_\infty&\leq& s^{-1}{\|x-x^s\|_1\over 1-2\kappa}+
{\widehat{\rho}+\nu(H)\over 1-2\kappa}.
\end{array}
\ee{therrnew}
In particular,
\begin{equation}\label{alsonew}
1\leq p\leq\infty \Rightarrow \Risk_p(\wh{x}_{\reg}(\cdot)|\e,\sigma,s,\upsilon)\le
{(2s)^{1\over p}\over 1-2\kappa}(s^{-1}\upsilon+\widehat{\rho}+\nu(H)).
\end{equation}
\end{proposition}

\subsection{Penalized $\ell_1$ Recovery}\label{sectpenrec}
Now consider the {\sl penalized $\ell_1$-recovery} $\wh{x}_\pen$ as follows:
\be
\wh{x}_\pen(y)\in \Argmin\limits_{v\in \bR^n} \{\|v\|_1+\theta s\|H^T(Av-y)\|_\infty\},
 \ee{penalized}
 where $y $ is as in \rf{obs}, and an integer $s\leq n$, a positive $\theta$, and a matrix $H$ are parameters of the construction.
 \begin{proposition}\label{prop1}
 Given an $m\times n$ sensing matrix $A$, an integer $s\leq n$,  a matrix $H=[h_1,...,h_n]\in\bR^{m\times n}$ and positive reals
 $\gamma_i$, $1\le i\le n$, satisfying the condition $\bH(\gamma)$, and a $\theta>0$, assume that
 \begin{equation}\label{gammasmall}
\widehat{\gamma}_s:=\|\gamma\|_{s,1}<\half
 \end{equation}
and
\begin{equation}\label{theta}
(1-\widehat{\gamma}_s)^{-1}<\theta<(\widehat{\gamma}_s)^{-1}
\end{equation}
Further, let $\sigma\geq0$, $\epsilon\in(0,1)$, and let
\begin{equation}\label{nui}
\nu_i=\nu_{\epsilon,\sigma,\cU}(h_i),\,i=1,...,n,\,\nu(H)=\max\limits_i\nu_i.
  \end{equation}
  Consider the penalized recovery $\widehat{x}_\pen(\cdot)$ associated with $H,s,\theta$.
  There exists a set $\Xi\subset\bR^m$, $\Prob\{\xi\in\Xi\}\geq1-\epsilon$, of ''good'' realizations of $\xi$ such that
  \par {\rm (i)} Whenever $\xi\in\Xi$, for every signal $x\in\bR^n$ and every $u\in\cU$ one has
\begin{equation}\label{r200}
\begin{array}{lrcl}
(a)&\|\widehat{x}_\pen(Ax+\sigma\xi+u)-x\|_1&\leq& {2\|x-x^s\|_1+2s\theta\nu(H)\over \min[\theta(1-\widehat{\gamma}_s)-1,1-\theta \widehat{\gamma}_s]}\vspace{0.2cm}\\
(b)&\|\widehat{x}_\pen(Ax+\sigma\xi+u))-x\|_\infty&\leq&
\left({1\over s\theta}+\widehat{\gamma}\right)\|\widehat{x}_\pen(Ax+\sigma\xi+u)-x\|_1+2\nu(H)\vspace{0.2cm}\\
&&
\le&2\left({1\over s\theta}+\widehat{\gamma}\right){\|x-x^s\|_1\over \min[\theta(1-\widehat{\gamma}_s)-1,1-\theta \widehat{\gamma}_s]}
+2\nu(H)\left[{1+s\theta\widehat{\gamma}\over\min[\theta(1-\widehat{\gamma}_s)-1,1-\theta \widehat{\gamma}_s]}+1\right],\\
\end{array}
\end{equation}
where, as in Lemma \ref{allp}, $\widehat{\gamma}=\max\limits_i\gamma_i$.
\par
{\rm (ii)}
 When $\theta=2$ and $\widehat{\gamma}<{1\over 2s}$, one has for every $x\in \bR^n$, $u\in \cU$ and $\xi\in \Xi$:
\begin{equation}\label{r200a}
\begin{array}{lrcl}
(a)&\|\widehat{x}_\pen(Ax+\sigma\xi+u)-x\|_1&\leq& 2{\|x-x^s\|_1\over 1-2s\widehat{\gamma}}+4s{\nu(H)\over 1-2s\widehat{\gamma}}\\
(b)&\|\widehat{x}_\pen(Ax+\sigma\xi+u))-x\|_\infty&\leq&
2s^{-1}{\|x-x^s\|_1\over 1-2s\widehat{\gamma}}
+4{\nu(H)\over1-2s\widehat{\gamma}},\\
\end{array}
\end{equation}
whence for every
  $\upsilon\geq0$ and $1\le p\le \infty$:
\begin{equation}\label{risks}
\Risk_p(\widehat{x}_\pen(\cdot)|\e,\sigma,s,\upsilon)
\leq {2s^{{1\over p}}\over 1-2s\widehat{\gamma}}(s^{-1}\upsilon+2\nu(H)).
\end{equation}
\end{proposition}
 The next statement is in the same relation to Proposition \ref{prop1} as Proposition \ref{prop101} is to Proposition \ref{prop10}{ and Lemma \ref{allp}.
 \begin{proposition}\label{prop11}
 Given an $m\times n$ sensing matrix $A$, noise intensity $\sigma$, uncertainty set $\cU$ and  a tolerance $\epsilon\in(0,1)$, let the matrix
$H=[h_1,...,h_n]$ from {\rm (\ref{penalized})} satisfy the condition $\bH_{s,\infty}(\kappa)$ for some $\kappa<1/2$, and let $\theta=2$.  Then there exists a set $\Xi\subset\bR^m$, $\Prob\{\xi\in\Xi\}\geq1-\epsilon$, of ''good'' realizations of $\xi$ such that
whenever $\xi\in\Xi$, for every $x\in\bR^n$ and every $u\in\cU$ one has
\be
\begin{array}{rcl}
\|\wh{x}_{\pen}(Ax+\sigma\xi+u)-x\|_1&\leq& 2{\|x-x^s\|_1\over
1-2\kappa}+4s{\nu(H)\over
1-2\kappa};\\
\|\wh{x}_{\pen}(Ax+\sigma\xi+u)-x\|_\infty&\leq& 2s^{-1}{\|x-x^s\|_1\over 1-2\kappa}+
4{\nu(H)\over 1-2\kappa}.
\end{array}
\ee{therrnewnew}
In particular,
\begin{equation}\label{alsonewnew}
1\leq p\leq\infty \Rightarrow \Risk_p(\wh{x}_{\pen}(\cdot)|\e,\sigma,s,\upsilon)\le
{2s^{{1\over p}}\over 1-2\kappa}(s^{-1}\upsilon+2\nu(H)).
\end{equation}
\end{proposition}
Note that under the premise of Proposition \ref{prop101}, the smallest possible values of $\rho_i$ are the quantities $\nu_i$, which results in $\widehat{\rho}=\nu(H)$; with this choice of $\rho_i$, the risk bound for the regular recovery, as given by the right hand side in (\ref{alsonew}), coincides within factor 2 with the risk bound for the penalized recovery with $\theta=2$ as given by (\ref{alsonewnew}); both bounds assume that $H$ satisfies $\bH_{s,\infty}(\kappa)$ with $\kappa<1/2$ and imply that
\begin{equation}\label{taketheform}
 1\leq p\leq\infty \Rightarrow \Risk_p(\wh{x}(\cdot)|\e,\sigma,s,\upsilon)\le
2{s^{{1\over p}}\over 1-2\kappa}(s^{-1}\upsilon+2\nu(H)).
\end{equation}
When $\upsilon=0$, the latter bound admits a quite transparent interpretation: everything is {\sl as if} we were observing the sum of an unknown {\sl $s$-dimensional} signal and an observation error of the uniform norm $O(1)\nu(H)$.
 \section{Efficient construction of the contrast matrix $H$}\label{sectoptH}\label{BuildingH}
\label{matrixh}
In what follows, we fix $A$, the ``environment parameters'' $\epsilon,\,\sigma,\,\cU$ and the ``level of sparsity''  $s$ of signals $x$ we intend to recover, and are interested in building the contrast matrix $H=[h_1,...,h_n]$ resulting in as small as possible error bound (\ref{taketheform}). All we need to this end is to answer the following question (where we should specify the norm $\varphi(\cdot)$ as $\nu_{\epsilon,\sigma,\cU}(\cdot)$):
\begin{quote}
(?) {\sl Let $\varphi(\cdot)$ be a norm on $\bR^m$, and $s$ be a positive integer. What is the domain $G_s$ of pairs
$(\omega,\kappa)\in\bR^2_+$ such that $\kappa<1/2$ and  there exists matrix $H=[h_1,...,h_n]\in\bR^{m\times n}$   satisfying the condition $\bH_{s,\infty}(\kappa)$ and the relation $\varphi(H):=\max\limits_i\varphi(h_i)\leq\omega$? How to find such an $H$, provided it exists?}
\end{quote}
Invoking Lemma \ref{newlemma}, we can reformulate this question as follows:
\begin{quote}
(??) {\sl Let $\varphi(\cdot)$ and $s$ be as in (?). Given $(\omega,\kappa)\in\bR^2_+$, how to find vectors $h_i\in\bR^m$, $1\leq i\leq n$, satisfying
$$
(a):\ \varphi(h_i)\leq \omega;\ \&\  (b):\ |x_i|\leq |h_i^TAx|+s^{-1}\kappa \|x\|_1\,\forall x\in\bR^n\eqno{(\cP_i)}
$$
for every $i$, or to detect correctly that no such collection of vectors exists?}
\end{quote}
Indeed, by Lemma \ref{newlemma}, if $H'$ satisfies $\bH_{s,\infty}(\kappa)$ and $\varphi(H')\leq\omega$, then there exists $H=[h_1,...,h_n]$
such that $h_i$ satisfy $(\cP_i.b)$ for all $i$ and $\varphi(H)\leq\varphi(H')\leq\omega$, so that $h_i$ satisfy $(\cP_i.a)$ for all $i$ as well. Vice versa,
if $h_i$ satisfy $(\cP_i)$, $1\leq i\leq n$, then the matrix $H=[h_1,...,h_n]$ clearly satisfies $\bH_{s,\infty}(\kappa)$, and $\varphi(H)\leq\omega$.
\par
The answer to (??) is given by the following
\begin{lemma}\label{arikn} Given  $\kappa>0$, $\omega\geq0$, and a positive integer $s$, let $\gamma=\kappa/s$. For every $i\leq n$, the following  three properties are equivalent to each other:
\begin{itemize}
\item[{\rm (i)}] There exists $h=h_i$ satisfying $(\cP_i)$;
\item[{\rm (ii)}] The optimal value in the optimization problem
$$
{\rm \Opt}_i(\gamma)=\min\limits_h\left\{\varphi(h): \,\|A^Th-e_i\|_\infty\leq \gamma\right\}\eqno{(P_i^\gamma)}
$$
where $e_i$ is $i$-th standard basic orth in $\bR^n$, is $\leq\omega$;
\item[{\rm (iii)}] One has
\begin{equation}\label{omehas12}
\forall x\in\bR^n: |x_i|\leq \omega \varphi_*(Ax)+\gamma\|x\|_1,
\end{equation}
where $\varphi_*(u)=\max\limits_{\varphi(v)\leq1}u^Tv$ is the norm on $\bR^m$ conjugate to $\varphi(\cdot)$.
\end{itemize}
Whenever one (and then -- all) of these properties take place, problem $(P_i^\gamma)$ is solvable, and its optimal solution $h_i$ satisfies $(\cP_i)$.
\end{lemma}
\subsection{Optimal contrasts for regular and penalized recoveries}
As an immediate consequence of Lemma \ref{arikn}, we get the following description of the domain $G_s$ associated with the norm $\varphi(\cdot)=\nu_{\e,\sigma,\cU}(\cdot)$:
\begin{equation}\label{domainG}
\begin{array}{lrcl}
(a)&G_s&=&\left\{(\omega,\kappa)\geq0: s^{-1}\kappa\geq \gamma_*,\,\omega\geq \omega_*(s^{-1}\kappa)\right\},\\
\multicolumn{4}{l}{\hbox{where}}\\
(b)&\gamma_*&=&\max\limits_{1\leq i\leq n}\min\limits_h\|A^Th-e_i\|_\infty
=\max\limits_i\max\limits_x\left\{x_i:\|x\|_1\leq1,Ax=0\right\},\\
(c)&\omega_*(\gamma)&=&\max\limits_{1\leq i\leq n}\Opt_i(\gamma)\\
\end{array}
\end{equation}
where $\phi(\cdot)$ in $(P_i^\gamma)$ is specified as $\nu_{\e,\sigma,\cU}(\cdot)$. Note that
the second equality in $(b)$ is given by Linear Programming duality.
Indeed, by $(b)$, $\gamma_*$ is the smallest $\gamma$ for which all problems $(P_i^\gamma)$, $i=1,...,n$, are feasible, and thus, by Lemma \ref{arikn}, $(\gamma,\kappa)\in G_s$ if and only if
$\kappa/s\geq\gamma_*$ and $\omega\geq\omega_*(\kappa/s)$.
 \par Note that the quantity $\gamma_*$  depends solely of $A$, while $\omega_*(\cdot)$ depends on $\epsilon,\sigma,\cU$, as on parameters, but is independent of $s$.
\par
The outlined results suggest the following scheme of building the contrast matrix $H$:
\begin{itemize}
\item we compute $\gamma_*$ by solving $n$ Linear Programming problems in (\ref{domainG}.$b$); if $s\gamma_*\geq\half$, then $G_s$ does not contain points $(\omega,\kappa)$ with $\kappa<1/2$, so that our recovery routines are not applicable (or, at least, we cannot justify them theoretically);
\item when $s\gamma_*<\half$, the set $G_s$ is nonempty, and its Pareto frontier (the set of pairs $(\omega,\kappa)\in\bR^2_+$ such that $(\omega,\kappa)\geq(\omega',\kappa)\in G_s$ is possible if and only if $\omega'=\omega$) is the curve $(\omega_*(\gamma),s\gamma)$, $\gamma_*\leq\gamma<{1\over 2s}$. We choose a ``working point'' on this curve, that is, a point $\bar{\gamma}\in[\gamma_*,{1\over 2s}]$ and compute $\omega_*(\bar{\gamma})$ by solving the convex optimization programs $(P^{\bar{\gamma}}_i)$, $i=1,...,n$, with $\phi(\cdot)$ specified as $\nu_{\e,\sigma,\cU}(\cdot)$. $\omega_*(\bar{\gamma})$ is nothing but the maximum, over $i$, of the optimal values of these problems, and the optimal solutions $h_i$ to the problems induce the matrix $H=H(\bar{\gamma})=[h_1,...,h_n]$ which satisfies $\bH_{s,\infty}(s\bar{\gamma})$ and has $\nu(H)\leq \omega_*(\bar{\gamma})$. By reasoning which led us to (??),
    $$\nu(H(\bar{\gamma}))=\omega_*(\bar{\gamma}) =
    \min_{H'}\left\{\nu(H'): H'\hbox{\ satisfies\ } \bH_{s,\infty}(s\bar{\gamma})\right\},$$
that is, $H=H(\bar{\gamma})$ is the best for our purposes contrast matrices satisfying $\bH_{s,\infty}(s\bar{\gamma})$. With this contrast matrix, the error bound (\ref{taketheform}) for regular/penalized $\ell_1$ recoveries (in the former, $\rho_i=\nu_i(h_i)$, in the latter, $\theta=2$) read
\begin{equation}\label{risksrisks}
 1\leq p\leq\infty \Rightarrow \Risk_p(\wh{x}(\cdot)|\e,\sigma,s,\upsilon)\le
2{s^{{1\over p}}\over 1-2s\bar{\gamma}}(s^{-1}\upsilon+2\omega_*(\bar{\gamma})).
\end{equation}
\end{itemize}
The outlined strategy does not explain how to choose $\bar{\gamma}$. This issue could be resolved, e.g., as follows. We choose an upper bound on the sensitivity of the risk (\ref{risksrisks}) to $\upsilon$, i.e., to the $\|\cdot\|_1$-deviation of a signal to be recovered from the set of $s$-sparse signals. This sensitivity is proportional to ${1\over 1-2s\bar{\gamma}}$, so that an upper bound on the sensitivity translates into an upper bound $\gamma^+<{1\over 2s}$ on $\bar{\gamma}$. We can now choose $\bar{\gamma}$ by minimizing the remaining term in the risk bound over $\gamma\in[\gamma_*,\gamma^+]$, which amounts to solving the optimization problem
$$
\max\left\{\tau: \;\tau \omega_*(\gamma)\leq 1-2s\gamma, \;\gamma_*\leq \gamma \leq \gamma^+\right\}.
$$
Observing that $\omega_*(\cdot)$ is, by its origin, a convex function, we can solve the resulting problem efficiently by bisection in $\tau$. A step of this bisection requires solving a {\sl univariate} convex feasibility problem with efficiently computable constraint and thus is easy, at least for moderate values of $n$.
\section{Range of feasibility of condition $\bH_{s,\infty}(\kappa)$}\label{boundingomega}\label{LimitsofPerformance}

We address the crucial question of what can be said about the magnitude of the
quantity $\omega_*(\cdot)$, see (\ref{domainG}) and the risk bound (\ref{risksrisks}). One way to answer it is just to compute the (efficiently computable!) quantity $\omega_*(\gamma)$ for a desired value of $\gamma$.
Yet it is natural to know theoretical upper bounds on $\omega_*$ in some ``reference'' situations. Below, we provide three results of this type.
\par
At this point, it makes sense to express in the notation that $\omega_*(\gamma)$ depends, as on parameters, on the sensing matrix $A$ and the ``environment parameters'' $\epsilon,\sigma,\cU$, so that in this section  we write $\omega_*(\gamma|A,\epsilon,\sigma,\cU)$ instead of $\omega_*(\gamma)$.
\subsection{Bounding $\omega_*(\cdot)$ via mutual incoherence}\label{secmuinc}\label{sectmuinc}
Recall that for an $m\times n$ sensing matrix $A=[A_1,...,A_n]$ with no zero columns, its {\sl mutual incoherence} is defined as
$$
\mu(A)=\max\limits_{1\leq i\neq j\leq n} {|A_i^TA_j|\over A_i^TA_i}.
$$
Compressed Sensing literature contains numerous mutual-incoherence-related results (see, e.g.,  \cite{DE,EB,GN} and references therein). To the best of our knowledge, all these results state that  if $s$ is a positive integer and $A$ is a sensing matrix such that ${s\mu(A)\over\mu(A)+1}<\half$, then $\ell_1$-based sparse recovery is well suited for recovering $s$-sparse signals (e.g., recovers them exactly when there is no observation noise, admit explicit error bounds when there is noise and/or the signal is only nearly $s$-sparse, etc.). To the best of our knowledge, all these results, up to the values of absolute constant factors in error bounds, are covered by the risk bounds (\ref{risksrisks}) combined with the following immediate observation:
\begin{observation}\label{mutualincoh} Whenever $A=[A_1,...,A_n]$ is an $m\times n$ matrix with no zero columns and $s$ is a positive integer, the matrix $H(A)={1\over\mu(A)+1}[A_1/A_1^TA_1,A_2/A_2^TA_2,...,A_n/A_n^TA_n]$ satisfies the condition $\bH_{s,\infty}\left({s\mu(A)\over\mu(A)+1}\right)$.
\end{observation}
Verification is immediate: the diagonal entries in the matrix $Z=I-H^TA$ are equal to $\gamma:=1-{1\over \mu(A)+1}={\mu(A)\over\mu(A)+1}$, while the magnitudes of the off-diagonal entries in $Z$  do not exceed $\gamma$. Therefore
$$
\begin{array}{l}
x\in\bR^n\Rightarrow \gamma\|x\|_1\geq \|Zx\|_\infty=\|x-H^TAx\|_\infty\geq \|x\|_\infty-\|H^TAx\|_\infty\Leftrightarrow \|x\|_\infty\leq \|H^TAx\|_\infty+\gamma \|x\|_1\\
\Leftrightarrow
\hbox{\ $H$ satisfies $\bH_{s,\infty}(s\gamma)$}.\\
\end{array}
$$
Observe that the Euclidean norms of the columns in $H(A)$ do not exceed $\left[\min\limits_i\|A_i\|_2\right]^{-1}$, whence $\nu(H(A))\leq r(\cU)+{\sigma\sqrt{2\ln(n/\epsilon)}\over\min_i\|A_i\|_2}$, where $r(\cU)=\max\limits_{u\in\cU}\|u\|_2$.  In the notation from Section \ref{matrixh}, our observations can be summarized as follows:
\begin{corollary}\label{cormuin} For every $m\times n$ matrix $A$ with no zero columns, one has $\gamma_*\leq \gamma:={\mu(A)\over\mu(A)+1}$ and $\omega_*(\gamma|A,\epsilon,\sigma,\cU)\leq\nu(H(A))\leq r(\cU)+{\sigma\sqrt{2\ln(n/\epsilon)}\over\min_i\|A_i\|_2}$. In particular,
$$
s\leq {\mu(A)+1\over 3\mu}\Rightarrow \omega_*({1\over 3s}|A,\epsilon,\sigma,\cU)\leq r(\cU)+{\sigma\sqrt{2\ln(n/\epsilon)}\over\min_i\|A_i\|_2}.
$$
\end{corollary}
 It should be added that as $m,n$ grow in such a way that $\ln(n)\leq O(1)\ln m$, realizations $A$ of ``typical'' random $m\times n$ matrices (e.g., those with independent $\cN(0,1/m)$ entries or with independent entries taking values $\pm1/\sqrt{m}$) with overwhelming probability satisfy $\mu(A)\leq O(1)\sqrt{\ln(n)/m}$ and $\|A_i\|_2\leq 0.9$ for all $i$. By Corollary \ref{cormuin}, it follows that for these $A$ the condition
 $\bH_{s,\infty}(\kappa)$ with, say, $\kappa=1/3$ can be satisfied for $s$ as large as $O(1)\sqrt{m/\ln(n)}$ merely by the choice $H=H(A)$, which ensures that $\nu(H)\leq O(1)[r(\cU)+\sigma\sqrt{2\ln(n/\epsilon)}]$; in particular, in the indicated range of values of $s$ one has $\omega_*({1\over 3s})\leq O(1)[r(\cU)+\sigma\sqrt{2\ln(n/\epsilon)}]$.
\subsection{The case of $A$ satisfying the Restricted Isometry Property}
 \begin{proposition} \label{propRIP} Let $A$ satisfy $\RIP(\delta,k)$ with some $\delta\in(0,1)$ and with $k>1$. Then there exists matrix $H(A)$ which, for every positive integer $s$, satisfies  the condition $\bH_{s,\infty}(s\gamma(\delta,k))$, with
\begin{equation}\label{k}
\gamma(\delta,k)={\sqrt{2}\delta\over(1-\delta)\sqrt{k-1}},
\end{equation}
and is such that $\nu(H(A))\leq
\left[r(\cU)+\sigma\sqrt{2\ln(n/\epsilon)}\right]/\sqrt{1-\delta}$. In particular,
\begin{equation}\label{GammaRip}
s\leq {1-\delta\over 3\sqrt{2}}\sqrt{k-1}\Rightarrow \omega_*({1\over 3s}|A,\epsilon,\sigma,\cU)\leq{1\over \sqrt{1-\delta}}\left[r(\cU)+\sigma\sqrt{2\ln(n/\epsilon)}\right].
\end{equation}
\end{proposition}

\subsection{Oracle inequality}
Here we assume that $A\in\bR^{m\times n}$ possesses the following property (where $S$ is a positive integer and $\varphi>0$):
\begin{quote} $\bO(S,\omega)$:  For every $i\in\{1,...,n\}$ and every $S$-element subset $I\ni i$ of $\{1,...,n\}$ there exists a routine
   $\cR_{i,I}$ for recovering $x_i$ from a noisy observation
\[
y=Ax+u+\sigma e,\eqno{[e\sim\cN(0,I_m), u\in\cU]}
\]
of {\sl unknown} signal $x\in \bR^n$,
{\sl known} to be supported on $I$ such that for every such signal and every $u\in\cU$ one has
\[
\Prob\{|\cR_{i,I}(Ax+u+\sigma e)-x_i|\ge \omega\}\le \e.
\]
\end{quote}
\par
We intend to demonstrate that in this situation for all $s$ in certain range (which extends as $S$ grows and $\omega$ decreases) the uniform error of the regular and the penalized recoveries associated with properly selected contrast matrix is, with probability $\geq1-\epsilon$, ``close'' to $\omega$.
The precise statement is as follows:
\begin{proposition}\label{propOracle} Given $A$ and the ``environment parameters'' $\epsilon<1/16$, $\sigma$, $\cU$, assume that $A$ satisfies the condition $\bO(S,\gamma)$ with certain $S,\gamma$. Then for every integer $s$ from the range
\begin{equation}\label{range}
1\leq s\leq {\sigma\sqrt{2S\ln(1/\e)}\over 4\omega\|A\|}
\end{equation}
(here $\|\cdot\|$ is the standard matrix norm, the largest singular value) there exists a contrast matrix $H$ satisfying the condition $\bH_{s,\infty}(\four)$ and such that $\nu(H)\leq 2\sqrt{1+\ln(n)/\ln(1/\epsilon)}\omega$, so that in the outlined range of values of $H$ one has
$\omega_*({1\over 4s})\leq 2\sqrt{1+\ln(n)/\ln(1/\epsilon)}\omega$, and the associated with $H$ error bound (\ref{risksrisks}) for regular/penalized
$\ell_1$ recovery is
\begin{equation}\label{riskagain}
\Risk_p(\widehat{x}(\cdot)|\epsilon,\sigma,s,\upsilon)\leq
16s^{{1\over p}}\left[\omega\sqrt{1+{\ln n\over \ln(1/\e)}}+{\upsilon\over 4s}\right].
\end{equation}
\end{proposition}
Proposition \ref{propOracle} justifies to some extent, our approach; it says that if
there exists a routine which recovers $S$-sparse signals {\em with a priori known sparsity pattern} within certain accuracy
(measured component-wise), then our recovering routines exhibit ``close'' performance
without any knowledge of
the sparsity pattern, albeit in a smaller range of values of the sparsity parameter.
\subsection{Condition $\bH_{s,\infty}(\kappa)$: limits of performance}
Recall that when recovering $s$-sparse signals, the  condition $\bH_{s,\infty}(\kappa)$ helps {\sl only when $\kappa<1/2$}. Unfortunately, with these $\kappa$, the condition is feasible in a severely restricted range of values of $s$. Specifically, from \cite[Proposition 5.1]{JKKN} and Lemma \ref{newlemma} it immediately follows that
 \begin{quote}
 (*) {\sl If $A\in\bR^{m\times n}$ is not ``nearly square,'' that is, if $n>2(2\sqrt{m}+1)^2$, then the condition $\bH_{s,\infty}(\kappa)$ with $\kappa<1/2$ can{\em not} be satisfied when $s$ is ``large'', namely, when  $s>2\sqrt{2m}+1$}.
  \end{quote}
Note that from the discussion at the end of section \ref{secmuinc} we know that the ``$O(\sqrt{m})$ limit of performance'' of the condition $\bH_{s,\infty}(\cdot)$ stated in (*) is ``nearly sharp:'' -- when $s\leq O(1)\sqrt{m}$, the condition $\bH_{s,\infty}({1\over 3})$ associated with a typical randomly generated $m\times n$ sensing matrix $A$ is feasible and can be satisfies with a contrast matrix $H$ with quite moderate $\nu(H)$.
\par
  (*) says that unless $A$ is nearly square, the condition $\bH_{s,\infty}(\cdot)$ can validate $\ell_1$ sparse recovery only in a severely restricted
range $s\leq O(\sqrt{m})$ of values of the sparsity parameter. This is in sharp contrast with {\sl unverifiable} sufficient conditions for ``goodness'' of $\ell_1$ recovery, like RIP: it is well known that when $m,n$ grow, realizations of ``typical'' random $m\times n$ matrices, like those mentioned at the end of Section \ref{sectmuinc}, with overwhelming probability possess $\RIP(0.1,2s)$ with $s$ as large as $O(m/\ln(2n/m))$. As a result, ``unverifiable'' sufficient conditions, like RIP, can justify the validity of $\ell_1$ recovery routines in a much wider (and in fact -- the widest possible) range of values of the sparsity parameter $s$ than the ``fully computationally tractable'' condition $\bH_{s,\infty}(\cdot)$. This being said, note that this comparison is not completely fair. Indeed, aside of its tractability, the condition $\bH_{s,\infty}(\kappa)$ with $\kappa<1/2$ ensures the error bounds (\ref{risksrisks}) {\sl in the entire range $1\leq p\leq\infty$ of values of $p$}, which perhaps  is not the case with conditions like RIP. Specifically, consider the ``no nuisance'' case $\cU=\{0\}$, and let $A$ satisfy $\RIP(0.1,2S)$ for certain $S$. It is well known (see, e.g., the next section) that
in this case the Dantzig Selector recovery ensures for every $s\leq S$ and every $s$-sparse signal $x$ that
$$
\|\widehat{x}_{\ds}-x\|_p\leq O(1)\sigma\sqrt{\ln(n/\epsilon)}s^{1/p},\;1\leq p\leq 2,
$$
with probability $\geq1-\epsilon$. However, we are not aware of similar bounds (under whatever conditions) for ``large'' $s$ and $p>2$. For comparison: in the case in question, for ``small'' $s$, namely, $s\leq O(1)\sqrt{S}$, we have $\omega_*({1\over 3s})\leq O(1)\sigma\sqrt{\ln(n/\epsilon)}$ (by
Proposition \ref{propRIP}), whence for regular and penalized $\ell_1$ recoveries with appropriately chosen contrast matrix (which can be built efficiently!) one has for all $s$-sparse $x$
$$
\|\widehat{x}-x\|_p\leq O(1)\sigma\sqrt{\ln(n/\epsilon)}s^{{1\over p}}  \quad\forall p\in[1,\infty]
$$
with probability $\geq 1-\epsilon$ (see (\ref{risksrisks})). We wonder whether a similar (perhaps, with extra logarithmic factors) bound  can be obtained {\sl for large $s$} (e.g., $s\geq m^{{1\over 2}+\delta}$) for a {\sl whatever} $\ell_1$ recovery routine and a {\sl whatever} essentially nonsquare (say, $m<n/2$) $m\times n$ sensing matrix $A$ with columns of Euclidean length $\leq O(1)$.
\section{Extensions}\label{secextent}\label{Extensions}
We are about to demonstrate that the pivot element of the preceding sections --- the condition $\bH_{s,\infty}(\kappa)$ --- is the strongest (and seemingly the only verifiable one)
in a natural parametric series of conditions on a contrast matrix $H$; every one of these conditions validates the regular
and the penalized $\ell_1$ recoveries associated with $H$ in certain {\sl restricted} range of values of  $p$ in the error bounds (\ref{risksrisks}).
\subsection{Conditions $\bH_{s,q}(\kappa)$}
Let us fix an $m\times n$ sensing matrix $A$. Given a positive integer $s\leq m$, a $q\in[1,\infty]$ and a real $\kappa>0$, let us say that an $m\times n$ contrast matrix $H$ satisfies condition $\bH_{s,q}(\kappa)$, if
\begin{equation}\label{bHsq}
\forall x\in\bR^n: \|x\|_{s,q}\le s^{{1\over q}}\|H^TAx\|_\infty +\kappa s^{{1\over q}-1}\|x\|_1,
\end{equation}
where $\|x\|_{s,q}=\|x^s\|_q$ and $x^s$, as always, is the vector obtained from $x$ by zeroing all but the $s$ largest in magnitude entries. Observe that
\begin{itemize}
\item What used to be denoted $\bH_{s,\infty}(\kappa)$ before, is exactly what is called $\bH_{s,\infty}(\kappa)$ now;
\item If $H$ satisfies $\bH_{s,q}(\kappa)$, $H$ satisfies $\bH_{s,q'}(\kappa)$ for all $q'\in[1,q]$ (since for $s$-sparse vector $x^s$ we have $\|x^s\|_{q'}\leq s^{{1\over q'}-{1\over q}}\|x^s\|_{s,q}$).
\end{itemize}
Less immediate observations are as follows:
\begin{itemize}
\item Let $A$ be an $m\times n$ matrix and let $s\leq n$ be a positive integer. We say that $A$ is {\sl $s$-good} if for all $s$-sparse $x\in \bR^n$ the  $\ell_1$-recovery
\[
\wh{x}\in \Argmin_v \{\|v\|_1:Av=y\}
\]
is exact in the case of noiseless observation $y=Ax$. It turns out that feasibility of $\bH_{s,1}(\kappa)$ with $\kappa<\half$ is  intimately related to $s$-goodness of $A$:
\begin{lemma}\label{arik1} $A$ is $s$-good if and only if there exist $\kappa<\half$ and $H\in\bR^{m\times n}$ satisfying $\bH_{s,1}(\kappa)$.
\end{lemma}
\item The Restricted Isometry Property implies feasibility of $\bH_{s,2}(\kappa)$ with small $\kappa$:
\begin{lemma}\label{Dantzig} Let $A$ satisfy $\RIP(\delta,2s)$ with $\delta<{1\over 3}$. Then the matrix $H={1\over 1-\delta}A$ satisfies the condition $\bH_{s,2}(\kappa)$ with $\kappa={\delta\over 1-\delta}<\half$.
\end{lemma}
\end{itemize}
\subsection{Regular and penalized $\ell_1$ recoveries with contrast matrices satisfying $\bH_{s,q}(\kappa)$}
Our immediate goal is to obtain the following extension of the main results of Section \ref{ell1rout}, specifically,  Propositions \ref{prop101}, \ref{prop11}:
\begin{proposition}\label{arikwillnotlikeit}
Assume we are given an $m\times n$ sensing matrix $A=[a_1,...,a_n]$, an integer $s\leq m$, $\kappa<1/2$, a contrast matrix $H=[h_1,...,h_n]\in\bR^{m\times n}$, and $q\in[1,\infty]$ such that $H$ satisfies the condition $\bH_{s,q}(\kappa)$. Denote
$\nu_i=\nu_{\e,\sigma,\cU}(h_i)$, where the norm $\nu_{\e,\sigma,\cU}(\cdot)$ is defined in \rf{nunorm}, and $\nu(H)=\max_i \nu_i$. Let also noise intensity $\sigma$, uncertainty set  $\cU$ and tolerance $\epsilon\in(0,1)$ be given.
\par {\rm (i)} Consider the regular recovery \rf{regular} with the contrast matrix $H$ and the parameters $\rho_i$ satisfying the relations
$$
\rho_i\geq \nu_i,\,1\leq i\leq n,
$$
and let $\widehat{\rho}=\max\limits_i\rho_i$.
Then
\be
1\leq p\leq q\Rightarrow \Risk_p(\widehat{x}_{\reg}(\cdot)|\e,\sigma,s,\upsilon)\le (3s)^{1\over p}{\widehat{\rho}+\nu(H)+s^{-1}\upsilon\over 1-2\kappa}.
\ee{prop1i}
\par
{\rm (ii)} Consider the penalized recovery \rf{penalized} with  the contrast matrix $H$ and $\theta=2$. Then
\begin{equation}\label{proplii}
1\leq p\leq q\Rightarrow \Risk_p(\widehat{x}_{\pen}(\cdot)|\e,\sigma,s,\upsilon)\le 3s^{1\over p}{2\nu(H)+s^{-1}\upsilon\over 1-2\kappa}.
\end{equation}
\end{proposition}
\subsection{Error bounds for Lasso and Dantzig Selector under condition $\bH_{s,q}(\kappa)$}\label{secDanLas}
We are about to demonstrate that the feasibility of condition $\bH_{s,q}(\kappa)$ with $\kappa<\half$ implies some consequences for the performance of Lasso and Dantzig Selector when recovering $s$-sparse signals in $\|\cdot\|_p$ norms, $1\leq p\leq q$. This might look strange at the first glance, since neither Lasso nor Dantzig Selector use contrast matrices. The surprise, however, is eliminated by the following observation:
\begin{quote}
(!) {\sl Let $H$ satisfy $\bH_{s,q}(\kappa)$ and let $\widehat{\lambda}$ be the maximum of the Euclidean norms of columns in $H$. Then}
\begin{equation}\label{then}
\forall x\in \bR^n:\|x\|_{s,q}\leq \widehat{\lambda}s^{{1\over q}}\|Ax\|_2+\kappa s^{{1\over q}-1}\|x\|_1.
\end{equation}
\end{quote}
The fact that a condition like (\ref{then}) with $\kappa<1/2$ plays a crucial role in the performance analysis of Lasso and Dantzig Selector is neither surprising nor too novel. For example, the standard error bounds for the latter algorithms under the RIP assumption  are in fact based on the validity of
(\ref{then}) with $\widehat{\lambda}=O(1)$ for $q=2$ (see Lemma \ref{Dantzig}). Another example is given by the Restricted Eigenvalue \cite{tsyb1} and the Compatibility conditions \cite{sara2,sara1}. Specifically, the Restricted Eigenvalue condition $\RE(s,\rho,\varkappa)$ ($s$ is positive integer, $\rho>1$, $\varkappa>0$ states that
$$
\|x^s\|_2\leq {1\over\varkappa}\|Ax\|_2\hbox{\ whenever\ } \rho\|x^s\|_1\geq \|x-x^s\|_2,
$$
whence $\|x^s\|_1\leq {\sqrt{s}\over\varkappa}\|Ax\|_2$ whenever $(\rho+1)\|x^s\|_1\geq \|x\|_1$, so that
\begin{equation}\label{Tsybakov}
\forall x\in\bR^n: \|x\|_{s,1}\leq {s^{1/2}\over\varkappa}\|Ax\|_2+{1\over 1+\rho}\|x\|_1\an{;}{.}
\end{equation}
Further, the Compatibility condition of \cite{sara1} is nothing but (\ref{Tsybakov}) with $\rho=3$. We see that both Restricted Eigenvalue and Compatibility conditions imply
(\ref{then}) with $q=1$, $\widehat{\lambda}=(\varkappa\sqrt{s})^{-1}$ and certain $\kappa<1/2$.
\par
We are about to present a simple result on the performance of Lasso and Dantzig Selector algorithms in the case when $A$ satisfies the condition (\ref{then}). The result is as follows:
\begin{proposition}\label{LassoDantzig} Let $m\times n$ matrix $A=[a_1,...,a_n]$ satisfy (\ref{then}) with $\kappa<\half$ and some $q\in[1,\infty]$, and let $\beta=\max\limits_i\|a_i\|_2$. Let also the ``environment parameters'' $\sigma>0$, $\epsilon\in(0,1)$ be given, and let there be no nuisance: $\cU=\{0\}$.
\par
{\rm (i)} Consider the Dantzig Selector recovery
$$
\wh{x}_{\ds}(y)\in\Argmin_v\left\{\|v\|_1: \|A^T(Av-y)\|_\infty\leq \rho\right\},
$$
where
\begin{equation}\label{rhoDS}
\rho\geq \varrho:=\sigma\beta\sqrt{2\ln(n/\epsilon)}.
\end{equation}
Then
\begin{equation}\label{DantzigBound}
1\leq p\leq q\Rightarrow \Risk_p(\widehat{x}_{\ds}(\cdot)|\epsilon,\sigma,s,\upsilon)\leq {2(3s)^{1\over p}\over 1-2\kappa}
\left[{2s\widehat{\lambda}^2(\rho+\varrho)\over 1-2\kappa}+s^{-1}\upsilon\right].
\end{equation}
\par{\rm (ii)} Consider the Lasso recovery
$$
\wh{x}_{\las}(y)\in\Argmin_v\left\{\|v\|_1+\varkappa\|Av-y\|_2^2\right\},
$$
and let $\varkappa$ satisfy the relation
$$
2\kappa +2\varrho\varkappa<1,
$$
where $\varrho$ is given by (\ref{rhoDS}). Then
\begin{equation}\label{LassoBound}
1\leq p\leq q\Rightarrow \Risk_p(\widehat{x}_{\las}(\cdot)|\epsilon,\sigma,s,\upsilon)\leq {4s^{{1\over p}}\over 1-2\kappa-2\varrho\varkappa}
\left[{2s\widehat{\lambda}^2\over\varkappa}+s^{-1}\upsilon\right].
\end{equation}
In particular, with
\begin{equation}\label{varkappa}
\varkappa={1-2\kappa\over 4\varrho},
\end{equation}
one has
\begin{equation}\label{LassoBound1}
1\leq p\leq q\Rightarrow \Risk_p(\widehat{x}_{\las}(\cdot)|\epsilon,\sigma,s,\upsilon)\leq {8s^{{1\over p}}\over 1-2\kappa}
\left[{8s\varrho\widehat{\lambda}^2\over1-2\kappa}+s^{-1}\upsilon\right].
\end{equation}
\end{proposition}
\paragraph{Discussion.} Let us compare the error bounds given by Propositions \ref{arikwillnotlikeit}, \ref{LassoDantzig}. Assume that
there is no nuisance ($\cU=\{0\}$) and $A$ is such that the condition $H_{s,q}(\four)$ is satisfied by certain matrix $H$, the maximum of
Euclidean norms of the columns of $H$ being $\widehat{\lambda}$. Assuming that the penalized recovery uses $\theta=2$, and the regular recovery
uses $\widehat{\rho}=\nu(H)=\widehat{\lambda}\sigma\sqrt{2\ln(n/\epsilon)}$), the associated risk bounds
as given by Proposition \ref{arikwillnotlikeit} become
\begin{equation}\label{withcontrast}
\Risk_p(\widehat{x}(\cdot)|\e,\sigma,s,\upsilon)\le
O(1)s^{1\over p} \left[\widehat{\lambda}\sigma\sqrt{2\ln(n/\epsilon)}+s^{-1}\upsilon\right]\,\,1\leq p\leq q.
\end{equation}
Note that these bounds admit a transparent interpretation: in the range $1\leq p\leq q$ an $s$-sparse signal is recovered {\sl as if}
we were identifying correctly its support and estimating the entries with the uniform error $O(1)\widehat{\lambda}\sigma\sqrt{2\ln(n/\epsilon)}$.
\par
Now, as we have already explained, the existence of a matrix $H$ satisfying $\bH_{s,q}(\four)$ with columns in $H$ being of Euclidean lengths $\leq\widehat{\lambda}$
implies validity of (\ref{then}) with $\kappa=\four$. Assuming that in Dantzig Selector one uses $\rho=\varrho$, and that $\varkappa$ in Lasso
is chosen according to (\ref{varkappa}), the error bounds for Dantzig Selector and Lasso as given by Proposition \ref{LassoDantzig} become
\begin{equation}\label{withoutcontrast}
\Risk_p(\widehat{x}(\cdot)|\e,\sigma,s,\upsilon)\le
O(1)s^{1\over p} \left[[\beta\widehat{\lambda}]s\widehat{\lambda}\sigma\sqrt{2\ln(n/\epsilon)}+s^{-1}\upsilon\right]\,\,1\leq p\leq q.
\end{equation}
Observe that $\beta\widehat{\lambda}\geq O(1)$ (look what happens with (\ref{then}) when $x$ is the $i$-th basic orth). We see that the bounds
(\ref{withoutcontrast}) are worse than the bounds (\ref{withcontrast}), primarily due to the presence of the factor $s$ in the first bracketed term in
(\ref{withoutcontrast}). At this point it is unclear whether this drawback is an artifact caused
by poor analysis of the Dantzig Selector and Lasso algorithms or it indeed ``reflects reality.'' Some related numerical results presented in Section \ref{NumRes1} suggest that
the latter option could be the actual one.

Moreover, consider an example of the recovery problem with a $2\times 2$ matrix $A$ with unit columns and singular values $1$ and $\eps$. It can be easily seen that if $x$ is aligned with the second right singular vector of $A$ (corresponding to the singular value $\eps$) the error of the Dantzig Selector may be as large as $O(\eps^{-2}\sigma)$, while the error of ``$H$-conscious'' recovery  will be $O(\eps^{-1}\sigma)$  up to the logarithmic factor in $\e$ (indeed, choosing $H=A^{-1}$ results in $\lambda=\eps^{-1}$). This toy example suggests that the extra $\lambda$ factor in the bound \rf{withoutcontrast}, at least for Dantzig Selector, is not only due to our clumsy analysis.

This being said, it should be stressed that the comparison of regularized/penalised $\ell_1$
recoveries with Dantzig Selector and Lasso based solely on above the error bounds is somehow biased against Dantzig Selector and Lasso.
Indeed, in order for regular/penalized  $\ell_1$ recoveries to enjoy their ``good''  error bounds, we should specify the required contrast matrix, which is
not the case for Lasso and Dantzig Selector: the bounds (\ref{withoutcontrast}) require only existence of such a matrix\footnote{And even less than that, since
feasibility of $H_{s,q}(\kappa)$ is just a {\sl sufficient} condition for the validity of (\ref{then}), the condition which indeed underlies Proposition \ref{LassoDantzig}.}. Besides this, there is at least one case where error bounds for Dantzig Selector are as good as (\ref{withcontrast}), specifically,
the case when $A$ possesses, say, $\RIP(0.1,2s)$. Indeed, in this case, by Lemma \ref{Dantzig}, the matrix $H=O(1)A$ satisfies $\bH_{s,2}(\four)$, meaning that
Dantzig Selector with properly chosen $\rho$ is nothing but the regular recovery with contrast matrix $H$ and as such obeys the bounds (\ref{withcontrast}) with $q=2$.
\par
It is time to point out that the above discussion is somehow scholastic: when $q<\infty$ and $s$ is nontrivial, we do not know how to verify
efficiently the fact that the condition $\bH_{s,q}(\kappa)$ is satisfied by a given $H$, not speaking about efficient synthesis of $H$ satisfying this condition. One should
not think that these tractability issues concern only our algorithms which need a good contrast matrix. In fact, {\sl all} conditions which allow
to validate Dantzig Selector and Lasso beyond the scope of the ``fully tractable'' condition $H_{s,\infty}(\kappa)$ are, to the best of our knowledge,
unverifiable -- they cannot be checked efficiently, and thus we never can be sure that Lasso and Dantzig Selector (or any other known computationally efficient technique
for sparse recovery) indeed work well for a {\em given sensing matrix}. As we have seen in Section \ref{matrixh}, the situation improves dramatically when passing from unverifiable
conditions $\bH_{s,q}(\kappa)$, $q<\infty$, to the efficiently verifiable condition $\bH_{s,\infty}(\kappa)$,  although in a severely restricted
range of values of $s$.

\section{Numerical examples}\label{sectnumres}\label{Numresults}
We present here a small simulation study.
\subsection{Regular/penalized recovery vs. Lasso: no-nuisance case}\label{NumRes1}
To illustrate the discussion in Section \ref{secDanLas}, we compare numerical performance of Lasso and penalized recovery in the observation model \rf{nun} without nuisance:
\[
y=Ax+\sigma\xi, \;\;\;\xi\sim N(0,I_m),
\]
where $\sigma>0$ is known. The sensing matrix  $A$ is specified by selecting at random $m=120$ rows of the $128\times 128$ Hadamard matrix\footnote{The $k$-th Hadamard matrix $H^k$ is given by the recurrence $H^0=1$, $H^{p+1}=[H^p;H^p;H^p,-H^p]$. It is a $2^k\times 2^k$ matrix with orthogonal rows and all entries equal to $\pm1$.}, and ``suppressing'' the first of the selected rows by multiplying it by 1.e-3. The resulting $120\times 128$ sensing matrix has orthogonal rows; 119 of its 120 singular values are equal to $8\sqrt{2}$, and the remaining singular value is $0.008\sqrt{2}$.
 \par
We have processed $A$ as explained in Section \ref{matrixh} (a reader is referred to this section for the description of entities involved).\footnote{It is worth to mention that when $A$ is comprised of (perhaps, scaled) rows of an Hadamard matrix (and in fact, of scaled rows of any other Fourier transform matrix associated with a finite Abelian group) the synthesis described in Section \ref{matrixh} simplifies dramatically due to the fact that all problems $(P_i^\gamma)$ turn out to be equivalent to each other, and their optimal solutions are obtained from each other by simple linear transformations. As a result, we can work with a {\sl single} problem $(P_1^\gamma)$ instead of working with $n$ of them.} We started with computing $\gamma_*$, which turned out to be
    0.0287, meaning that the level of $s$-goodness of $A$ is at least 17. In our experiment, we aimed at recovering signals with at most $s=10$ nonzero entries and with no nuisance ($\cU=\{0\}$). The synthesis of the corresponding ``optimal'' contrast matrix $H=H_*$ as outlined in Section \ref{matrixh} results in $\bar{\gamma}=0.294$, $\omega_*(\bar{\gamma})=0.0899\sqrt{2\ln (n/\epsilon)}$. Note that we are in the case of $\cU=\{0\}$, and in this case the optimal $H$ is independent of the values of $\sigma$ and $\epsilon$.
\par
We compare the penalized $\ell_1$-recovery with the contrast matrix $H_*$ and $\theta=2$ with the Lasso recovery on randomly generated signals $x$ with 10 nonzero entries. We consider two choices of the penalty $\varkappa$ in Lasso: the ``theoretically optimal'' choice (\ref{varkappa}) and the ``ideal'' choice, where we scanned the ``fine grid'' $(1.05)^k$, $k=0,\pm1,\pm2,...$ of values of $\varkappa$ and selected the value for which the Lasso recovery was at the smallest $\|\cdot\|_1$-distance from the true signal. The confidence parameter $\epsilon$ in (\ref{varkappa}) was set to $0.01$.

The results of a typical experiment are presented in Table \ref{tabHadam}. We see that as compared to the penalized $\ell_1$ recovery,
the accuracy of Lasso with the theoretically optimal choice of the penalty is nearly 10 times worse. With the ``ideal'' (unrealistic!) choice of penalty, Lasso is never better than the penalized $\ell_1$ recovery, and for the smallest value of $\sigma$ is nearly 4 times worse than the latter routine.
\begin{table}
\begin{center}
{\footnotesize
\begin{tabular}{|l|c|c|c|c|c|}
\cline{4-6}
\multicolumn{3}{c}{}&\multicolumn{3}{|c|}{$\|\widehat{x}-x\|_p$}\\
\hline
Recovery&$\sigma$&$\varkappa$&$p=1$&$p=2$&$p=\infty$\\
\hline\hline
Penalized&&N/A&2.1e-4&6.5e-5&3.8e-5\\
\hline
Lasso&1.e-4&3.74e-3$^*$&2.1e-4&5.2e-4&3.9e-5\\
\hline
Lasso&&4.01e-2$^\dag$&1.6e-3&5.2e-4&2.0e-4\\
\hline\hline
Penalized&&N/A&2.2e-5&6.0e-6&2.7e-6\\
\hline
Lasso&1.e-5&4.78e-4$^*$&3.1e-5&8.1e-6&3.4e-6\\
\hline
Lasso&&4.01e-3$^\dag$&1.8e-4&5.8e-5&2.1e-5\\
\hline\hline
Penalized&&N/A&2.1e-6&6.2e-7&2.5e-7\\
\hline
Lasso&1.e-6&1.10e-4$^*$&8.8e-6&1.6e-6&5.9e-7\\
\hline
Lasso&&4.01e-4$^\dag$&1.8e-5&5.4e-6&1.9e-6\\
\hline\hline
\end{tabular}}
\caption{\label{tabHadam}. Lasso vs. penalized $\ell_1$ recovery. Choice of $\varkappa$: $^*$ -- ``ideal'' choice; $^\dag$ -- theoretical choice.}
\end{center}
\end{table}
\subsection{The nuisance case}\label{NumRes2}
In the second experiment we study the behavior of recovery procedures in the situation when
an ``input nuisance'' is present:
\[
y=A(x+v)+\sigma\xi,
\]
 where  $x\in \bR^n$ is an unknown sparse signal, $v\in {\cal V}$ with known ${\cal V}\subset \bR^n$, $\sigma$ is known and $\xi\in \bR^m$ is standard normal $\xi\sim N(0,I_m)$; in terms of \rf{obs}, $u=Av$ and $\cU=A\cV$. We compare the performance of the regular and penalized recoveries to that of the Lasso and Dantzig Selector algorithms. To handle the nuisance, the latter methods were modified as follows: instead of the standard Lasso estimator we use the estimator
\[
\wh{x}_{\las}(y)\in \Argmin_{x\in \bR^n}\min_{v\in\cV}\left\{
\|x\|_1+\varkappa \|A(x+v)-y\|_2^2\right\},
\]
where the penalization coefficient $\varkappa$ is chosen according to \cite[Theorem 4.1]{tsyb1}; in turn, the Dantzig Selector is substituted by
\[
\wh{x}_\ds(y)\in \Argmin_{x\in \bR^m}\min\limits_{v\in\cV}\left\{
\|x\|_1\,:|[A^T(A(x+v)-y)]_i|\le \varrho_i,\;\;i=1,...,m\right\}
\]
with $\varrho_i=\sigma\sqrt{2\ln (n/\epsilon)}\|A_i\|_2$, where $A_i$ are the columns of $A$ and $\e$ is given (in what follows $\e=0.01$).
\par
We present below the simulation results for two setups with $n=256$:
\begin{enumerate}
\item {\em Gaussian setup}: a ${161\times 256}$ sensing matrix $A_{\rm Gauss}$ with independent $N(0,1)$ entries is generated, then its columns are normalized. The nuisance set $\cV=\cV(L)\subset \bR^{256}$ is as follows:
\[
\cV(L)=\{v\in \bR^{256},\;
|v_{i+1}-2v_i+v_{i-1}|\le L,\;\mbox{for $i=2,...,255$},\;v_2=v_1=0\},
\]
where $L$ is a known parameter; in other words, we observe the sum of a sparse signal and ``smooth background.''
\item {\em Convolution setup}: a $240\times 256$ sensing matrix  $A_{\rm conv}$ is constructed as follows: consider a signal $x$ ``living'' on ${\mathbf{Z}}^2$ and supported on the $16\times 16$ grid $\Gamma=\{(i,j)\in{\mathbf{Z}}^2: \;0\leq i,j\leq 15\}$. We subject such a signal to discrete time convolution with a kernel supported on the set $\{(i,j)\in{\mathbf{Z}}^2: -7\leq i,j\leq 7\}$, and then restrict the result on the $16\times 15$ grid  $\Gamma_+=\{(i,j)\in\Gamma: 1\leq j\leq 15\}$. This way we obtain a linear mapping $x\mapsto A_{\rm  conv}x: \;\bbr^{256}\to\bbr^{240}$. The nuisance set $\cV=\cV(L)\subset \bR^{256}$ is composed of zero-mean signals $u$ on $\Gamma$ which satisfy
    \[
    |[D^2u]_{i,j}|\le L,
    \]
    where $D$ is the discrete (periodic) homogeneous Laplace operator:
\[
[Du]_{i,j}={1\over 4}\left(u_{i,\overline{j-1}}
+u_{\overline{i-1},j}
+u_{i,\overline{j+1}}
+u_{\overline{i+1},j}-4u_{i,j}\right),\;\;i,j=1,...,16,
\]
with $\overline{i}=i\mod 16,\;\overline{j}=j\mod 16$.
\end{enumerate}
In the simulations we acted as follows: given the sensing matrix $A$, the nuisance set $\cU=A\cV$ and the values of $s$ and $\sigma$, we compute the contrast matrix $H$ by choosing a ``reasonable'' value $\bar{\gamma}>\gamma_*$ of $\gamma$ and specifying $H$ as the matrix satisfying $\bH_{s,\infty}(s\bar{\gamma})$ and such that $\nu(H)=\omega_*(\bar{\gamma})$, see Section \ref{matrixh}}. Then $N$ samples of {\em random} signal $x$, random nuisance $v\in \cV$ and random perturbation $\xi$ were generated, and the corresponding observations were processed by every one of the algorithms we are comparing\footnote{Randomness of the sparse signal $x$ is important. Using the techniques of \cite{JNCS}, one can verify that in the convolution setup there are signals with only 3 non-vanishing components which cannot be recovered by $\ell_1$ minimization even in the noiseless case $\cV=\{0\}$, $\sigma=0$. In other words, the $s$-goodness characteristic of the corresponding matrix $A$ is equal to $2$.}. The plots below present the average, over these $N=100$ experiments, $\ell_\infty$ and $\ell_1$ recovery errors.
All recovery procedures were using {\tt Mosek} optimization software \cite{mosek}.
\par
\begin{figure}[ht]
  \centering
\begin{tabular}{cc}
    \includegraphics[width=0.45\textwidth]{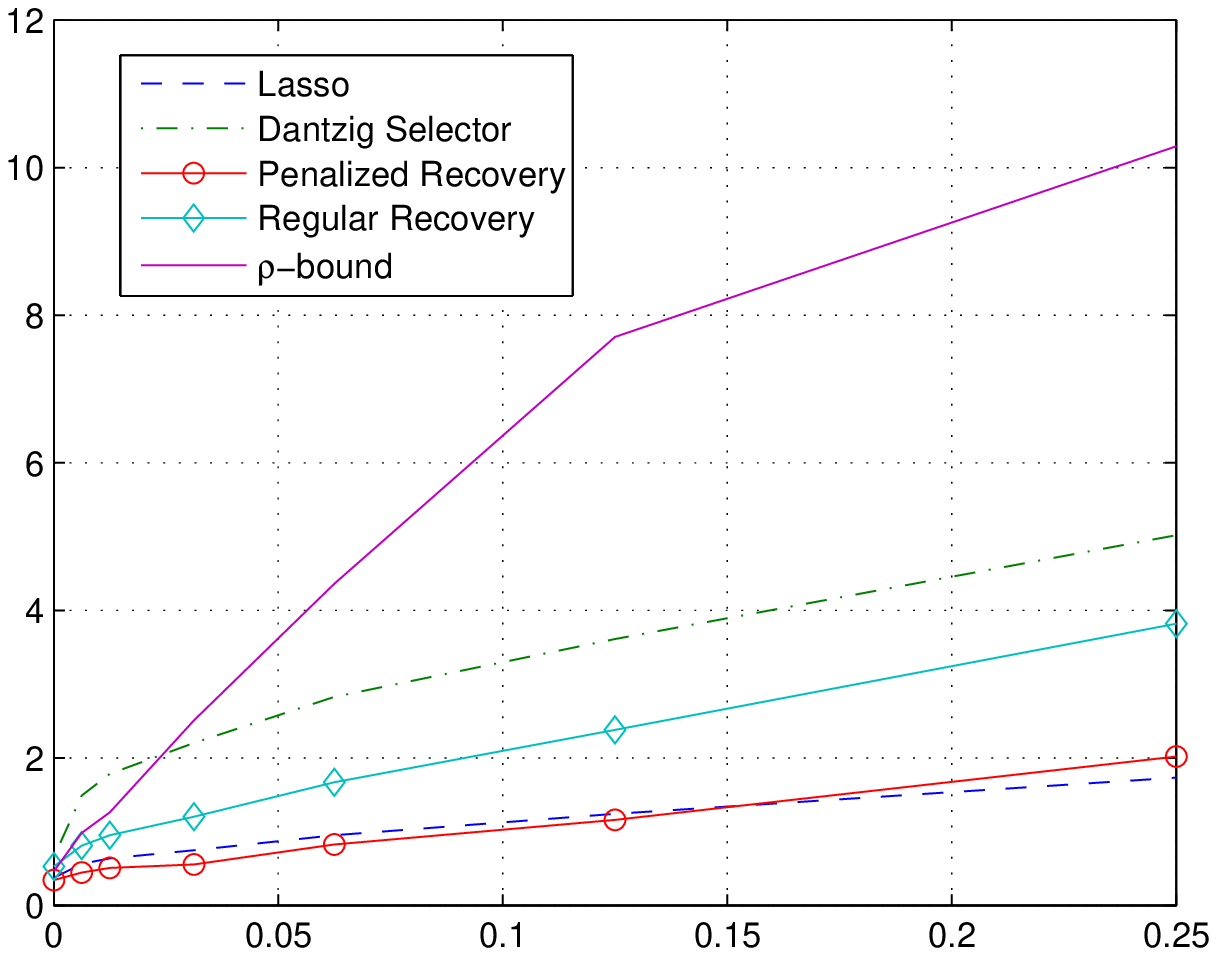} &
    \includegraphics[width=0.45\textwidth]{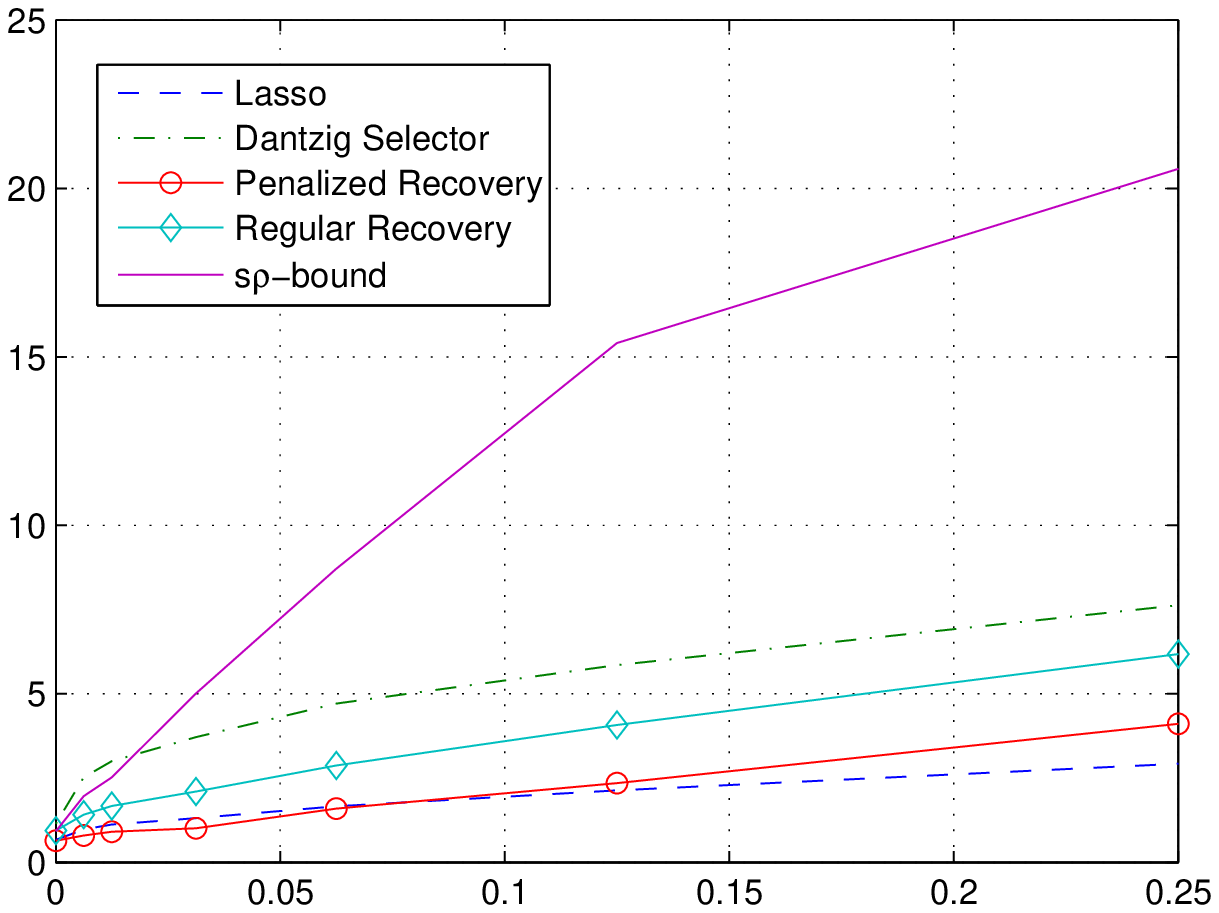} \\
        $\ell_\infty$-error & $\ell_1$-error
 \end{tabular}
\caption{\label{fig:gauss1} Mean recovery error as a function of the nuisance magnitude $L$. Gaussian setup parameters: $\sigma=0.1$, $s=2$, $\mu=0.1$, $\|x\|_1=10$.}
\end{figure}
\begin{figure}[ht]
  \centering
\begin{tabular}{cc}
    \includegraphics[width=0.45\textwidth]{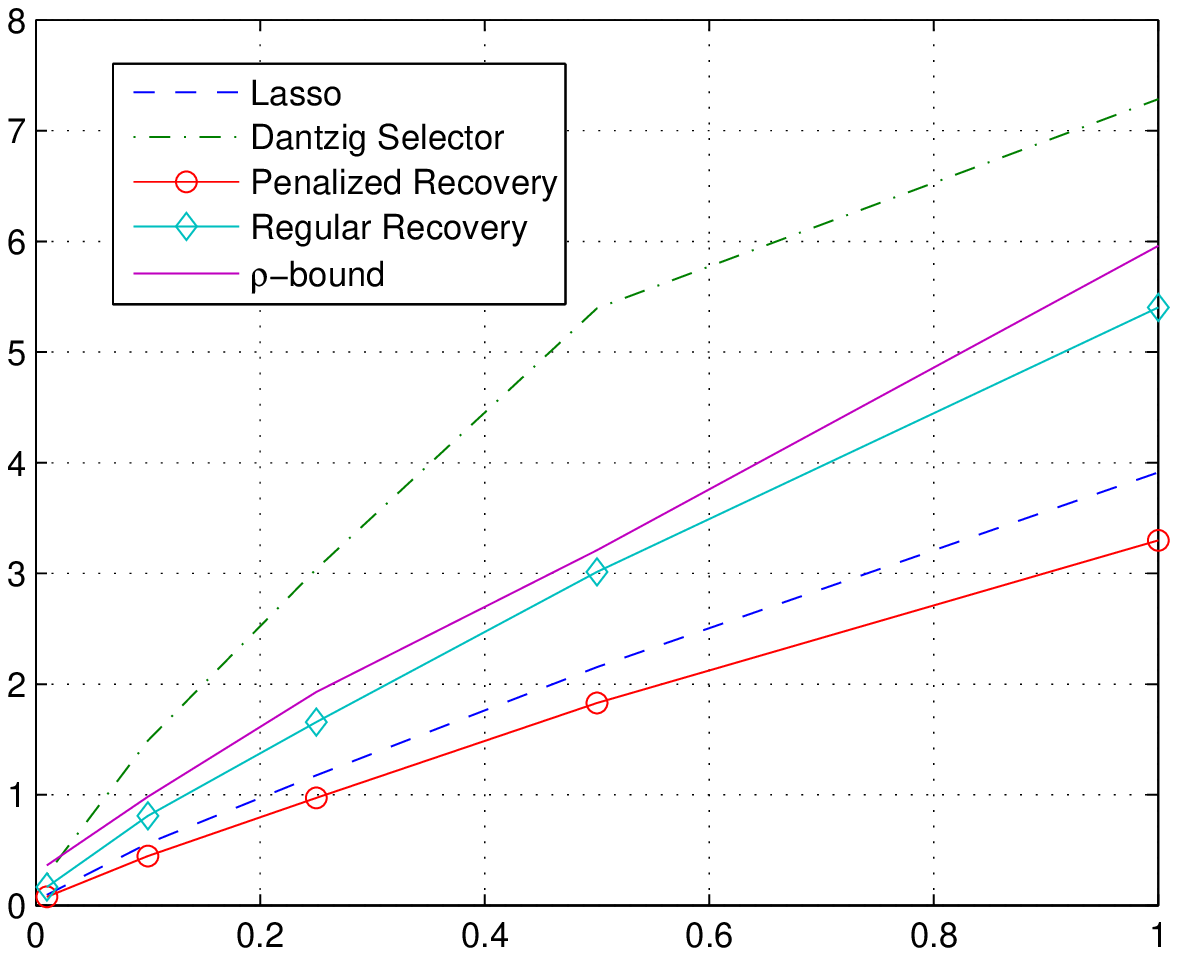} &
    \includegraphics[width=0.45\textwidth]{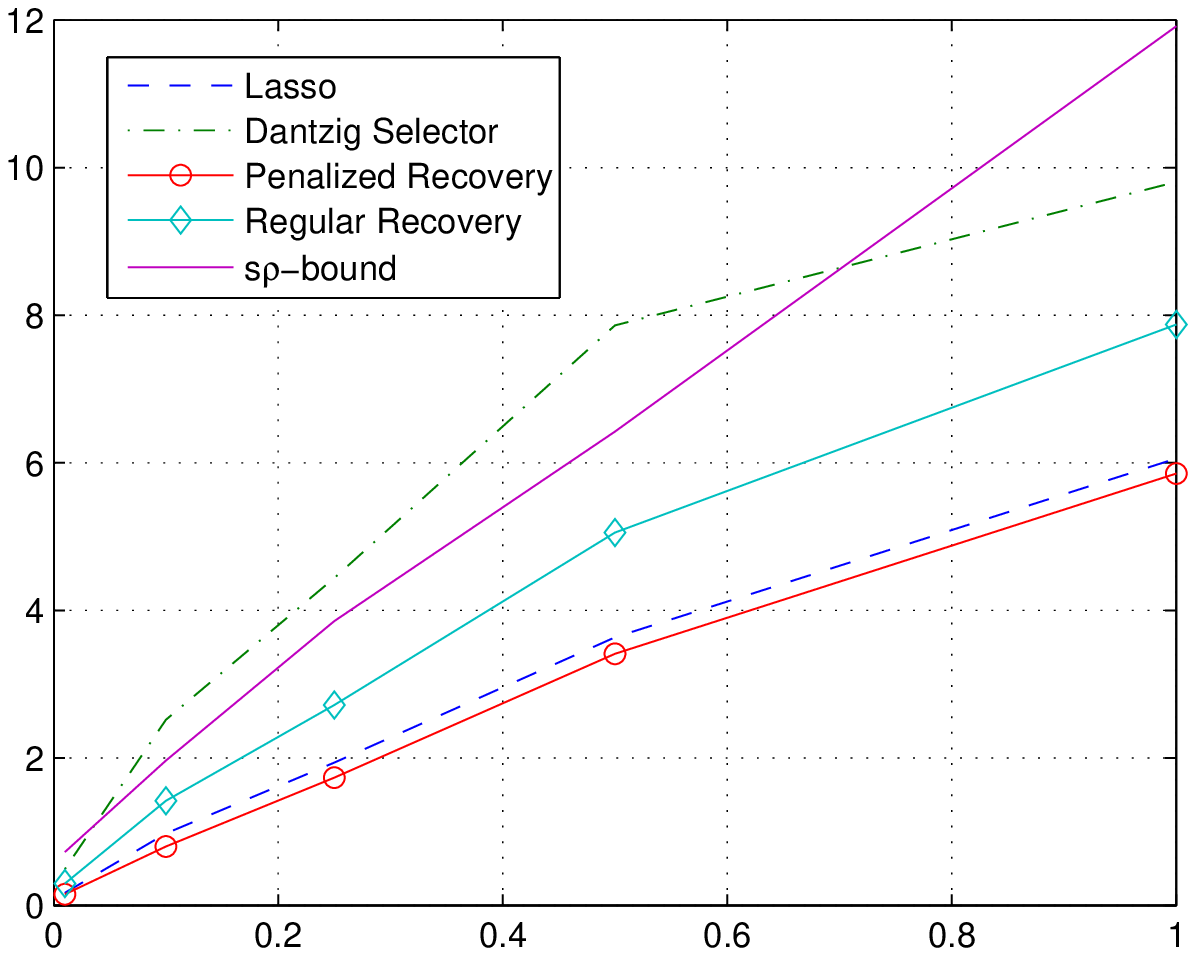} \\
        $\ell_\infty$-error & $\ell_1$-error
 \end{tabular}
  \caption{Mean recovery error as a function of the noise StD $\sigma$. Gaussian setup parameters: $L=0.01$, $s=2$, $\mu=0.1$, $\|x\|_1=10$. \label{fig:gauss2}
  }
\end{figure}
\begin{figure}[ht]
  \centering
\begin{tabular}{cc}
    \includegraphics[width=0.45\textwidth]{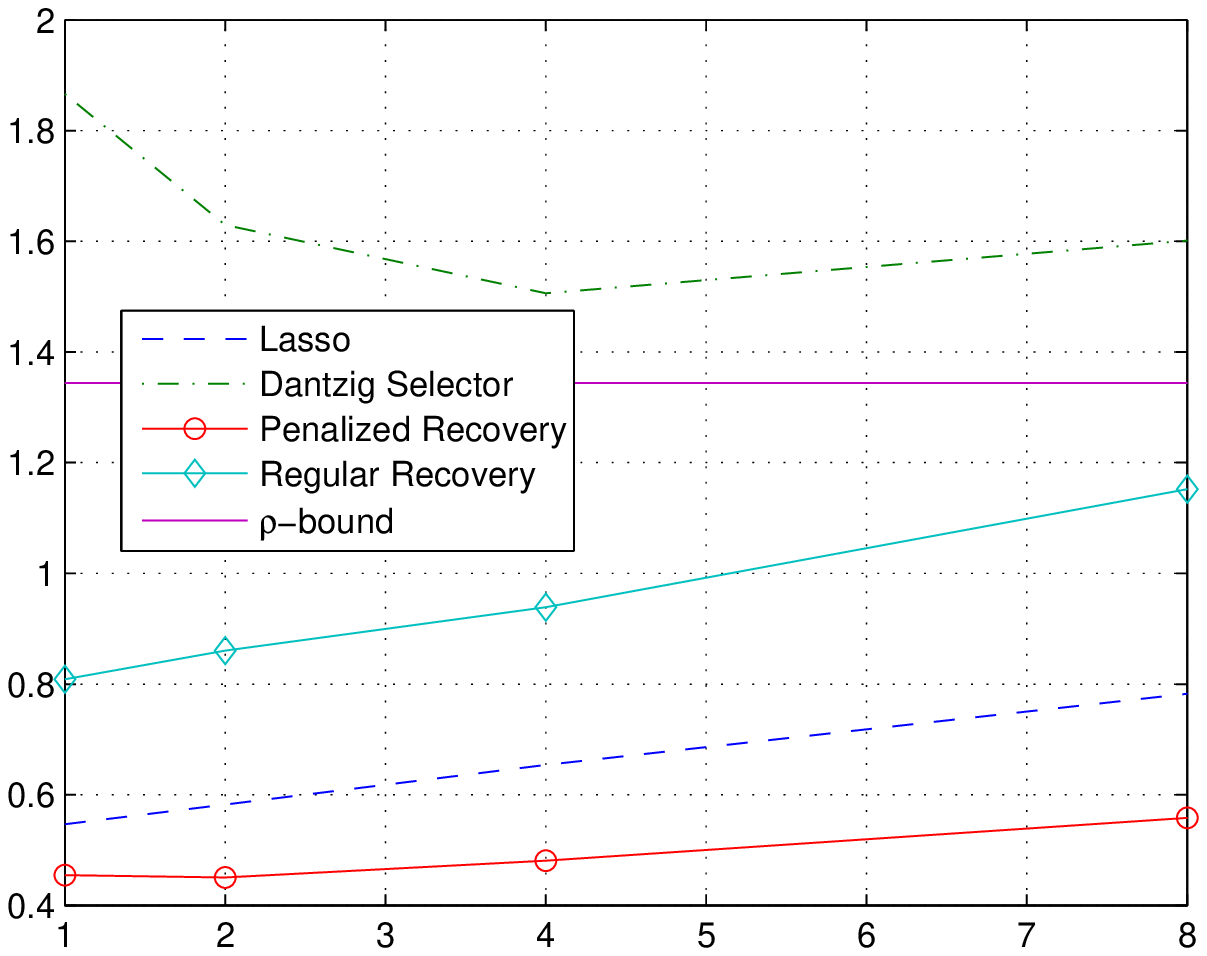} &
    \includegraphics[width=0.45\textwidth]{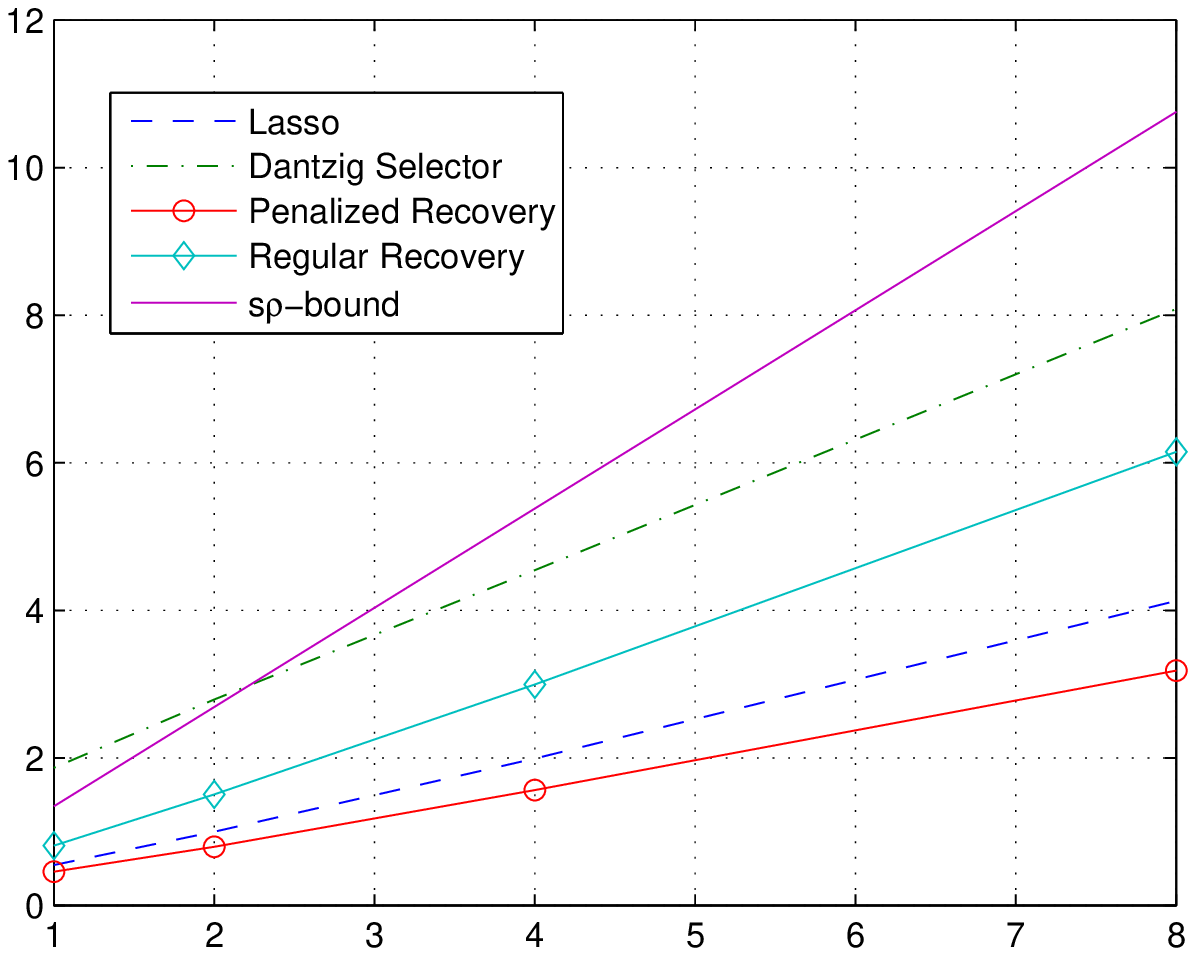} \\
        $\ell_\infty$-error & $\ell_1$-error
 \end{tabular}
  \caption{Mean recovery error as a function of the number $s$ of nonzero entries in the signal. Gaussian setup parameters: $L=0.01$, $\sigma=0.1$, $\bar{\gamma}=0.1$, $\|x\|_1=5s$. \label{fig:gauss3}
  }
\end{figure}

We start with Gaussian setup in which the signal $x$ has $s=2$ non-vanishing components, randomly drawn, with $\|x\|_1=10$. For the penalized and the regular recovery algorithms the contrast matrix $H$ was computed using $\bar{\gamma}=0.1$.
On Figure \ref{fig:gauss1} we plot the average recovery error as a function of the value of the parameter $L$ of the nuisance set $\cV$, for fixed $\sigma=0.1$, and on Figure \ref{fig:gauss2} --- as a function of
$\sigma$ for fixed $L=0.01$.
In the next experiment we fix the ``environmental parameters'' $\sigma$, $L$ and vary the number $s$ of nonzero entries in the signal $x$ (of norm $\|x\|_1=5s$). On Figure \ref{fig:gauss3} we present the recovery error as a function of $s$.
\par
\begin{figure}[ht]
  \centering
\begin{tabular}{cc}
    \includegraphics[width=0.45\textwidth]{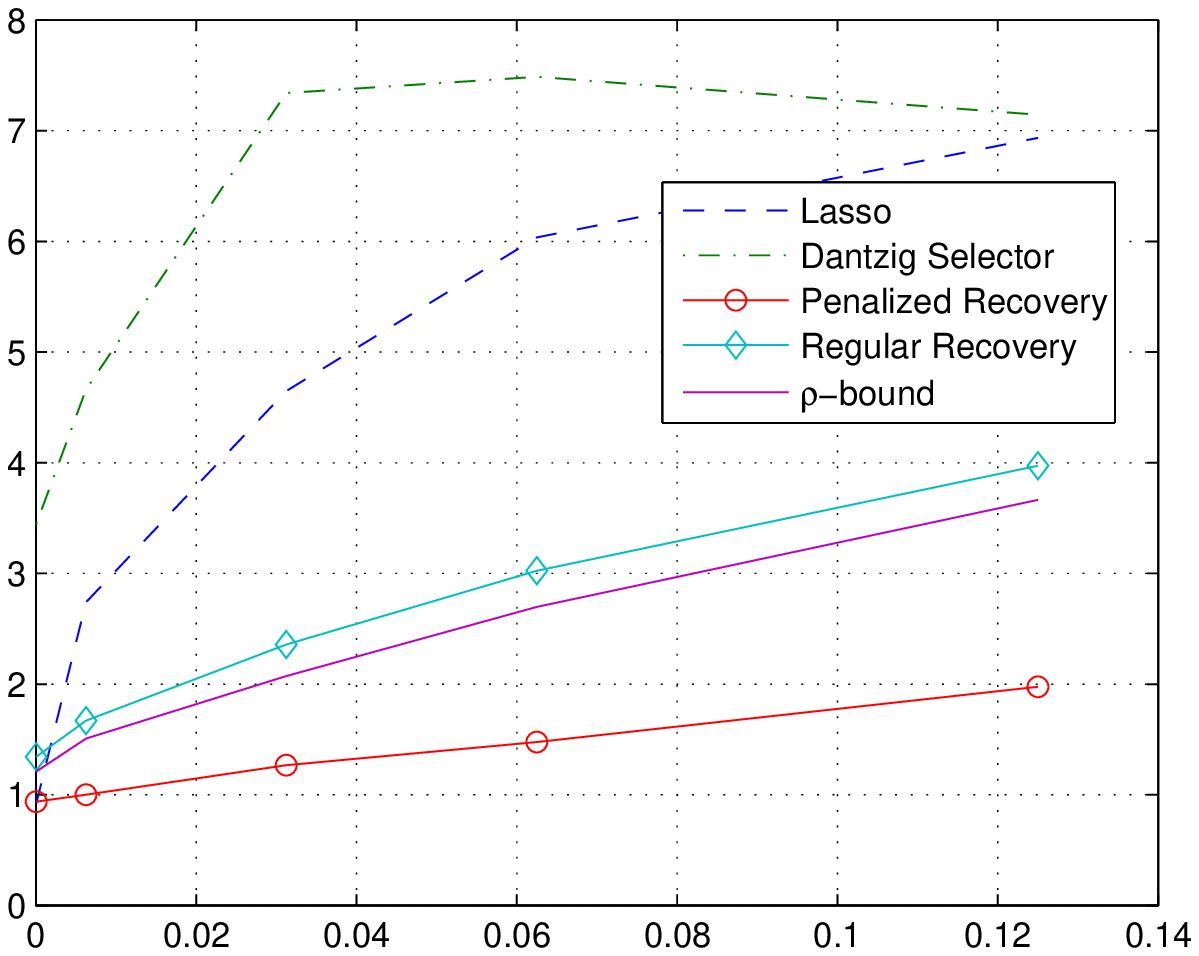} &
    \includegraphics[width=0.45\textwidth]{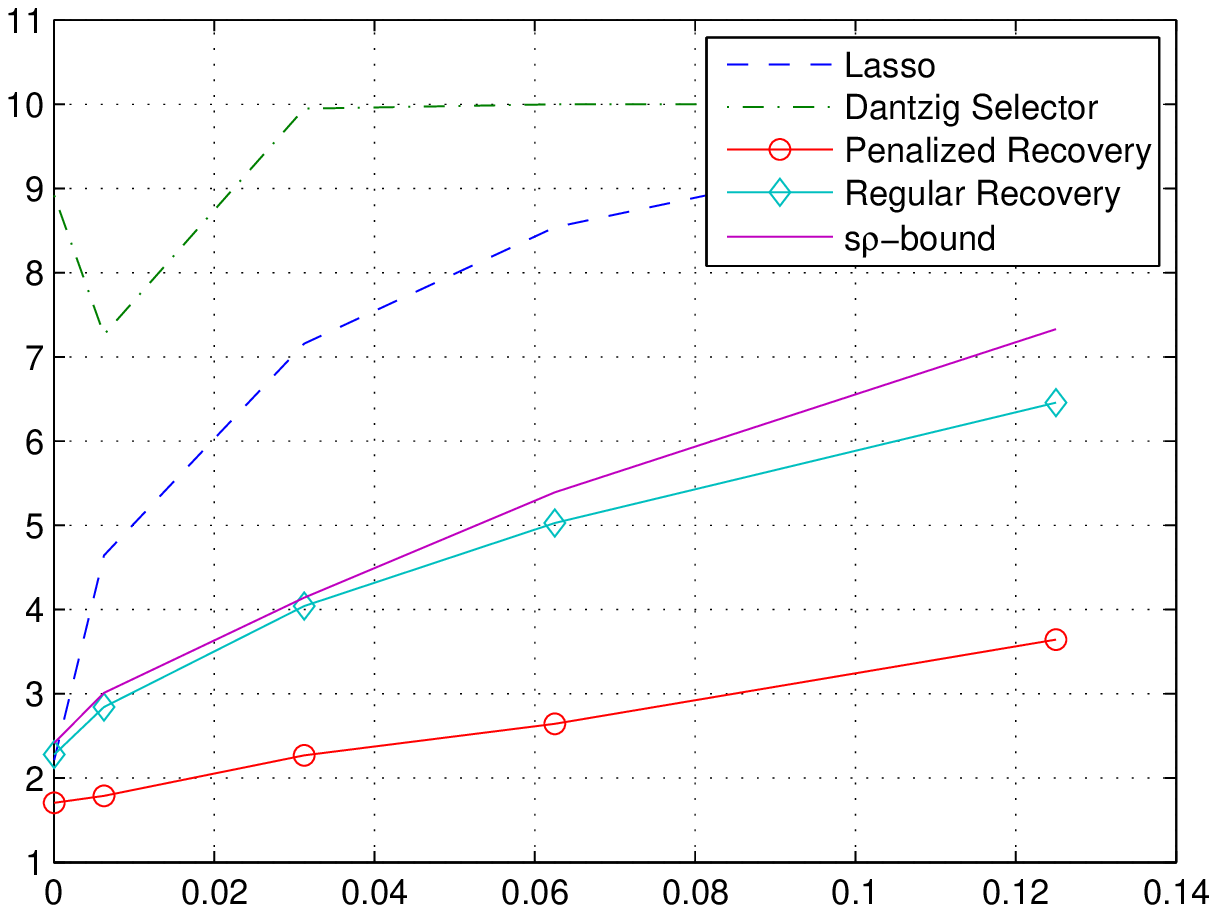} \\
        $\ell_\infty$-error & $\ell_1$-error
 \end{tabular}
\caption{Mean recovery error as a function of the nuisance magnitude $L$. Convolution setup parameters: $\sigma=0.1$, $s=2$, $\bar{\gamma}=0.2$, $\|x\|_1=10$. \label{fig:conv1}}
\end{figure}
\begin{figure}[ht]
  \centering
\begin{tabular}{cc}
    \includegraphics[width=0.45\textwidth]{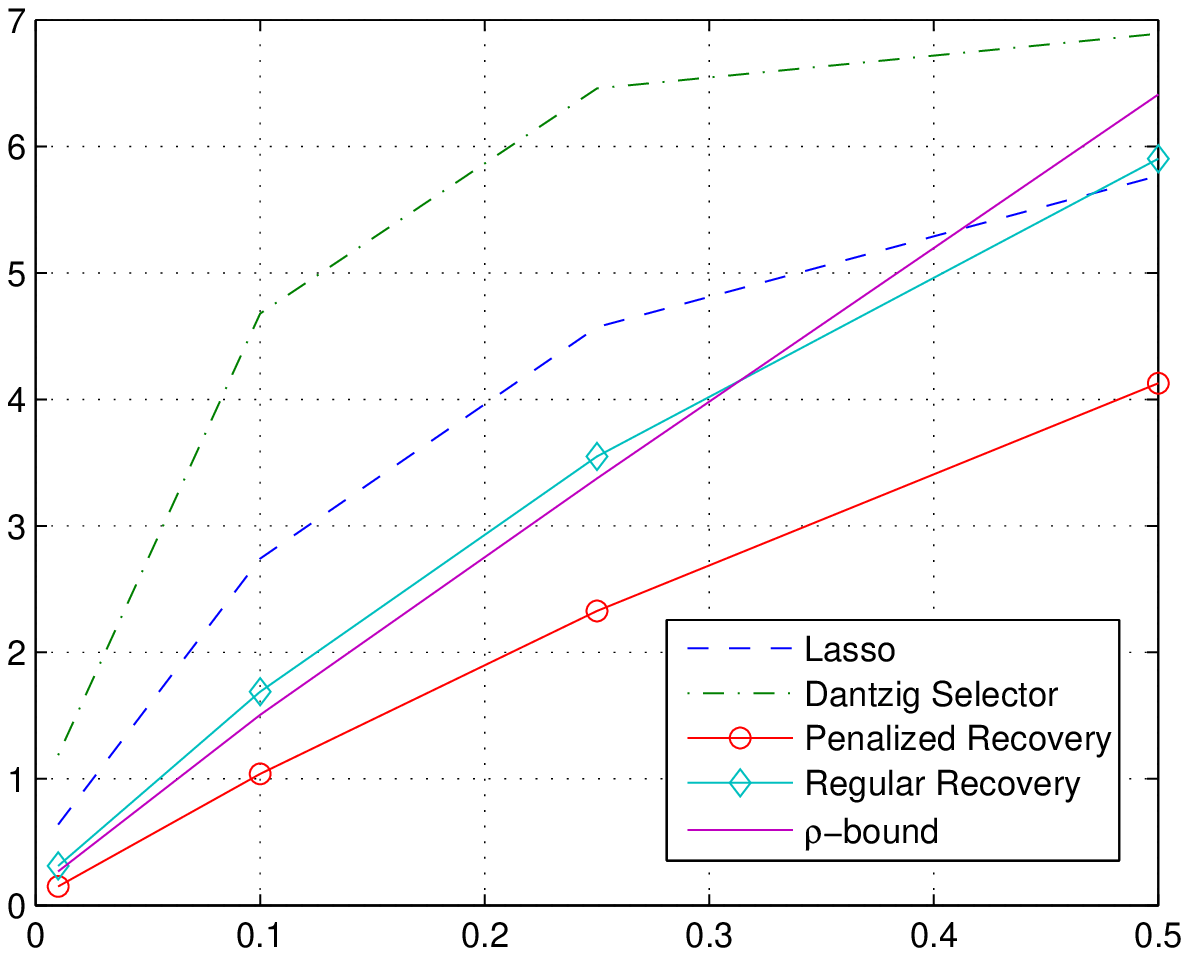} &
    \includegraphics[width=0.45\textwidth]{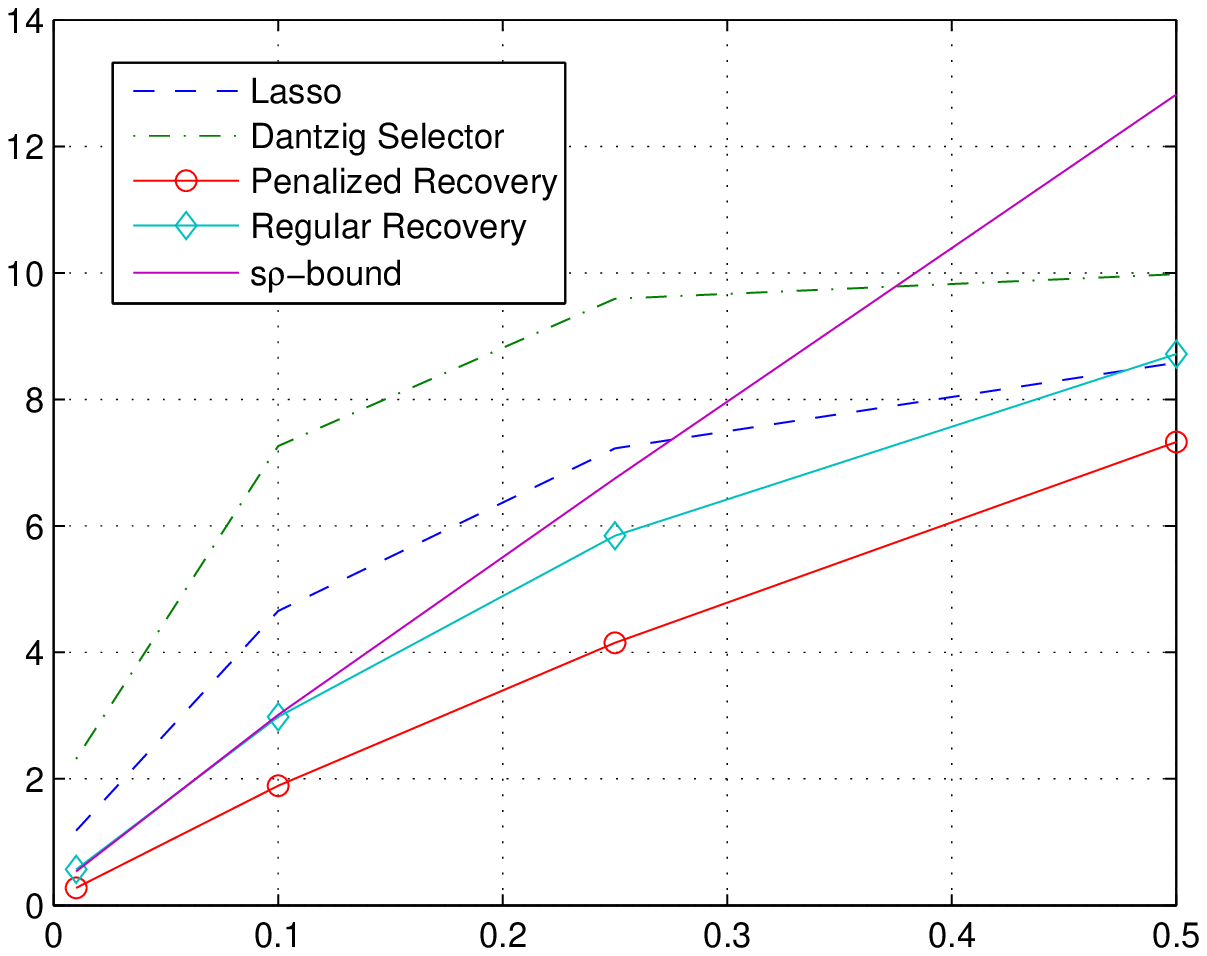} \\
        $\ell_\infty$-error & $\ell_1$-error
 \end{tabular}
  \caption{Mean recovery error as a function of the noise StD $\sigma$. Convolution setup parameters: $L=0.01$, $s=2$, $\bar{\gamma}=0.2$, $\|x\|_1=10$. \label{fig:conv2}
  }
\end{figure}
\begin{figure}[ht]
  \centering
\begin{tabular}{cc}
    \includegraphics[width=0.45\textwidth]{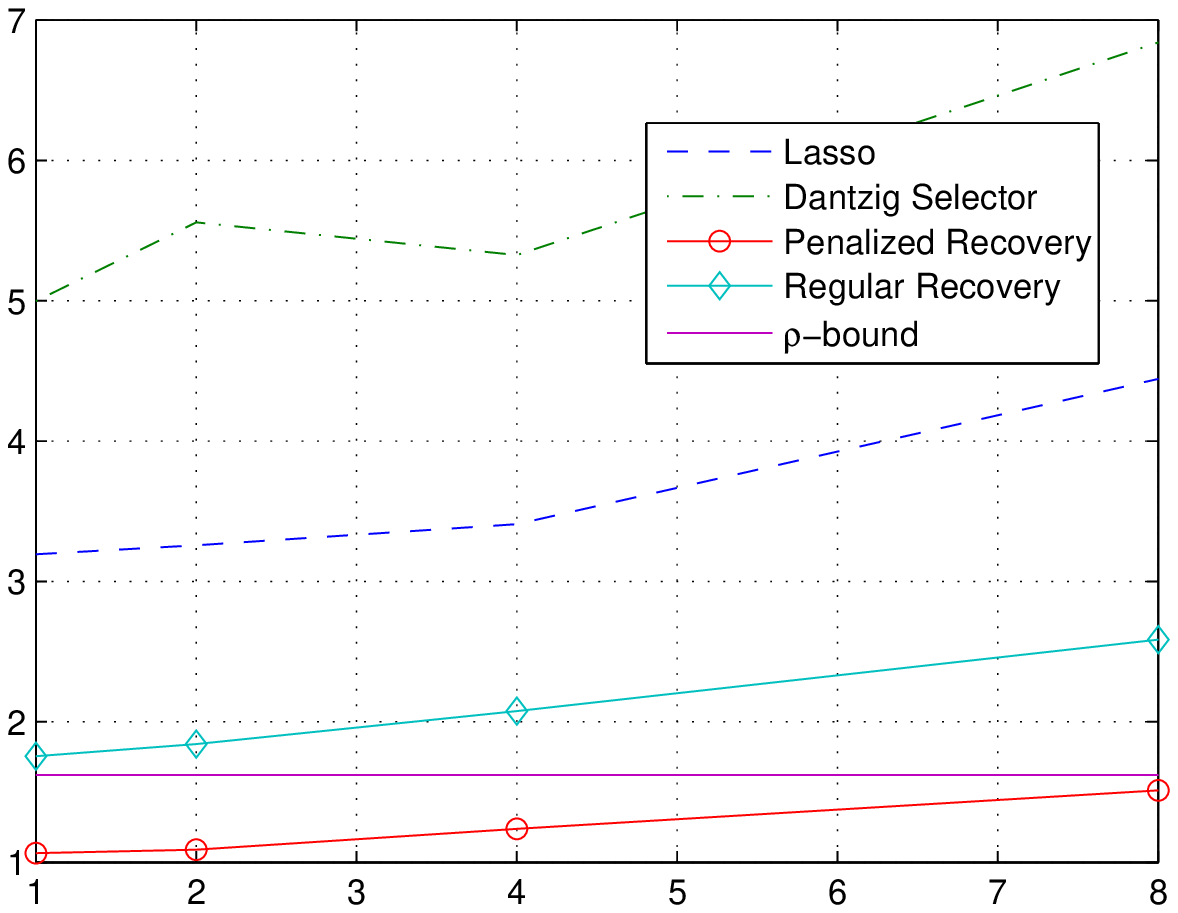} &
    \includegraphics[width=0.45\textwidth]{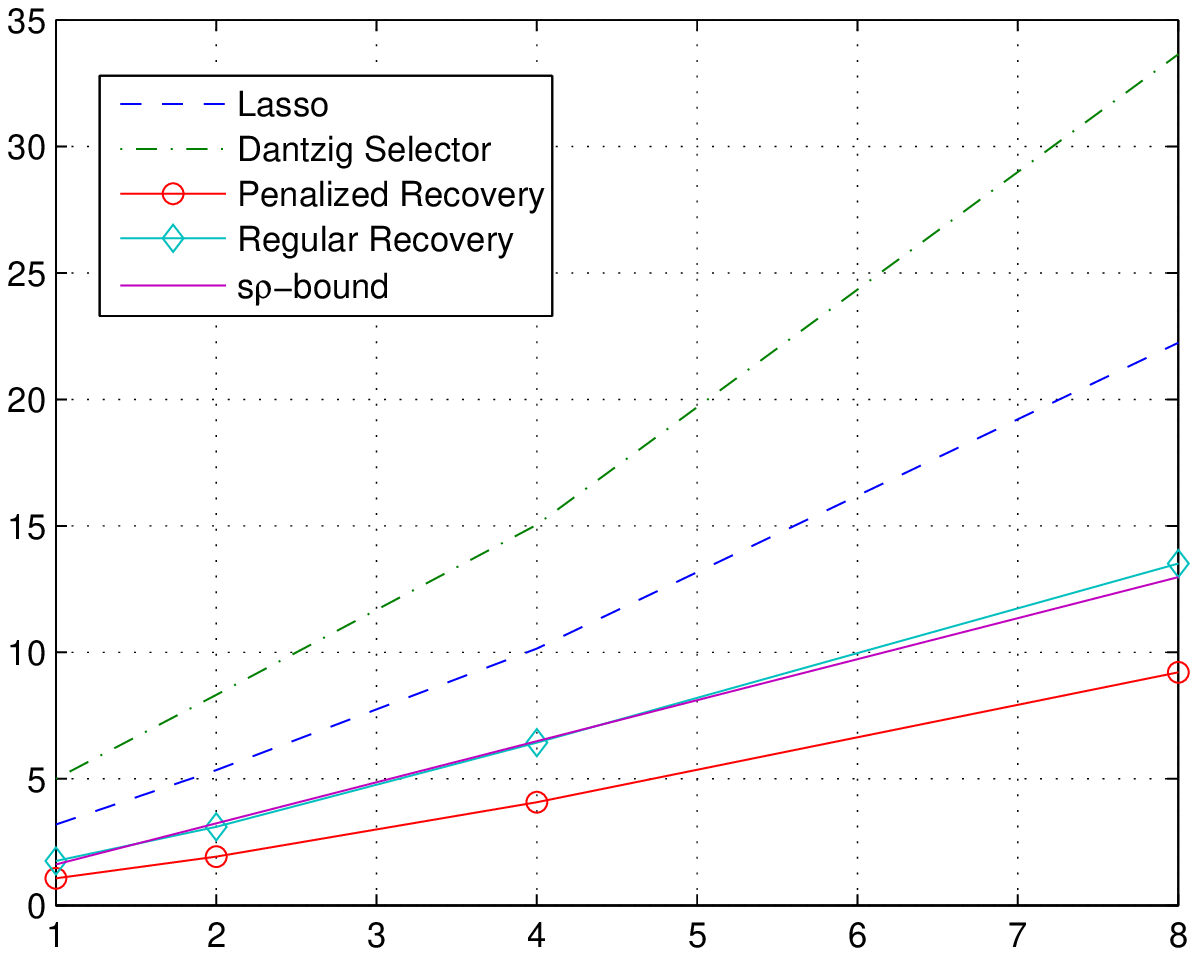} \\
        $\ell_\infty$-error & $\ell_1$-error
 \end{tabular}
  \caption{Mean recovery error as a function of the number $s$ of nonzero entries in the signal. Convolution setup parameters: $L=0.01$, $\sigma=0.1$, $\bar{\gamma}=0.2$, $\|x\|_1=5s$. \label{fig:conv3}
  }
\end{figure}

 We run the same simulations in the convolution setup. The contrast matrix $H$ for the penalized and the regular recoveries  is computed using $\bar{\gamma}=0.2$. On Figure \ref{fig:conv1} we plot the average recovery error as a function of the ``size'' $L$ of the nuisance set $\cV$ for fixed $\sigma=0.1$, on Figure \ref{fig:conv2} --- as a function of
$\sigma$ for fixed $L=0.01$, and on Figure \ref{fig:conv3} --- as a function of $s$.

\par
We observe quite different behavior of the recovery procedures in our two setups. In the Gaussian setup the nuisance signal $v\in \cV$ does not mask the true signal $x$, and the performance of the Lasso and Dantzig Selector is quite good in this case. The situation changes dramatically in the convolution setup, where the performance of the Lasso  and Dantzig Selector degrades rapidly when the parameter $L$ of the nuisance set increases.\footnote{The error plot for these estimators on Figure \ref{fig:conv1} flatters for higher values of $L$ simply because they always underestimate the signal, and the error of recovery is always less than the corresponding norm of the signal.} The conclusion suggested by the outlined numerical results is that the penalized $\ell_1$ recovery, while sometimes losing slightly to Lasso, in some of the experiments outperforms significantly all other algorithms we are comparing.
\begin{figure}[ht]
  \centering
\begin{tabular}{cc}
    \includegraphics[width=0.45\textwidth]{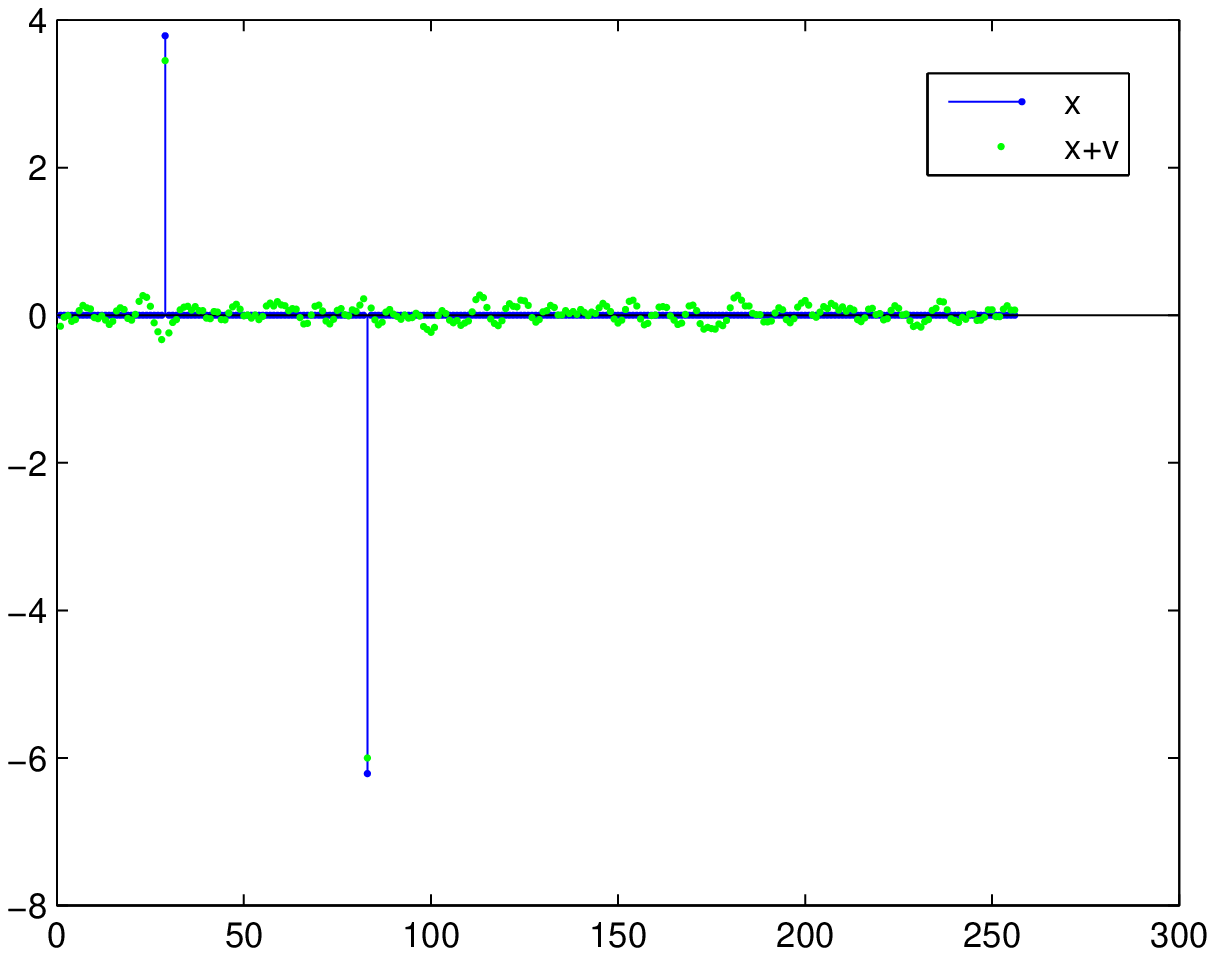} &
    \includegraphics[width=0.45\textwidth]{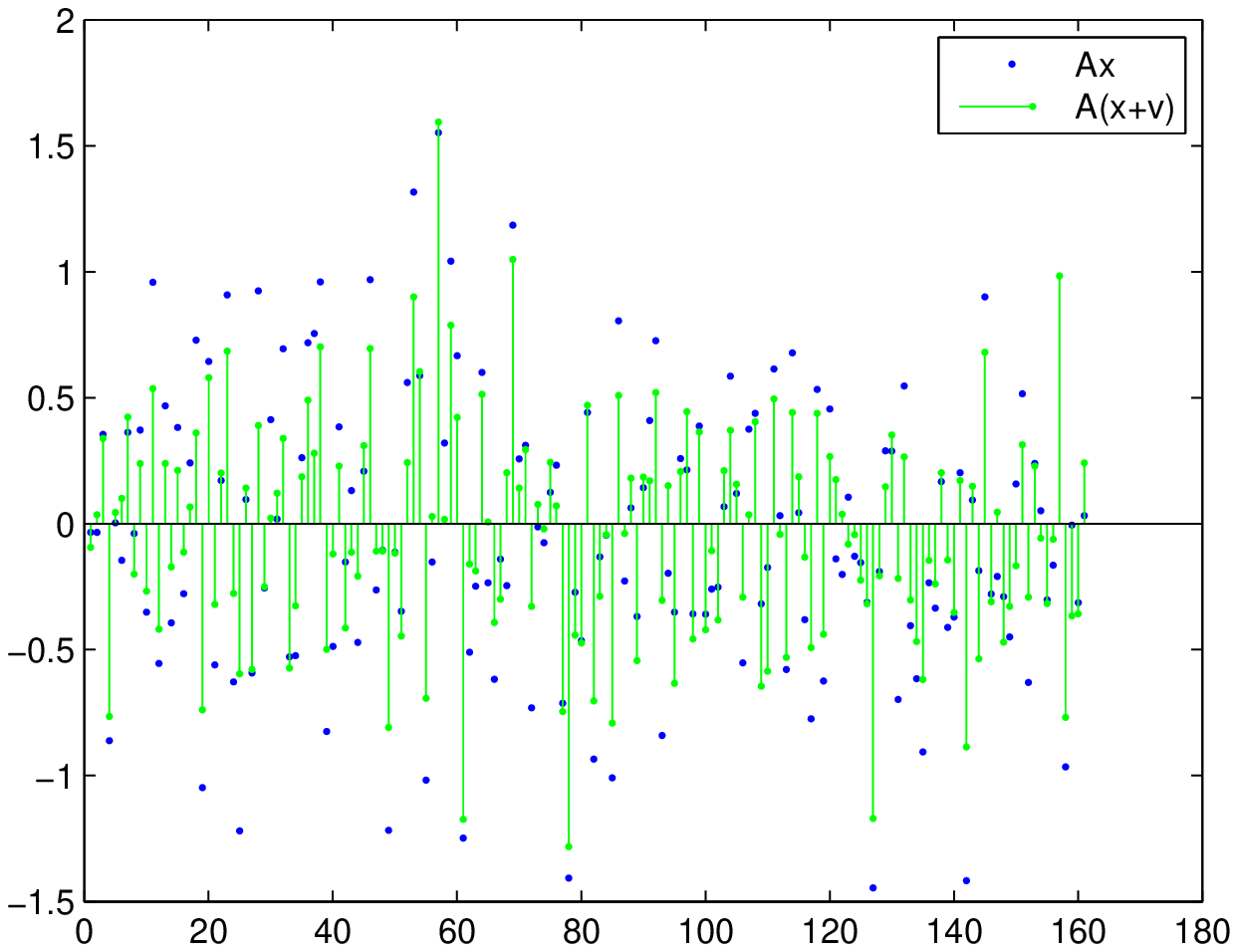} \\
        signal $x$ and nuisance $v$  & image $Ax$ and $A(x+v)$\\
 \end{tabular}
  \caption{A typical signal/worst Lasso nuisance. Gaussian setup with parameters: $L=0.05$, $\sigma=0.1$, $s=2$, $\|x\|_1=10$, $\bar{\gamma}=0.1$. \label{fig:conv4}
  }
\end{figure}

\begin{figure}[ht]
  \centering
\begin{tabular}{cc}
    \includegraphics[width=0.45\textwidth]{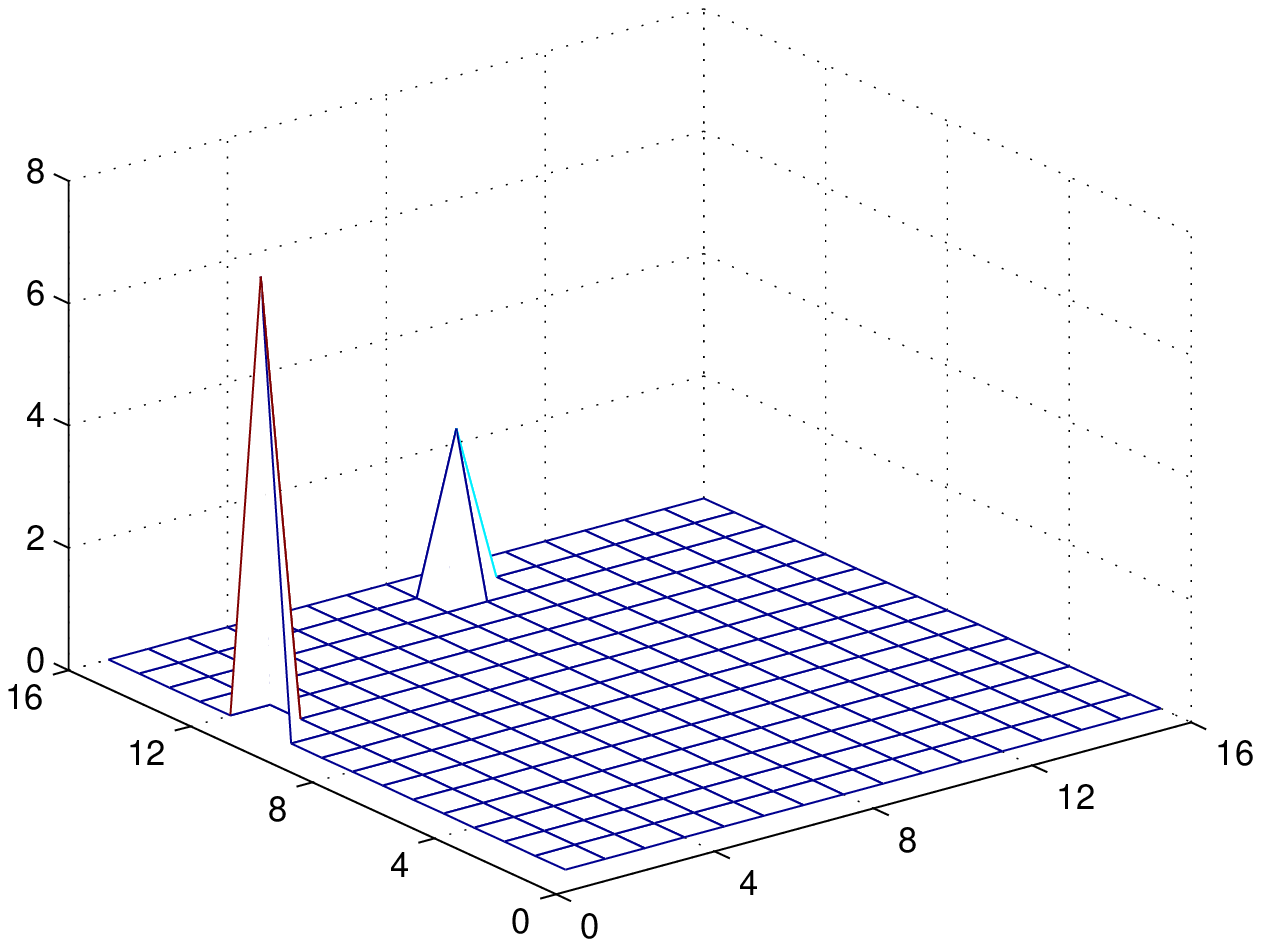} &
    \includegraphics[width=0.45\textwidth]{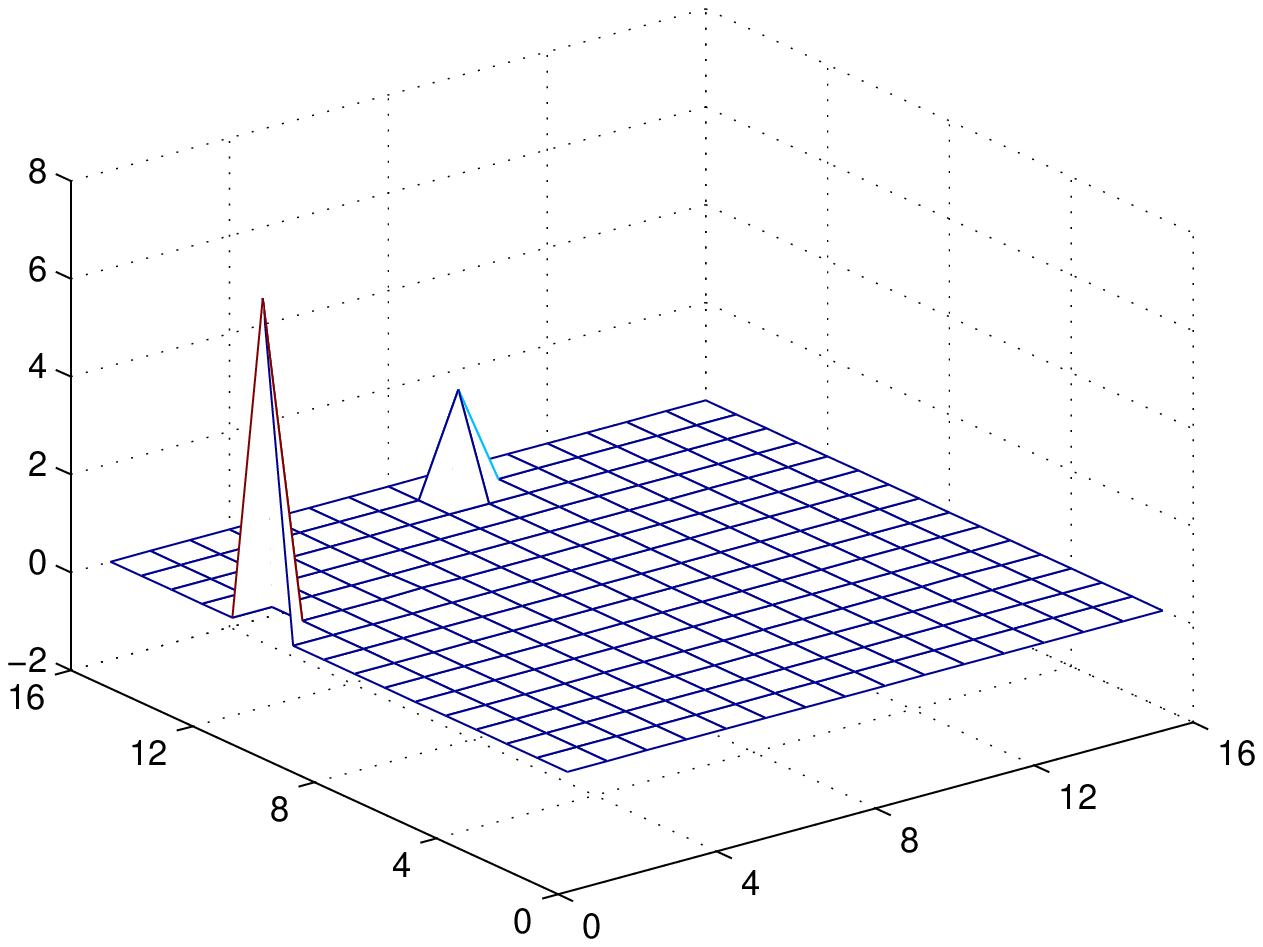} \\
        signal $x$ & recovery $\wh x_\pen$\\
    \includegraphics[width=0.45\textwidth]{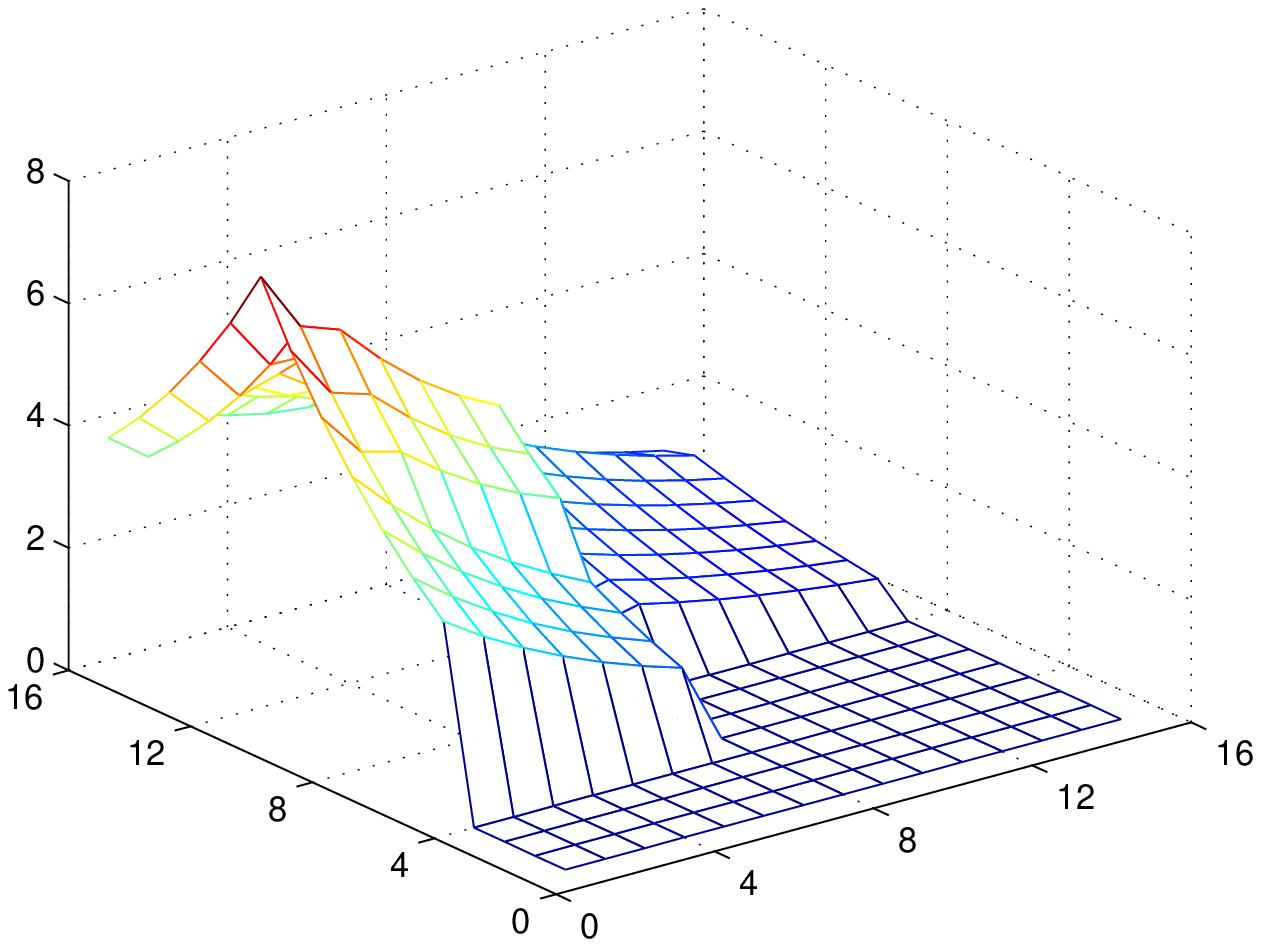} &
    \includegraphics[width=0.45\textwidth]{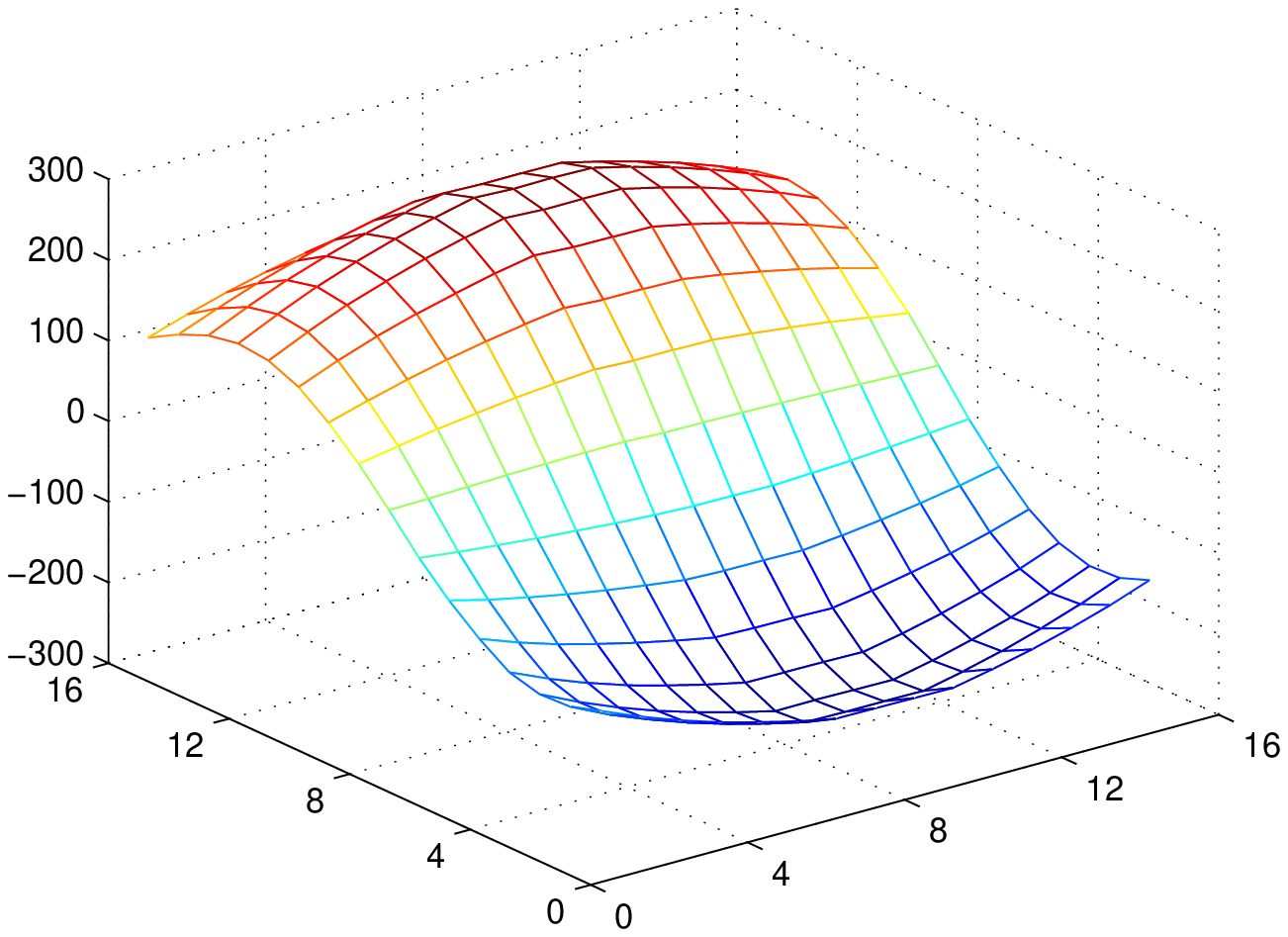} \\
        signal image $Ax$  & observation $y$\\
 \end{tabular}
  \caption{A typical signal/recovery in Convolution setup. Parameters: $L=0.025$, $\sigma=0.1$, $s=2$, $\|x\|_1=10$, $\bar{\gamma}=0.2$. \label{fig:conv40}
  }
\end{figure}
\section{Non-Euclidean matching pursuit algorithm}\label{sectmatpur}\label{MPursuit}
\label{sec:mp}
The Matching Pursuit algorithm for sparse recovery  is motivated by the desire to provide a reduced complexity alternative to the algorithms using $\ell_1$-minimization. Several implementations of Matching Pursuit has been proposed in the Compressive Sensing literature (see, e.g., \cite{EB,DE,GN}). They are based on successive Euclidean projections of the signal and the corresponding performance
results rely upon the bounds on mutual incoherence parameter $\mu(A)$ of the sensing matrix. We are about to show how the construction of Section \ref{matrixh} can be used to design a specific version of the Matching Pursuit algorithm which we refer to as {\em Non-Euclidean Matching Pursuit (NEMP) algorithm}. The NEMP algorithm can be an interesting option if the $\ell_1$-recovery is to be used repeatedly on the observations obtained with the same sensing matrix $A$; the numerical complexity of the pursuit algorithm for a given matrix $A$ may only be a fraction of that of the recovery, especially when used on high-dimensional data.

Suppose that we have in our disposal $\bar{\gamma} \geq 0$ such that the condition $\bH(\bar{\gamma}[1;...,1])$ is feasible; invoking Lemma \ref{arikn}, in this case we can find efficiently a contrast matrix $H=[h_1,...,h_n]$ such that
\be
| [I-H^TA]_{ij}|\leq \bar{\gamma},\,\,\nu(H)=\omega_*(\bar{\gamma}),
\ee{feq99}
where, as always, $\nu(H)=\max\limits_i\nu(h_i)$ with $\nu(h)=\nu_{\epsilon,\sigma,\cU}(h)$ given by \rf{nunorm}.
\par
Consider a signal $x\in \bR^n$ such that $\|x-x^s\|_1\le \upsilon$, where, as usual, $x^s$ is the vector obtained from $x$ by replacing all but the $s$ largest in magnitude entries in $x$ with zeros, and let $y$ be an observation as in \rf{obs}.
\par
Suppose that $s\bar{\gamma}<1$, and let $\upsilon\geq0$ be given. Consider the following iterative procedure:
\begin{algorithm}\label{MP}
$~$
\begin{enumerate}
\item \underline{Initialization:} Set $v^{(0)}=0$, \[
    \alpha_0={\|H^Ty\|_{s,1}+s\nu(H)+\upsilon\over 1-s\bar{\gamma}}.
\]
\item
\underline{Step $k$, $k=1,2,...$:} Given $v^{(k-1)}\in \bR^n$ and $\alpha_{k-1}\geq 0$, compute
\begin{enumerate}
\item $u=H^T(y-Av^{(k-1)})$ and vector
$\Delta\in\bR^n$ with the entries
\[
\Delta_i=\sign(u_i)[|u_i|-\bar{\gamma} \alpha_{k-1}-\nu(H)]_+,\;\;1\le i\le n
\]
(here $[a]_+=\max[0,a]$).
\item
Set  $v^{(k)} = v^{(k-1)} + \Delta$ and
\begin{equation}\label{finitedif}
\alpha_k=2s\bar{\gamma}\alpha_{k-1}+2s\nu(H)+\upsilon.
\end{equation}
and loop to step $k+1$.
\end{enumerate}
\item The approximate solution found after $k$ iterations is $v^{(k)}$.
\end{enumerate}
\end{algorithm}

\begin{proposition}\label{Prop_NEMP} Assume that $s\bar{\gamma}<1$ and an $\upsilon\geq0$ is given.
Then there exists a set $\Xi\subset\bR^m$, $\Prob\{\xi\in\Xi\}\geq1-\e$, of ''good'' realizations of $\xi$ such that whenever $\xi\in\Xi$, for every $x\in\bR^n$ satisfying $\|x-x^s\|_1\leq\upsilon$ and every $u\in\cU$,
the approximate solution $v^{(k)}$ and the value $\alpha_k$ after the $k$-th step of   Algorithm \ref{MP} satisfy
\bse
&(a_{k})&\mbox{for all $i$}\;\;v^{(k)}_i\in\Conv\{0;x_i\}\\
&(b_{k})&\|x-v^{(k)}\|_1
\leq\alpha_{k}\;\;\mbox{and}\;\; \|x-v^{(k+1)}\|_\infty\le 2\bar{\gamma}\alpha_{k}+2\omega_*(\bar{\gamma}).
\ese
\end{proposition}
\par
Note that if $2s\bar{\gamma}<1$ then also $s\bar{\gamma}<1$ and  Proposition \ref{Prop_NEMP} holds true. Furthermore, by (\ref{finitedif}) the sequence $\alpha_k$ converges exponentially fast to the limit $\alpha_\infty:={2s\omega_*(\bar{\gamma})+\upsilon\over 1-2s\bar{\gamma}}$:
\bse
\|v^{(k)}-x\|_1\le \alpha_k=(2s\bar{\gamma})^k[\alpha_0-\alpha_\infty] +\alpha_\infty.
\ese
Along with the second inequality of $(b_k)$ this implies  the bounds:
\bse
\|v^{(k)}-x\|_\infty\le 2\bar{\gamma} \alpha_{k-1}+2\omega_*(\bar{\gamma})\le {\alpha_{k}\over s},
\ese
and, since $\|x\|_p\le \|x\|_1^{1\over p}\|x\|_\infty^{p-1\over p}$ for $1\le p\le \infty$,
\[
\|v^{(k)}-x\|_p\le s^{1-p\over p}\left((2s\bar{\gamma})^k[\alpha_0-\alpha_\infty] +\alpha_\infty\right).
\]
The bottom line here is as follows:
\begin{corollary}\label{MPB} Let $\bar{\gamma}<1/(2s)$ be such that the condition $\bH(\bar{\gamma}[1;...;1])$ is feasible, so that we can find efficiently a contrast matrix $H$ satisfying (\ref{feq99}). With Algorithm \ref{MP} associated with $H$ and some $\upsilon\geq0$, one ensures
 that for every $t=1,2,...$, the approximate solution $v^{(t)}$ found after $t$ iterations satisfies
\[
\begin{array}{c}
\Risk_p(v^{(t)}|\e,\sigma,s,\upsilon)\leq s^{1\over p}\left({2\omega_*(\bar{\gamma})+s^{-1}\upsilon \over1-2s\bar{\gamma}}+(2s\bar{\gamma})^{t}
\left[{\omega_*(\bar{\gamma})+s^{-1}(\|H^Ty\|_{s,1}+\upsilon)\over1-s\bar{\gamma}}
-{2\omega_*(\bar{\gamma})+s^{-1}\upsilon\over 1-2s\bar{\gamma}}\right]\right),\,\,
1\leq p\leq\infty.\\
\end{array}
\]
(cf. (\ref{risksrisks})).
\end{corollary}
To put this result in proper perspective, note that the mutual incoherence based condition
$$
{\mu(A)\over 1+\mu(A)}<{1\over 2s}
$$
underlying typical convergence results for the Matching Pursuit algorithms as applied to recovery of $s$-sparse signals (see, e.g. \cite{EB,DE,GN}) definitely is sufficient for convergence of the NEMP algorithm with $\bar{\gamma}={\mu(H)\over \mu(H)+1}$, see Section \ref{sectmuinc}. It follows that the scope of NEMP is at least as wide as that of ``theoretically valid'' Matching Pursuit algorithms known from the literature; in the situation in question Corollary \ref{MPB} recovers some results from
\cite{DE,EB,GN}.

\appendix
\section{Proofs}
\subsection{Proofs for section \ref{ell1rout}}
\subsubsection{Proof or Lemma \ref{newlemma}}
The first claim is evident. Now let $H=[h_1,...,h_n]$ satisfy $\bH_{s,\infty}(\kappa)$, and let $h_{n+i}=-h_i$, $1\leq i\leq n$.  Then for every $i\leq n$ and every $x\in\bR^n$ with $\|x\|_1\leq1$ we have
$$
{\kappa\over s}\|H^TAx\|_\infty={\kappa\over s}+\max\limits_{\lambda\in\bR^{2n}_+:\sum_i\lambda_i=1}\left[{\sum}_{i=1}^{2n}\lambda_ih_i\right]^TAx\geq x_i,
$$
or, which is the same,
$$
\min\limits_{x:\|x\|_1\leq 1}\max\limits_{\lambda\in\bR^{2n}_+:\sum_i\lambda_i=1}\left[{\kappa\over s}+\left[{\sum}_{i=1}^{2n}\lambda_ih_i\right]^TAx-x_i\right]\geq0.
$$
By von Neumann lemma, this is the same as
$$
\max\limits_{\lambda\in\bR^{2n}_+:\sum_i\lambda_i=1}\min\limits_{x:\|x\|_1\leq 1}\left[{\kappa\over s}+\left[{\sum}_{i=1}^{2n}\lambda_ih_i\right]^TAx-x_i\right]\geq0,
$$
and the outer $\max$ clearly is achieved, meaning that there exists $\lambda^i\geq0$, $\sum_{j=1}^{2n}\lambda^i_j=1$, such that with $h^\prime_i=\sum_{j=1}^{2n}\lambda^i_jh_j$ one has
$$
{\kappa\over s}+[h^\prime_i]^TAx-x_i\geq0\,\,\forall x:\|x\|_1\leq1,
$$
so that for every $x$ with $\|x\|_1\leq1$ one has ${\kappa\over s}+|[h^\prime_i]^TAx|-x_i\geq0$; applying the latter inequality to $-x$ in the role of $x$, we get
${\kappa\over s}+|[h^\prime_i]^TAx|\geq|x_i|$ whenever $\|x\|_1\leq1$, whence, of course,
${\kappa\over s}\|x\|_1+|[h^\prime_i]^TAx|\geq|x_i|$ for all $x$. We conclude that the matrix $H'=[h_1^\prime,...,h_n^\prime]$ satisfies $\bH({\kappa\over s}[1;...;1])$. It remains to note that by construction the columns of $H'$ are convex combinations of the columns of $H$ and $-H$, and that building $H'$ reduces to solving $n$ matrix games and thus can be carried out efficiently. \epr
\subsubsection{Proof of Proposition \ref{prop10}}\label{proofprop10}
Let
\[
\Xi=\{\xi:\;|h_i^Ty|\leq\sqrt{2\ln(n/\e)}\|h_i\|_2\;1\le i\leq n\},
 \]
 so that $\Prob\{\xi\in\Xi\}\geq1-\epsilon.$ Let us fix $\xi\in\Xi$,
 a set $I=\{1,...,n\}\backslash J\subset\{1,...,n\}$ satisfying (\ref{suchthat}), a signal $x\in\bR^n$ and a realization $u\in\cU$ of the nuisance, and let $\widehat{x}$ be the value of the estimate (\ref{regular}) at the observation $y=Ax+u+\sigma\xi$. We are about to verify that $\widehat{x}$ satisfies (\ref{r100}), which, of course, will complete the proof. \par
  Observe that because of $\xi\in \Xi$ we have
$$
|h_i^T(Ax-y)|= |h_i^T(u+\sigma\xi)|\leq \max_{u'\in\cU}|h_i^Tu'|+\sigma\sqrt{2\ln(n/\epsilon)}\|h_i\|_2=\nu(h_i)=\nu_i,\;\;1\le i\le n.$$
Now, $\rho_i\geq\nu_i$ by (\ref{rhoi}), whence $|h_i^T(y-Ax)|\leq \rho_i$ for all $i$, and thus $x$ is a feasible solution to the optimization problem in (\ref{regular}) and thus $\|x\|_1\geq\|\widehat{x}\|_1$. Setting $z=\widehat{x}-x$, we now have $\|\widehat{x}_J\|_1=\|\widehat{x}\|_1-\|\widehat{x}_I\|_1\leq \|x\|_1-\|\widehat{x}_I\|_1\leq \|x\|_1-\|x_I\|_1+\|z_I\|_1=\|x_J\|_1+\|z_I\|_1$,
whence $\|z_J\|_1\leq \|\widehat{x}_J\|_1+\|x_J\|_1\leq 2\|x_J\|_1+\|z_I\|_1$. It follows that
\begin{equation}\label{eqeqeq}
\|z\|_1\leq 2\|z_I\|_1+2\|x_J\|.
\end{equation}
Further, $|h_i^TA(\widehat{x}-x)|\le |h_i^T(A\widehat{x}-y)|+|h_i^T(Ax-y)|$. Since $\widehat{x}$ is feasible for the optimization problem in (\ref{regular}), we have
$|h_i^T(A\widehat{x}-y)|\leq\rho_i$, and we have already seen that $|h_i^T(Ax-y)|\leq\nu_i$, hence
\be
|h_i^TAz|
\le \rho_i+\nu_i
\ee{delta}
for all $1\le i\le n$. Applying \rf{cond00} we now get
\bse
\|z_I\|_1&=&\sum_{i\in I}|z_i|\le \sum_{i\in I}[|h_i^TAz|+\gamma_i \|z\|_1]
\le \sum_{i\in I}(\rho_i+\nu_i)+\left[\sum_{i\in I}\gamma_i\right]\left[\|z_I\|_1+\|z_{J}\|_1\right]\\
&=& \rho_I+\nu_I+\gamma_I\left[\|z_I\|_1+\|z_{J}\|_1\right]
\leq \rho_I+\nu_I+2\gamma_I\left[\|z_I\|_1+\|x_{J}\|_1\right],
\ese
where the concluding $\leq$ is given by (\ref{eqeqeq}). Taking into account that $\gamma_I<\half$, we get
$$
\|z_I\|_1\leq {\rho_I+\nu_I+2\gamma_I\|x_{J}\|_1\over1-2\gamma_I}.
$$
Invoking (\ref{eqeqeq}) once again, we finally get
$$
\|z\|_1=\|z_I\|_1+\|z_{J}\|_1\leq 2\|z_I\|_1+2\|x_{J}\|_1\leq 2{
\rho_I+\nu_I+2\gamma_I\|x_{J}\|_1\over1-2\gamma_I}+2\|x_{J}\|_1,
$$
and we arrive at (\ref{r100}.$a$).
\par To prove (\ref{r100}.$b$), we  apply \rf{cond00} to $z=\widehat{x}-x$, thus getting
$$
|z_i|\leq |h_i^TAz|+\gamma_i\|z\|_1.
$$
As we have already seen, $|h_i^TAz|\leq\rho_i+\nu_i$, and the first ``$\leq$'' in (\ref{r100}.$b$) follows; the second ``$\leq$'' in (\ref{r100}.$b$) is then readily given by (\ref{r100}.$a$). Now (ii) and (iii) are immediate consequences of (\ref{r100}) and the fact that $\widehat{\gamma}_s\leq s\widehat{\gamma}$. \epr
\subsubsection{Proof of Lemma \ref{allp}}\label{proofcorallp}  In what follows, we use the notation from Proposition \ref{prop10}. For  $x\in X(s,\upsilon)$, denoting by  $I$ the support of $x^s$, we have
$$
\|x_{J}\|_1\leq\upsilon,\,\rho_I\leq\widehat{\rho}_s\leq s\widehat{\rho},\,\nu_I\leq\widehat{\nu}_s\leq s\nu(H),\,\gamma_I\leq \widehat{\nu}_s\leq s\widehat{\gamma}.
$$
Assuming $\widehat{\gamma}_s<\half$, for $\xi\in\Xi$ (which happens with probability $\geq1-\epsilon$), \rf{therr} implies that for all $u\in\cU$ it holds
\bse
\|\widehat{x}_{\reg}(y)-x\|_1\leq \underbrace{{2\over1-2\widehat{\gamma}_s}[\upsilon+\widehat{\rho}_s+\widehat{\nu}_s]}_{P},\;\;\mbox{and}\;\;
\|\widehat{x}_{\reg}(y)-x\|_\infty\leq \underbrace{\widehat{\rho}+\nu(H)+2\widehat{\gamma}{\upsilon+\widehat{\rho}_s+\widehat{\nu}_s\over1-2\widehat{\gamma}_s}}_{Q},
\ese
which combines with the standard bound $\|z\|_p\leq\|z\|_1^{{1\over p}}\|z\|_\infty^{{p-1\over p}}$ to imply (\ref{rl1}). When $s\widehat{\gamma}<\half$, we clearly have
\bse
P\leq {2\over 1-2s\widehat{\gamma}}[\upsilon+s(\widehat{\rho}+\nu(H))],\,\,
Q\leq \widehat{\rho}+\nu(H)+{2\widehat{\gamma}\over 1-2s\widehat{\gamma}}[\upsilon+s(\widehat{\rho}+\nu(H))]={2\over 1-2s\widehat{\gamma}}[\widehat{\gamma}\upsilon+\half[\widehat{\rho}+\nu(H)]],
\ese
and (\ref{also}) follows due to $\|\widehat{x}_{\reg}(y)-x\|_p=P^{{1\over p}}Q^{{p-1\over p}}$. \epr
\subsubsection{Proof of Proposition \ref{prop101}}
The proof is obtained by minor modifications from the one of Proposition \ref{prop10}. Same as in the latter proof, let $\Xi=\{\xi:\;|h_i^Ty|\leq\sqrt{2\ln(n/\e)}\|h_i\|_2\;1\le i\leq n\}$, where $h_i$ are the columns of $H$, so that $\Prob\{\xi\in\Xi\}\geq1-\epsilon.$ Let us fix $\xi\in\Xi$, $x\in \bR^n$, $u\in\cU$, let $\eta=\sigma\xi+u$,  $y=Ax+\eta$, $\widehat{x}=\widehat{x}_{\reg}(y)$, $z=\widehat{x}-x$. Finally, let $I$ be the support of $x^s$.
\par
Due to $\xi\in\Xi$, we have
$$
|h_i^T(Ax-y)|= |h_i^T(u+\sigma\xi)|\leq \max_{u'\in\cU}|h_i^Tu'|+\sigma\sqrt{2\ln(n/\epsilon)}\|h_i\|_2=\nu(h_i)=\nu_i,\;\;1\le i\le n,$$
whence, by (\ref{rhoi}), $x$ is a feasible solution to the optimization problem in (\ref{regular}) and thus $\|x\|_1\geq\|\widehat{x}\|_1$. The latter, exactly as in the proof of Proposition \rf{prop10}, implies the validity of (\ref{eqeqeq}):
\begin{equation}\label{eqeqeqeq}
\|z\|_1\leq 2\|z_I\|1+2\|x_J\|_1.
\end{equation}
Besides this, the same reasoning as in the proof of Proposition \ref{prop10} results in (\ref{delta}), whence
\begin{equation}\label{deltanew}
\|H^TAz\|_\infty\leq \widehat{\rho}+\nu(H).
\end{equation}
Applying (\ref{Hsinfty}) to $z$, we get
$$
\|z_I\|_1\leq s\|z\|_\infty\leq s\|H^TAz\|_\infty+\kappa \|z\|_1\leq s(\widehat{\rho}+\nu(H))+\kappa \|z\|_1,
$$
which combines with (\ref{eqeqeqeq}) to imply that
\begin{equation}\label{normzone}
\|z\|_1\leq {1\over 1-2\kappa}\left[2s(\widehat{\rho}+\nu(H))+2\|x_J\|_1\right],
\end{equation}
which is nothing but the first relation in (\ref{therrnew}).
Applying to $z$ (\ref{Hsinfty}) once again,
we get
$$
\|z\|_\infty\leq \|H^TAz\|_\infty+s^{-1}\kappa\|z\|_1,
$$
which combines with (\ref{normzone}) to imply the second relation in (\ref{therrnew}). Relation (\ref{therrnew}) combines with the Moment inequality  to imply (\ref{alsonew}).
\epr

\subsubsection{Proof of Proposition \ref{prop1}}\label{proofprop1}
\paragraph{(i):} Given $\sigma$, $\epsilon$, let, same as in the proof of Proposition \ref{prop10},  $\Xi=\{\xi:|h_i^T\xi|\leq \sqrt{2\ln(n/\epsilon)}\|h_i\|_2,\;1\le i\leq n\},$ so that $\Prob\{\xi\in\Xi\}\geq1-\epsilon$.
Let us fix $\xi\in\Xi$, $u\in\cU$ and a signal $x\in\bR^n$, and let us prove that for these data (\ref{r200}) takes place; this clearly will prove (i). Let us set $y=Ax+\sigma\xi+u$, $\widehat{x}=\widehat{x}_{\pen}(y)$, $z=\widehat{x}-x$, $\eta=u+\sigma\xi$. Let also $I$ be the support of $x^s$.
\par
 Observe that by the origin of $\widehat{x}$, we have
\begin{equation}\label{eq1}
\|\widehat{x}\|_1+s\theta\|H^T(A\wh{x}-y)\|_\infty\leq \|x\|_1+s\theta\|H^T(Ax-y)\|_\infty=\|x\|_1+s\theta\|H^T\eta\|_\infty,
\end{equation}
and $$\|H^T(A\wh{x}-y)\|_\infty=\|H^T(Az+Ax-y)\|_\infty\geq\|H^TAz\|_\infty-\|H^T(Ax-y)\|_\infty
=\|H^TAz\|_\infty-\|H^T\eta\|_\infty.$$ Combining the resulting inequality with (\ref{eq1}), we get
\begin{equation}\label{eq2}
\|\widehat{x}\|_1+s\theta\|H^TAz\|_\infty\leq \|x\|_1+2s\theta\|H^T\eta\|_\infty\leq \|x\|_1+2s\theta\nu(H),
\end{equation}
where the concluding $\leq$ is due to $\xi\in\Xi$ combined with (\ref{nui}).
Further,
$$\|\widehat{x}\|_1=\|x+z\|_1=\|x_I+z_I\|_1+\|x_J+z_{J}\|_1\geq \|x_I\|_1-\|z_I\|_1+\|z_{J}\|_1-\|x_{J}\|_1,
 $$
 which combines with (\ref{eq2}) to imply that
$$
\|x_I\|_1-\|z_I\|_1+\|z_{J}\|_1-\|x_{J}\|_1+s\theta\|H^TAz\|_\infty\leq \|x\|_1+2s\theta\nu(H),
$$
or, which is the same,
\begin{equation}\label{eq3}
\|z_{J}\|_1-\|z_I\|_1+s\theta\|H^TAz\|_\infty\leq 2\|x_{J}\|_1+2s\theta\nu(H).
\end{equation}
By \rf{cond00}, we have
\begin{equation}\label{eq4a}
\forall i: |z_i|\leq \|H^TAz\|_\infty+\gamma_i\|z\|_1,
\end{equation} whence $\|z_I\|_1\leq s\|H^TAz\|_\infty+\widehat{\gamma}_s\|z\|_1$ and therefore

\[(1-\widehat{\gamma}_s)\|z_I\|_1-\widehat{\gamma}_s\|z_{J}\|_1-s\|H^TAz\|_\infty\leq0.
\]
Multiplying the latter inequality by $\theta$ and summing up with (\ref{eq3}), we get
$$
[\theta(1-\widehat{\gamma}_s)-1]\|z_I\|_1+(1-\theta \widehat{\gamma}_s)\|z_{J}\|_1\leq 2\|x_{J}\|_1+2s\theta\nu(H).
$$
In view of condition (\ref{theta}), the coefficients in the left hand side are positive, and (\ref{r200}.$a$) follows.
\par
 To prove (\ref{r200}.$b$), note that from (\ref{eq2}) it follows that
$$
\|H^TAz\|_\infty\leq {1\over s\theta }[\|x\|_1-\|\widehat{x}\|_1]+2\nu(H)\leq {1\over s\theta} \|z\|_1+2\nu(H),
$$
which combines with (\ref{eq4a}) to imply that
$$
\|z\|_\infty\leq {1\over s\theta}\|z\|_1+2\nu(H)+\widehat{\gamma}\|z\|_1.
$$
Recalling that $z=\widehat{x}-x$ and invoking (\ref{r200}.$a$), (\ref{r200}.$b$) follows.
\paragraph{(ii)--(iii):} (\ref{r200a}) is an immediate consequence of (\ref{r200}) due to $\widehat{\gamma}_s\leq s\widehat{\gamma}$. Assuming that $x\in X(s,\upsilon)$ and $\theta=2$ and taking into account that $\widehat{\gamma}_s\leq s\widehat{\gamma}$, we obtain from \rf{r200a} that uniformly on $\xi\in \Xi$ and $u\in \cU$
\bse
\|\bar{x}_{\reg}(y)-x\|_1\leq \underbrace{\left[{2\upsilon+4s\nu(H)\over1-2s\widehat{\gamma}}\right]}_{P},\;\;
\|\bar{x}_{\reg}(y)-x\|_\infty\leq \underbrace{\left[{(s^{-1}+2\widehat{\gamma})\upsilon+4\nu(H)\over 1-2s\widehat{\gamma}}\right]}_{Q}.
\ese
Using, as in the proof of Lemma \ref{allp}, the standard bound $$\|z\|_p\leq\|z\|_1^{{1\over p}}\|z\|_\infty^{{p-1\over p}}\le P^{1\over p}Q^{{p-1\over p}}
$$
we come to \rf{risks}. \epr
\subsubsection{Proof of Proposition \ref{prop11}}
The proof is obtained by minor modifications from the one of Proposition \ref{prop10}. Same as in the latter proof, let  $\Xi=\{\xi:|h_i^T\xi|\leq \sqrt{2\ln(n/\epsilon)}\|h_i\|_2,\;1\le i\leq n\},$ so that $\Prob\{\xi\in\Xi\}\geq1-\epsilon$.\par
Let us fix $\xi\in\Xi$, $u\in\cU$ and a signal $x\in\bR^n$. Let us set $y=Ax+\sigma\xi+u$, $\widehat{x}=\widehat{x}_{\pen}(y)$, $z=\widehat{x}-x$, $\eta=u+\sigma\xi$. Let also $I$ be the support of $x^s$.
\par
 Observe that by the origin of $\widehat{x}$ and due to $\theta=2$ we have
\begin{equation}\label{eq1new}
\|\widehat{x}\|_1+s\theta\|H^T(A\wh{x}-y)\|_\infty\leq \|x\|_1+s\theta\|H^T(Ax-y)\|_\infty=\|x\|_1+2s\|H^T\eta\|_\infty,
\end{equation}
and $$\|H^T(A\wh{x}-y)\|_\infty=\|H^T(Az+Ax-y)\|_\infty\geq\|H^TAz\|_\infty-\|H^T(Ax-y)\|_\infty
=\|H^TAz\|_\infty-\|H^T\eta\|_\infty.$$ Combining the resulting inequality with (\ref{eq1new}), we get
\begin{equation}\label{eq2new}
\|\widehat{x}\|_1+2s\|H^TAz\|_\infty\leq \|x\|_1+4s\|H^T\eta\|_\infty\leq \|x\|_1+4s\nu(H),
\end{equation}
where the concluding $\leq$ is due to $\xi\in\Xi$ combined with the definition of $\nu(H)$.
Further,
$$\|\widehat{x}\|_1=\|x+z\|_1=\|x_I+z_I\|_1+\|x_J+z_{J}\|_1\geq \|x_I\|_1-\|z_I\|_1+\|z_{J}\|_1-\|x_{J}\|_1,
 $$
 which combines with (\ref{eq2new}) to imply that
$$
\|x_I\|_1-\|z_I\|_1+\|z_{J}\|_1-\|x_{J}\|_1+2s\|H^TAz\|_\infty\leq \|x\|_1+4s\nu(H),
$$
or, which is the same,
\begin{equation}\label{eq3new}
\|z_{J}\|_1-\|z_I\|_1+2s\|H^TAz\|_\infty\leq 2\|x_{J}\|_1+4s\nu(H).
\end{equation}
By (\ref{Hsinfty}) we have
\begin{equation}\label{eq4anew}
\|z\|_\infty\leq \|H^TAz\|_\infty+{\kappa\over s}\|z\|_1,
\end{equation} whence $\|z_I\|_1\leq s\|H^TAz\|_\infty+\kappa\|z\|_1$ and therefore

\[(1-\kappa)\|z_I\|_1-\kappa\|z_{J}\|_1-s\|H^TAz\|_\infty\leq0.
\]
Multiplying the latter inequality by 2 and summing up with (\ref{eq3new}), we get
$$
(1-2\kappa)\|z\|_1\leq 2\|x_{J}\|_1+4s\nu(H),
$$
and the first relation in (\ref{therrnewnew}). The second relation in (\ref{therrnewnew}) is readily given by the first one combined with (\ref{Hsinfty}). We have proved that (\ref{therrnewnew})
holds true whenever $\xi\in\Xi$; since $\Prob\{\xi\in\Xi\}\geq1-\epsilon$, (\ref{alsonewnew}) follows. \epr
\subsection{Proofs for sections \ref{matrixh}, \ref{LimitsofPerformance}}
\subsubsection{Proof of Lemma \ref{arikn}}\label{proofarikn}
(i)$\Rightarrow$(iii): If $h_i$ satisfies $(\cP_i)$, then  for every $x$ we have
\[|x_i|\leq |h_i^TAx|+\gamma\|x\|_1\leq \omega\varphi_*(Ax)+\gamma\|x\|_1,
 \]where the first and the second inequalities are given by $(\cP_i.b)$ and $(\cP_i.a)$, respectively. \qed\\
(iii)$\Rightarrow$(ii): Assume that (iii) takes place; then, by homogeneity, $\omega\varphi_*(Ax)+\gamma\geq x_i$ for every $x$ with $\|x\|_1\leq1$, or, which is the same, the optimal value in the
conic problem \[
\min_x\left\{\omega\varphi_*(Ax)-x_i:\|x\|_1\leq1\right\}
 \]
 is $\geq-\gamma$. The problem clearly is strictly feasible and bounded, so that by Conic Duality Theorem the dual problem is solvable with the same optimal value. Now, the dual problem reads \[
 \max\limits_{g,h,s}\{-s:\;\varphi(h)\leq\omega,\,\|g\|_\infty\leq s,\,A^Th+g=e_i\},
  \]and the fact that it is solvable with the optimal value $\geq-\gamma$ means that there exist $h,g$ such that $\varphi(h)\leq\omega$, $\|g\|_\infty\leq\gamma$ and $A^Th+g=e_i$, whence $h$ is a feasible solution to $(P_i^\gamma)$ with the value of the objective $\leq\omega$. \qed
\medskip\\
(ii)$\Rightarrow$(i): If $(P_i^\gamma)$ is feasible, it clearly is solvable; thus, in the case of (ii) there exists $h$ with $\varphi(h)\leq\omega$ and $\|A^Th-e_i\|_\infty\leq\gamma$. From the latter inequality it follows that $|e_i^Tx-h^TAx|\leq\gamma\|x\|_1$ for every $x$, so that $|x_i|-|h^TAx|\leq\gamma\|x\|_1$ for all $x$. We see that $h$ satisfies $(\cP_i)$, and thus (i) takes place. This reasoning shows also that whenever $(P_i^\gamma)$ is feasible with optimal value $\leq\omega$, it is solvable, and its optimal solution satisfies $(\cP_i)$. \epr

\subsubsection{Proof of Proposition \ref{propRIP}}\label{proofpropRIP}
Let  $\gamma=\gamma(\delta,k)$,  $\lambda=r(\cU)+\sigma\sqrt{2\ln(n/\epsilon)}$, so that what we need to prove is that there exists a matrix $H$ satisfying
$\bH_{s,\infty}(\gamma)$ and such that $\nu(H)\leq\lambda$. Invoking Lemma \ref{arikn}, all we need to this end is to show that
\begin{equation}\label{goal}
\forall x\in\bR^n: \|x\|_\infty\leq {1\over\sqrt{1-\delta}}\lambda\nu_*(Ax)+\gamma\|x\|_1
\end{equation}
Now, we clearly have $\nu(h)\leq\max\limits_{u\in\cU} u^Th+\sigma\sqrt{2\ln(n/\epsilon)}\|h\|_2\leq \lambda\|h\|_2$ for all $h$, whence $\phi_*(\eta)\geq \lambda^{-1}\|\eta\|_2$ for all $\eta$. Therefore all we need in order to justify (\ref{goal}) is to prove that
\begin{equation}\label{goal1}
\forall s\in\bR^n: \|x\|_\infty\leq {1\over\sqrt{1-\delta}}\|Ax\|_2+\gamma\|x\|_1.
\end{equation}
Let  $x\in\bR^n$.
Setting $s=\hbox{floor}(k/2)$,  let vectors $x^1,...,x^q$ be obtained from $x$ by the procedure as follows: $x^1$ is obtained by zeroing all but the $s$ largest in magnitude entries of $x$; $x^2$ is obtained by the same procedure from $x-x^1$, $x^3$ is obtained by the same procedure from $x-x^1-x^2$, and so on, until the step $q$ where we get $x=x^1+...+x^q$.
 We clearly have
$\|x^j\|_\infty\leq s^{-1}\|x^{j-1}\|_1$, $2\leq j\leq q$,  whence also $\|x^j\|_2\leq s^{-1/2}\|x^{j-1}\|_1$, $2\leq j\leq q$, since the vectors $x^j$ are $s$-sparse.
Setting
$\|Ax\|_2=\alpha$ and $\|Ax^1\|_2=\beta$, we have
$$
\alpha\beta= \|Ax\|_2\|Ax^1\|_2\geq    (Ax)^TAx^1
=[x^1]^TA^TAx^1+\sum_{j=2}^q[x^1]^TA^TAx^j
\geq \beta^2-\sum_{j=2}^q \delta\|x^1\|_2\|x^j\|_2,
$$
where the last $\geq$ is given by the following well-known fact:
\cite{Candes08}: \begin{quote}
(!) {\sl If $A$ is $\RIP(\delta,k)$ and $u,v$ are supported on a
common set of indices $I$ of cardinality $k$ and are orthogonal, then $|u^TA^TAv|\leq \delta\|u\|_2\|v\|_2$.}
\end{quote}
It follows that
$$\alpha\beta\geq\beta^2-\delta \|x^1\|_2\sum_{j=2}^q\|x^j\|_2
\geq\beta^2-\delta s^{-1/2}\|x^1\|_2\sum_{j=2}^q\|x^{j-1}\|_1
\geq\beta^2-\delta s^{-1/2}\|x^1\|_2\|x\|_1.
$$
Hence
\bse
\beta\le \alpha+{\delta\|x\|_1\|x^1\|_2\over \sqrt{s}\beta}\le
 \alpha+{\delta\|x\|_1\over \sqrt{s(1-\delta)}},
\ese
where the second inequality is due to the fact  that $\|x^1\|_2/\beta\leq 1/\sqrt{1-\delta}$ by RIP. Thus,
\[
\|x\|_\infty\leq\|x^1\|_2\le {\beta\over \sqrt{1-\delta}}\le  {\alpha\over \sqrt{1-\delta}}+{\delta\|x\|_1\over (1-\delta)\sqrt{s}}\leq
{\alpha\over\sqrt{1-\delta}}+\gamma\|x\|_1,
\]
where the concluding inequality is due to $s\geq (k-1)/2$ and $\gamma=\gamma(\delta,k)$. Recalling that
$\alpha=\|Ax\|_2$, (\ref{goal1}) follows. \epr
\subsubsection{Proof of Proposition \ref{propOracle}}
{\bf Proof.} We start with analysis of $\bO(S,\omega)$. Let $i\in\{1,...,n\}$, and let $I \ni i$ be a subset of
$\{1,...,n\}$ of cardinality $S$. Let $\bR^{(S)}$ be the linear space of all vectors from $\bR^n$ supported on $I$,  and let $X_R=\{x\in \bR^{(S)}:\|x\|_2\leq R\}$. Assume that we are given a noisy observation $y=Ax+u+\sigma\xi$ of a signal
$z=(x,u)\in (X_R,\cU)$, and that we want to recover from this observation the linear form $x_i$ of the signal. From
 $\bO(S,\omega)$ it follows that there exists a recovering routine such that for every $x\in X_R$ and $u\in \cU$ the probability of recovering error
to be $\geq \omega$ is $\leq \e$. Assuming $\epsilon\leq 1/16$ and applying the celebrated result of Donoho \cite{Don95}, there exists a {\sl linear} estimate $\phi_R^Ty$ such that for every $x\in X_R$ and $u\in \cU$ the probability
for the error of this estimate to be $\ge 1.22\omega$ is $\leq\epsilon$. Moreover (cf. Proposition 4.2 of \cite{YuN}), one can pick $\phi_R$ such that
$$
\forall x\in X_R, \;u\in \cU,\;\;\left\{
\begin{array}{lclr}
\Prob\{\phi_R^T[u+\sigma\xi+A_Ix]-x_i>1.22\omega\}&\leq&{\epsilon/2},&(a)\\
\Prob\{\phi_R^T[u+\sigma\xi+A_Ix]-x_i<-1.22\omega\}&\leq&{\epsilon/2},&(b)\\
\end{array}
\right.$$
where $A_I$ is the matrix obtained from $A$ by zeroing columns with indexes not belonging to $I$.
Let $p(R)=\max\limits_{u\in\cU}|\phi_R^Tu|$ and $r(R)=\|A_I^T\phi_R-e_i\|_2$, where $e_i$ is the $i$-th basic orth (so that
$x_i=e_i^Tx$). Specifying $\bar{x}$ as the vector from $X_R$ such that $\bar{x}^T(A_I^T\phi_R-e_i)=Rr(R)$, and $\bar{u}$ as the vector from $\cU$ such that $\phi_R^T\bar{u}=p(R)$ (the required $\bar{x}$, $\bar{u}$ clearly exist) and applying $(a)$ to the pair $(x,u)=(\bar{x},\bar{u})$, and $(b)$ to the pair $(x,u)=(-\bar{x}),-\bar{u})$, we get
$$
\begin{array}{lcl}
\Prob\{\sigma \phi_R^T\xi>1.22\omega-Rr(R)-p(R)\}&\leq&\epsilon/2,\\
\Prob\{\sigma \phi_R^T\xi<-1.22\omega+Rr(R)+p(R)\}&\leq&\epsilon/2.\\
\end{array}
$$
Hence, denoting by $\ErfInv(\e)$ the value of the inverse error function at $\e$, we obtain
$$
\ErfInv\left({\e\over 2}\right)\sigma\|\phi_R\|_2\leq 1.22\omega-Rr(R)-p(R).
$$
It follows that as $R\to\infty$, $\phi_R$ remains bounded and $r(R)=\|e_i-A_I^T\phi_R\|_2\to0$. Thus, there exists a sequence
$R_k\to+\infty$, of values of $R$ such that $\phi_{R_k}$ goes to a limit $\phi$ as $k\to\infty$, and this limit
satisfies the relations
$$
\ErfInv\left({\e\over 2}\right)\sigma\|\phi\|_2\leq 1.22\omega,\;\;\mbox{and}\;A_I^T\phi=e_i.
$$
Taking into account that $\ErfInv\left({\e\over 2}\right)\geq 0.92\sqrt{\ln(1/\epsilon)}$ when $\epsilon\leq 1/16$, we
arrive at the following result:
\begin{lemma}\label{lem2} Under assumption $\bO(S,\omega)$, for every $i\leq n$ and every $S$-element subset $I\ni i$ of $\{1,...,n\}$
there exists $\phi\in\bR^m$ such that $\phi^Ta_i=1$, $\phi^Ta_j=0$ for all $j\in I$, $j\neq i$ (here $a_1,...,a_n$ are the columns of $A$), and
\begin{equation}\label{phi}
\max\limits_{u\in\cU}|u^T\phi|+\sigma\sqrt{\ln(1/\epsilon)}\|\phi\|_2\leq \sqrt{2}\omega.
\end{equation}
\end{lemma}
We claim that in this case for all $x\in \bR^n$ it holds:
\be
\|x\|_\infty\leq\bar{\omega} \nu_*(Ax)+\overbrace{{\omega\|A\|\over\sigma\sqrt{2S\ln(1/\epsilon)}}}^{\widehat{\gamma}}\|x\|_1.
\ee{claim}
Taking this claim for granted, and invoking Lemma \ref{arikn}, we immediately arrive at the desired conclusion. Indeed, given $s$ satisfying (\ref{range}), we have ${1\over 4s}\geq\widehat{\gamma}$, so that (\ref{claim}) implies that
$$
\forall x\in\bR^n: \|x\|_\infty\leq \bar{\omega}\nu_*(Ax)+{1\over 4s}\|x\|_1,
$$
whence, by Lemma \ref{arikn}, there exists $H$ satisfying the condition $\bH_{s,\infty}(\four)$ and such that $\nu(H)\leq\bar{\omega}$, which is exactly what Proposition \ref{propOracle} states.\par
It remains to prove (\ref{claim}). Let us fix $x\in\bR^n$, and let $I$ be set of indices of the $S$ largest in magnitude entries in $x$. Denoting by $J$ the complement of $I$ in
$\{1,...,n\}$, we have $\|x_J\|_\infty\leq S^{-1}\|x_I\|_1$, whence
\begin{equation}\label{printer}
\|x_J\|_2\leq\|x_J\|_\infty^{1/2}\|x_J\|_1^{1/2}\leq
S^{-1/2}\|x_I\|_1^{1/2}\|x_J\|_1^{1/2}\leq {1\over 2}S^{-1/2}\|x\|_1.
 \end{equation} Let $i_*\in I$ be the index of the largest in magnitude entry of $x$. By Lemma \ref{lem2} there exists $\phi\in\bR^m$ satisfying (\ref{phi}) and such that $\phi^Ta_{i_*}=\sign(x_{i_*})$, $\phi^Ta_i=0$ for $i\in I\backslash \{i_*\}$. We have
\begin{equation}\label{boundphi}
\begin{array}{l}
\nu(\phi)\equiv\nu_{\epsilon,\sigma,\cU}(\phi)=\max\limits_{u\in\cU}u^T\phi+\sigma\sqrt{2\ln(n/\epsilon)}\|\phi\|_2\\
\leq \sqrt{2[1+\ln(n)/\ln(1/\epsilon)]}
\left[\max\limits_{u\in\cU}u^T\phi+\sigma\sqrt{\ln(1/\epsilon)}\|\phi\|_2\right]\leq2 \sqrt{1+\ln(n)/\ln(1/\epsilon)}\omega,
\end{array}
\end{equation}
where the concluding $\leq$ is given by (\ref{phi}). Now,
\bse
\nu(\phi)\nu_*(Ax)&\geq &\phi^TAx=\phi^TAx_I+\phi^TAx_J=|x_{i_*}|+\phi^TAx_J= \|x\|_\infty+\phi^TAx_J\\
&\geq& \|x\|_\infty-\|\phi\|_2\|Ax_J\|_2\geq \|x\|_\infty -\|\phi\|_2\|A\|\|x_J\|_2\geq \|x\|_\infty-{1\over
2}\|A\|S^{-1/2}\|\phi\|_2\|x\|_1,
\ese
with the concluding $\leq$ given by (\ref{printer}).
The resulting inequality, in view of (\ref{boundphi}) and the bound $\|\phi\|_2\leq {\sqrt{2}\omega\over\sigma\sqrt{\ln(1/\epsilon)}}$ given by (\ref{phi}) implies (\ref{claim}).\epr

\subsection{Proofs for section \ref{Extensions}}
\subsubsection{Proof of Lemma \ref{arik1}}\label{prooflemma1}
Recall (cf., e.g.,  Theorem 2.1 in \cite{JNCS}) that a necessary and sufficient condition for an $m\times n$  matrix $A$ to be $s$-good is the {\sl nullspace property} as follows: there exists $\kappa<1/2$ such that
\begin{equation}\label{nullspace}
\|x\|_{s,1}\leq \kappa\|x\|_1\,\,\forall (x\in\bR^n,Ax=0).
\end{equation}
Assume that this condition is satisfied, and let $u\in\bR^n$ be a vector with $s$ nonzero coordinates, equal to $\pm1$.  (\ref{nullspace}) says that the optimal value in the Linear Programming problem
$$
\max\limits_x\left\{u^Tx:Ax=0,\|x\|_1\leq1\right\}
$$
is at most $\kappa$; passing to the dual problem, we conclude that there exist $h_u$ and $g_u$ such that $A^Th_u+g_u=u$ and $\|g_u\|_\infty\leq\kappa$, whence for every $x\in\bR^n$ it holds
$$
u^Tx=x^TA^Th_u+g_u^Tx\leq \|h_u\|_1\|Ax\|_\infty +\kappa\|x\|_1.
$$
Since the set $U$ of the outlined vectors $u$ is finite, the quantity $L=\max\limits_{u\in U}\|h_u\|_1$ is finite, and
$$
\|x\|_{s,1}=\max\limits_{u\in U}u^Tx \leq L\|Ax\|_\infty+\kappa\|x\|_1\,\,\forall x,
$$
meaning that the condition $\bH_{s,1}(\kappa)$ holds true for $H=[LI_m,0_{m\times n-m}]$. Vice versa, the existence of $\kappa<1/2$ and $H$ satisfying $\bH_{s,q}(\kappa)$  clearly implies the validity of (\ref{nullspace}) with the same $\kappa$ and this implies the $s$-goodness of $A$. \epr

\subsubsection{Proof of Lemma \ref{Dantzig}}\label{proofDantzig}
Let $A$ satisfy $\RIP(\delta,2s)$ with $\delta<1$; we want to prove that then the matrix ${1\over 1-\delta}A$ satisfies the condition $\bH_{s,2}({\delta\over 1-\delta})$. Indeed, let $x\in\bR^n$. Let vectors $x^1,\,x^2,\,...,\,x^q$ be obtained from $x$ as follows: $x^1$ is obtained by zeroing all but the $s$ largest in magnitude entries of $x$ and keeping the latter entries intact, then $x^2$ is obtained by applying  the same procedure to $x-x^1$, and so on. We stop at step $q$ where we get $x=x^1+...+x^q$.
Observe that $\|x^j\|_\infty\leq {1\over s}\|x^{j-1}\|_1$, whence also $\|x^j\|_2\leq \|x^{j-1}\|_1s^{-1/2}$ (since $x^j$ is $s$-sparse). We now have
\bse
\sqrt{s}\|x^1\|_2\|A^TAx\|_\infty&\geq& \|x^1\|_1\|A^TAx\|_\infty\geq [x^1]^TA^TAx=[x^1]^TA^TAx^1+\sum_{j=2}^q[x^1]^TA^TAx^j\\
&\geq& (1-\delta)\|x^1\|_2^2-\delta\sum_{j=2}^q\|x^1\|_2\|x^j\|_2\;\;\;\;~\;\;\;\;(*)\\
&\geq& (1-\delta)\|x^2\|_2^2-\delta s^{-1/2}\sum_{j=2}^q\|x^1\|_2\|x^{j-1}\|_1\geq (1-\delta)\|x^2\|_2^2-\delta s^{-1/2}\|x^1\|_2\|x\|_1\\
\Rightarrow\;\;\; \|x\|_{s,2}&=&\|x^1\|_2\leq {1\over 1-\delta}\|A^TAx\|_\infty+{\delta\over 1-\delta}\|x\|_1
\ese
(in the above chain, step $(*)$ is valid due to $[x^1]^TA^TAx^1\geq (1-\delta)\|x^1\|_2$ (since $A$ is $\RIP(\delta,2s)$) and the statement (!), see the proof of Proposition \ref{propRIP}).
The concluding inequality in the chain says that ${1\over1-\delta}A$ satisfies $\bH_{s,2}({\delta\over1-\delta})$. \epr

\subsubsection{Proof of Proposition \ref{arikwillnotlikeit}}\label{proofproparikwillnotlikeit}
We present here the proof of (i), which is a straightforward modification of the proof of Proposition \ref{prop101}. The proof of (ii) can be obtained by equally straightforward modification of the proof of Proposition \ref{prop11}.
\par
Thus, suppose we are under the premise of (i), and let $\Xi$ be defined exactly as in the proof of Proposition \ref{prop10}, so that $\Prob\{\xi\in\Xi\}\geq1-\epsilon$ and  $|(\sigma\xi+u)^Th_i|\leq\nu_{\epsilon,\sigma,\cU}(h_i)\leq\nu(H)$ for all $\xi\in \Xi$, $u\in\cU$ and all $i$. Let us fix $x\in\bR^n$, $\xi\in\Xi$ and $u\in\cU$, let $I$ be the set of indices of the $s$ largest in magnitude entries in $x$, and let $\eta=\sigma\xi+u$, $y=Ax+\eta$,  $\widehat{x}=\widehat{x}_{\reg}(y)$, and $z=\widehat{x}-x$.
\par
Since $\xi\in\Xi$ and $u\in\cU$, we have $|h_i^T(Ax-y)|\leq\nu_i\leq\rho_i$ for all $i$, whence $x$ is a feasible solution to the optimization problem defining $\widehat{x}$, whence, exactly as in the proof of Proposition \ref{prop10},
\be
(a)&\|z\|_1&\leq 2\|z_I\|_1+2\|x_{J}\|_1,\nn
(b)&|h_i^TAz|&
\le \rho_i+\nu_i,\,1\leq i\leq n\nn
\Rightarrow&\|H^TAz\|_\infty&\leq \widehat{\rho}+\widehat{\omega}.
\ee{neweq12}
Now, $H$ satisfies the condition $\bH_{s,q}(\kappa)$ and thus satisfies the condition $\bH_{s,1}(\kappa)$. Applying the latter condition, we get
$$
\|z_I\|_1\leq s\|H^TAz\|_\infty+\kappa\|z\|_1.
$$
Invoking (\ref{neweq12}), we conclude that
\begin{equation}\label{ell1norm}
\|z\|_1\leq 2s\|H^TAz\|_\infty+2\kappa\|z\|_1+2\|x_J\|_1\leq 2s[\widehat{\rho}+\widehat{\omega}+s^{-1}\|x_J\|_1]+2\kappa \|z\|_1,
\end{equation}
thus
\be
\|z\|_1\le {2s\over 1-2\kappa}[\widehat{\rho}+\widehat{\omega}+s^{-1}\|x_J\|_1].
\ee{neweqa}
Next, $H$ satisfies $\bH_{s,q}(\kappa)$, whence $\|z\|_{s,q}\leq \|H^T Az\|_\infty +\kappa\|z\|_1$. Therefore, we get from \rf{neweq12}:
\be
\|z\|_{s,q}&\leq &s^{{1\over q}}\|H^TAz\|_\infty+\kappa s^{{1\over q}-1}\|z\|_1
\leq s^{{1\over q}}[\widehat{\rho}+\widehat{\omega}]+
2\kappa s^{{1\over q}}{\widehat{\rho}+\widehat{\omega}+s^{-1}\|x_J\|_1 \over 1-2\kappa}\\
&\le &s^{1\over q}{\widehat{\rho}+\widehat{\omega}+2\kappa s^{-1}\|x_J\|_1\over 1-2\kappa}\le
s^{1\over q}{\widehat{\rho}+\widehat{\omega}+s^{-1}\|x_J\|_1\over 1-2\kappa}.
\ee{neweq13}
All we need in order to extract (i) from (\ref{ell1norm}) and (\ref{neweq13}) is to verify that
$$
1\leq p\leq q\;\Rightarrow\; \|z\|_p\leq (3s)^{1/p}\theta, \;\; \theta={\widehat{\rho}+\widehat{\omega}+s^{-1}\|x_J\|_1\over 1-2\kappa}.
$$
The desired inequality holds true when $p=1$ (see (\ref{ell1norm})), thus, invoking the H\"older inequality, all we need is to verify that
\be
\|z\|_q\leq (3s)^{1/q}\theta.
\ee{(!)}
When $q=\infty$, \rf{(!)} is implied by (\ref{neweq13}), so let us assume that $q<\infty$. Let $\lambda$ be the $(s+1)$-st largest of the magnitudes of entries in $z$. By (\ref{neweq13}) we have $\lambda^qs\leq \|z\|_{s,q}^q\leq s\theta^q$, and $\lambda\leq \theta$. Hence, setting $z'=z-z^s$, we get
\[
\|z'\|_q^q\leq\lambda^{q-1}\|z'\|_1\leq \theta^{q-1}\|z\|_1\leq \theta^{q-1}2s\theta,
 \]
 where the concluding inequality is given by (\ref{ell1norm}). Thus, $\|z'\|^q_q\leq 2s\theta^q$, while $\|z^s\|^q_q\leq s\theta^q$ by (\ref{neweq13}). We see that $\|z\|^q_q\leq \|z^s\|^q_q+\|z'\|^q_q\leq 3s\theta^q$, as required in \rf{(!)}.
\epr
\subsubsection{Proof of Proposition \protect{\ref{LassoDantzig}}}
\paragraph{(i):} Let  $\Xi=\{\xi\in\bR^m:|\xi^Ta_i|\leq \varrho,1\leq i\leq n\}$, so that $\Prob\{\xi\in \Xi\}\geq 1-\epsilon$. Let us fix $\xi\in\Xi$ and $x\in \bR^n$, and let $y=Ax+\sigma\xi$, $\widehat{x}=\widehat{x}_{\ds}(y)$. We have $\|A^T(Ax-y)\|_\infty=\|A^T\sigma\xi\|_\infty\leq\varrho\leq\rho$, so that $x$ is a feasible solution to the optimization problem
specifying $\widehat{x}_{\ds}(y)$ and therefore $\|\widehat{x}\|_1\leq \|x\|_1$. Denoting by $I$ the support of $x^s$, setting $z=\widehat{x}-z$ and acting exactly as when deriving (\ref{eqeqeq}), we arrive at
\begin{equation}\label{eqeqeqa}
 \|z\|_1\leq 2\|z_I\|_1+2\|x_J\|_1.
\end{equation}
Further,
\[
\|A^TAz\|_\infty=\|A^T(A\widehat{x}-y+\sigma\xi)\|_\infty\leq \|A^T\sigma\xi\|_\infty+\|A^T(A\widehat{x}-y)\|_\infty\leq \rho+\varrho,
\]
and therefore
\begin{equation}\label{end16}
\|Az\|_2^2=z^TA^TAz\le \|z\|_1\|A^TAz\|_\infty\le (\rho+\varrho)\|z\|_1.
\end{equation}
On the other hand, by (\ref{then}) we have
\[
\|z_I\|_1\le\|z\|_{s,1}\leq s^{1-{1\over q}}\|z\|_{s,q}\leq s\widehat{\lambda}\|Az\|_2+\kappa\|z\|_1\le s\widehat{\lambda}(\rho+\varrho)^{1/2}\|z\|_1^{1/2}+\kappa\|z\|_1.
\]
Substituting the above bound into (\ref{eqeqeqa}), we get
\bse
\|z\|_1\le 2\kappa\|z\|_1+2s\widehat{\lambda}(\rho+\varrho)^{1/2}\|z\|^{1/2}_1+2\|x_J\|_1,
\ese
whence by elementary calculations
\begin{equation}\label{elementary}
\tau:=\|z\|_1^{1/2}\leq a+b^{1/2},\;\;\mbox{where}\;\;a={2s\widehat{\lambda}\sqrt{\rho+\varrho}\over 1-2\kappa},\,\,b={2\|x_J\|_1\over 1-2\kappa}.
\end{equation}
Invoking (\ref{then}) and (\ref{end16}), we have
\be
\|z\|_{s,q}&\leq& s^{{1\over q}}\widehat{\lambda}\|Az\|_2+\kappa s^{{1\over q}-1}\|z\|_1\leq s^{{1\over q}-1}\left[s\widehat{\lambda}\sqrt{\rho+\varrho}\|z\|_1^{1/2}+\kappa \|z\|_1\right]=
s^{{1\over q}-1}\left[{(1-2\kappa)\over 2}a\tau+\kappa \tau^2\right]\nn
&\leq& s^{{1\over q}-1}\left[{1-2\kappa\over 2}{a^2}+{1-2\kappa\over 2}ab^{1/2}+\kappa [a+b^{1/2}]^2\right]\le s^{{1\over q}-1}\left[{a^2\over 2}+{1+2\kappa\over 2}ab^{1/2}+\kappa b\right]\nn
&\le & {s^{{1\over q}-1}\over 2}[a+b^{1/2}]^2,
\ee{wehave33}
where the last inequality of the chain is due to $\kappa<\half$.
Assuming for a moment that $1<q<\infty$ and denoting by $\mu$ the $(s+1)$-st largest magnitude of entries in $z$, we conclude from the latter inequality that $\mu\leq {1\over 2s}[a+b^{1/2}]^2$. Hence, when setting $z'=z-z^s$ we obtain (cf. the verification of \rf{(!)})
$\|z'\|^q_q\leq \mu^{{q-1}}\|z'\|_1\leq (2s)^{1-q}[a+b^{1/2}]^{2q}$. Invoking (\ref{wehave33}) one more time we get
$$
\|z\|_q^q\le\|z^s\|_q^q+\|z'\|^q_q\le {3\over 2} (2s)^{1-q} [a+b^{1/2}]^{2q}.
$$
The resulting inequality combines with (\ref{elementary}) and the H\"older inequality to imply that
\be
\|z\|_p\leq {3^{1\over p} s^{{1\over p}-1}\over 2} [a+b^{1/2}]^{2}\le 3^{1\over p} s^{{1\over p}-1}[a^2+b],\;\;1\leq p\leq q.
\ee{implythat}

Note that the derivation of (\ref{implythat}) was carried out under the additional assumption that $1<q<\infty$. This assumption can now be removed: when $q=1$, (\ref{implythat}) is readily given by (\ref{elementary}). When $q=\infty$, $A$ satisfies (\ref{then}) for $q=\infty$ and thus -- for every  value of $q$ from $[1,\infty]$, meaning that (\ref{implythat}) holds true for every $q<\infty$, whence (\ref{implythat}) holds true for $q=\infty$ as well.
\par
Recalling that relation (\ref{implythat}) is valid whenever $\xi\in\Xi$ and $x\in \bR^n$ and plugging in the values of $a$ and $b$,  we arrive at (\ref{DantzigBound}). (i) is proved.
\paragraph{(ii):}  Same as above, let $\Xi=\{\xi:|a_i^T\xi|\leq \sqrt{2\ln(n/\epsilon)} \|a_i\|_2,\,1\leq i\leq n\}$, so that $\Prob\{\xi\not\in\Xi\}\leq\epsilon$. Let us fix $x\in\bR^n$, $\xi\in\Xi$,  and let $I$ be the support of $x^s$.  Let also $y=Ax+\sigma\xi$, $\widehat{x}=\widehat{x}_{\las}(y)$, $z=\widehat{x}-x$. We have
$$
\|\widehat{x}\|_1+\varkappa\|A\widehat{x}-y\|_2^2\leq \|x\|_1+\varkappa\|Ax-y\|_2^2\leq \|x\|_1+\sigma^2\varkappa\xi^T\xi,
$$
or, which is the same due to $A\widehat{x}-y=Az-\sigma\xi$,
$$
\|\widehat{x}\|_1+\varkappa\|Az\|_2^2-2\sigma\varkappa\xi^TAz\leq \|x\|_1.
$$
It follows that
\bse
0&\geq& \|\widehat{x}\|_1-\|x\|_1+\varkappa(\|Az\|_2^2-2\sigma\xi^TAz)\\
&=&(\|x_I+z_I\|_1-\|x_I\|_1)+(\|x_J+z_J\|_1-\|x_J\|_1)+\varkappa(\|Az\|_2^2-2\sigma\xi^TAz)\\
&\geq&-\|z_I\|_1+\|z_J\|_1-2\|x_J\|_1+\varkappa(\|Az\|_2^2-2\sigma\xi^TAz)\\
&\geq&   -\|z_I\|_1+\|z_J\|_1-2\|x_J\|_1-2\varkappa\varrho\|z\|_1+\varkappa\|Az\|_2^2,\\
\ese
where the last $\geq$ is readily given by the fact that $\|A^T\sigma\xi\|_\infty\leq \varrho$ for $\xi\in\Xi$.  We conclude that
$$
\|z_J\|_1\leq \|z_I\|_1+2\varkappa\varrho\|z\|_1-\varkappa\|Az\|_2^2+2\|x_J\|_1,
$$
and therefore
\begin{equation}\label{lasso12}
\|z\|_1\leq 2\|z_I\|_1+2\varkappa\varrho\|z\|_1-\varkappa\|Az\|_2^2+2\|x_J\|_1.
\end{equation}
Now, we have
$$
\|z_I\|_1\leq s^{1-{1\over q}}\|z\|_{s,q}\leq s\widehat{\lambda}\|Az\|_2+\kappa\|z\|_1,
$$
where the concluding inequality is given by (\ref{then}). Combining the resulting inequality with (\ref{lasso12}),
we get
$$
\|z\|_1(1-2\kappa-2\varkappa\varrho)\leq 2s\widehat{\lambda}\|Az\|_2-\varkappa\|Az\|_2^2+2\|x_J\|_1.
$$
Combining this inequality with (\ref{lasso12}), we get the first inequality in the following chain:
\begin{eqnarray}\label{lassol2b}
\|z\|_1&\leq& 2(\kappa+\varkappa\varrho)\|z\|_1+2\|x_J\|_1+\varkappa\left(2\|Az\|_2{s\widehat{\lambda} \over\varkappa}-\|Az\|_2^2\right)\\
&\le&2(\kappa+\varkappa\varrho)\|z\|_1+2\|x_J\|_1+{s^2\widehat{\lambda}^2\over\varkappa}\nonumber
\end{eqnarray}
and since $2\kappa+2\varkappa\varrho<1$, we arrive at
\begin{equation}\label{arriveq12}
\|z\|_1\leq a:={1\over1-2(\kappa+\varkappa\varrho)}\left[{s^2\widehat{\lambda}^2\over\varkappa}+2\|x_J\|_1\right].
\end{equation}
Since $2(\kappa+\varrho\varkappa)< 1$, the first inequality in \rf{lassol2b} is  possible only if $$\|Az\|_2^2-{2s\widehat{\lambda}\over\varkappa}\|Az\|_2-{2\|x_J\|_1\over\varkappa}\le 0,$$ whence
\begin{equation}\label{Abound}
\|Az\|_2\le {2s\widehat{\lambda}\over\varkappa}+{\|x_J\|_1\over s\widehat{\lambda}}.
\end{equation}
Invoking (\ref{then}), we get $\|z\|_{s,q}\leq s^{{1\over q}-1}\left[s\widehat{\lambda}\|Az\|_2+\kappa\|z\|_1\right]$, which combines with (\ref{Abound}) and (\ref{arriveq12}) to imply that
\begin{equation}\label{get264}
\|z\|_{s,q}\leq 2s^{{1\over q}-1}a.
\end{equation}
Denoting by $\mu$ the $(s+1)$-st largest of the magnitudes of entries in $z$, we conclude from (\ref{get264}) that $\mu\leq 2s^{-1}a$, whence, denoting $z'=z-z^s$,
$$
\|z'\|_q\leq \mu^{{q-1\over q}}\|z'\|_1^{{1\over q}}\leq (2s^{-1}a)^{{q-1\over q}}\|z\|_1^{{1\over q}}\leq
2^{{q-1\over q}}s^{{1\over q}-1}a,
$$
(we have used (\ref{arriveq12})),
which combines with (\ref{get264}) to imply that
\begin{equation}\label{eq447}
\|z\|_q\leq 4s^{{1\over q}-1}a.
\end{equation}
Combining (\ref{eq447}), (\ref{arriveq12}) and the H\"older inequality, we get
\begin{equation}\label{LassoB}
1\leq p\leq q\Rightarrow \|z\|_p\leq 4s^{{1\over p}-1}a.
\end{equation}
Plugging in the value of $a$ (see (\ref{arriveq12})) and recalling that (\ref{LassoB}) takes place whenever $\xi\in\Xi$ with $\Prob\{\xi\in\Xi\}\geq1-\epsilon$, we arrive at (\ref{LassoBound}). \epr
\subsection{Proof of Proposition \ref{Prop_NEMP}}
\label{App:Prop_NEMP}
The proof below follows the lines of the proof of Proposition 7 of \cite{JKKN}.
 Given $\epsilon\in(0,1)$, let  $\Xi=\{\xi:|h_i^T\xi|\leq \sqrt{2\ln(n/\epsilon)}\|h_i\|_2,\;1\le i\leq n\},$ so that $\Prob\{\xi\in\Xi\}\geq1-\epsilon$. Let us fix $\xi\in\Xi$, $x\in \bR^n$ such that
 $\|x-x^s\|_1\leq\upsilon$, and $u\in\cU$. For $\eta=y-Ax=u+\sigma\xi$, by the definition \rf{nunorm} of the norm  $\nu$ and because of $\nu(h_i)\le \nu(H)$,  we have  $\|H^T\eta\|_\infty\le \nu(H)\leq\omega_*(\bar{\gamma})$.
\par
We intend to prove the relations $(a_k)$, $(b_k)$ by induction in $k$. First, let us show that $(a_{k-1},b_{k-1})$ implies $(a_{k},b_{k})$. Thus, assume that
$(a_{k-1},b_{k-1})$ holds true.
Let $z^{(k-1)}=x-v^{(k-1)}$. By $(a_{k-1})$, $z^{(k-1)}$ is supported on the support of $x$. Note that
\bse
z^{(k-1)}-u&=&x-v^{(k-1)}-H^T(y-Av^{(k-1)})=(I-H^TA)(x-v^{(k-1)})-H^T\eta\\
&=&(I-H^TA)z^{(k-1)}-H^T\eta,
\ese
 Then by \rf{feq99} for any $1\le i\le n$,
\bse
-\bar{\gamma}\sum_{j}|z^{(k-1)}_j|-\omega_*(\bar{\gamma})\le z^{(k-1)}_i-u_i\le
\bar{\gamma}\sum_{j}z^{(k-1)}_j +\omega_*(\bar{\gamma}),
\ese
consequently,
\be
-\gamma:=-\bar{\gamma} \alpha_{k-1}-\omega_*(\bar{\gamma})\leq z^{(k-1)}_i-u_i\leq \bar{\gamma}\alpha_{k-1}+\omega_*(\bar{\gamma}):=\gamma,
\ee{cover}
so that the segment $S_i=[u_i-\gamma,\,u_i+\gamma]$ of the width
$
\ell=2\bar{\gamma}\alpha_{k-1}+2\omega_*(\bar{\gamma}),
$
covers $z^{(k-1)}_i$, and the closest to zero point of this interval is
\bse
\widetilde{\Delta}_i=\left\{\begin{array}{ll}
[u_i-\gamma]_+,& u_i\ge 0,\\
-[|u_i|-\gamma]_+,& u_i< 0,
\end{array}\right.
\ese
that is, $\widetilde{\Delta}_i=\Delta_i$ for all $i$. Since the segment $S_i$ covers $z^{(k-1)}_i$ and $\Delta_i$ is the closest to 0 point in $S_i$, while the width of $S_i$ is at most $\ell$, we clearly have
\be\begin{array}{rlcrl}
(a)&\Delta_i\in \Conv\left\{0,z^{(k-1)}_i\right\},&
(b)&|z^{(k-1)}_i-\Delta_i|\le \ell.\\
\end{array}
\ee{almost}
Since $(a_{k-1})$ is valid, (\ref{almost}.a) implies that
\[
v^{(k)}_i=v^{(k-1)}_i+\Delta_i\in\left[ v^{(k-1)}_i+\Conv\left\{0,x_i-v^{(k-1)}_i\right\}\right]\subseteq\Conv\{0,x_i\},
\]
and $(a_k)$ holds.
Further, let $I$ be the support of $x^s$. Relation $(a_k)$ clearly implies that $|z^{(k)}_i|\le |x_i|$, and
 we can write due to (\ref{almost}.b):
$$
\|x-v^{(k)}\|_1=\sum_{i\in I} |x_i-[v^{(k-1)}_i+\Delta_i]|+\sum_{i\not\in I}|z^{(k)}_i|
\le \sum_{i\in I} |z^{(k-1)}_i-\Delta_i|+\sum_{i\not\in I}|x_i|\le s\ell+\upsilon=\alpha_k.
$$
Since by (\ref{almost}.b)
\bse
\|x-v^{(k)}\|_\infty&=&\|x-v^{(k-1)}-\Delta\|_\infty\le \ell=2\bar{\gamma}\alpha_{k-1}+2\omega_*(\bar{\gamma}),
\ese
we conclude that $(b_k)$ holds true. The induction step is justified.
\par
It remains to show that $(a_0,b_0)$ holds true. Since $(a_0)$ is evident, all we need is to justify $(b_0)$.  Let
\[
\alpha_*=\|x\|_1,
\]
and let $u=H^Ty$. Same as above (cf. \rf{cover}), we
have for all $i$:
$$
|x_i-u_i|\le \bar{\gamma}\alpha_*+\omega_*(\bar{\gamma}).
$$
Then
\bse
\alpha_*=\sum_{i\in I}|x_i|+\sum_{i\not \in I}|x_i|\leq \sum_{i\in I} [|u_i|+
\bar{\gamma}\alpha_*+\omega_*(\bar{\gamma})]+\upsilon\le \|u\|_{s,1}+
s\bar{\gamma}\alpha_*+s\omega_*(\bar{\gamma})+\upsilon.
\ese
Hence
\[
\alpha_*\leq \alpha_0={\|u\|_{s,1}+s\omega_*(\bar{\gamma})+\upsilon\over 1-s\bar{\gamma}},
\]
which implies $(b_0)$.   \epr

\end{document}